\documentclass[aps,pre,amsmath,amsfonts,amssymb,superscriptaddress,showpacs,12pt]{revtex4}
\usepackage{mathptmx,amsthm,subfigure}
\usepackage{graphicx}
\usepackage[latin1]{inputenc}
\usepackage{amscd}
\usepackage{soul}

\usepackage{hyperref}

\usepackage{color}

\usepackage{amsthm}
\usepackage{framed}
\usepackage{amssymb}
\usepackage{amsmath,mathrsfs}
\usepackage{hhline}
\usepackage{ulem}
\usepackage{fancyhdr}
\usepackage{tikz}

\everymath{\displaystyle}
\everymath{\displaystyle}

\everymath{\displaystyle}
\newcommand{\beq}{\begin{equation}}
\newcommand{\eeq}{\end{equation}}

\newcommand{\bi}{\begin{itemize}}
\newcommand{\ei}{\end{itemize}}

\newcommand{\U}{\mathrm{U}}

\def\RR{\mathbb{R}}

\newcommand{\Tau}{\mathcal{T}}
\newcommand{\Loop}{\mathcal{L}}
\newcommand{\hol}{\mathfrak{hol}}
\newcommand{\Ma}{\mathrm{M}}
\newcommand{\esse}{\mathcal{S}}

\newcommand{\Div}{\mathcal{D}}
\newcommand{\nihat}{\mathrm{X}}
\newcommand{\grad}{\mathrm{grad}}
\newcommand{\total}{\mathrm{d}}
\newcommand{\RC}{\mathrm{R}}

\newcommand{\Hol}{{\rm Hol}}
\newcommand{\Aut}{{\rm Aut}}
\newcommand{\End}{{\rm End}}

\newcommand{\Riemann}{\mathcal{R}}
\newcommand{\tangent}{{\rm T}}

\newcommand{\paralleltransport}{{\rm P}}
\newcommand{\Lag}{{\rm L}}
\newcommand{\hyper}{{\rm S}}
\newcommand{\hyperm}{{\rm H}}
\newcommand{\metric}{{\rm g}}
\newcommand{\btheta}{\boldsymbol{ \theta}}
\newcommand{\bxi}{\boldsymbol{ \xi}}
\newcommand{\bz}{\boldsymbol{ z}}

\newcommand{\Lc}{\overline{\nabla}_{\mbox{\small LC}}}
\newcommand{\bieta}{\boldsymbol{\eta}}

\newtheorem{definition}{Definition}[section]

\newtheorem{theorem}{Theorem}[section]    
\newtheorem{lemma}{Lemma}[section]   
\newtheorem{pro}{Proposition}[section]
    
\newtheorem{remark}{Remark}[section]

\newtheorem{corollary}{Corollary}

\begin{document}

\title{\textbf{Towards a Canonical Divergence within Information Geometry}}
\author{Domenico Felice}
\email{DFelice@sunypoly.edu}
\affiliation{SUNY Polytechinic Institute, Albany, New York, 12203, USA\\
Max Planck Institute for Mathematics in the Sciences\\
 Inselstrasse 22--04103 Leipzig,
 Germany}
\author{Nihat Ay}
\email{nihat.ay@tuhh.de}
\affiliation{Hamburg University of Technology, Hamburg, Germany\\
 Max Planck Institute for Mathematics in the Sciences\\
 Inselstrasse 22--04103 Leipzig,
 Germany\\
 Faculty of Mathematics and Computer Science, University of Leipzig, PF 100920, 04009, Leipzig, Germany\\
Santa Fe Institute, Santa Fe, NM 87501, USA}

\begin{abstract}
In Riemannian geometry geodesics are integral curves of the Riemannian distance gradient. We extend this classical result to the framework of Information Geometry. In particular, we prove that the rays of level-sets defined by a pseudo-distance are generated by the sum of two tangent vectors. By relying on these vectors, we propose a novel definition of {a} canonical divergence and its dual function. We prove that the new divergence allows to recover a given dual structure $(\metric,\nabla,\nabla^*)$ of {a dually convex set on} a smooth manifold $\Ma$. Additionally, we show that this divergence {coincides with} the canonical divergence proposed by Ay and Amari in the case of: (a) self-duality, (b) dual flatness, (c) statistical geometric analogue of the concept of symmetric spaces in Riemannian geometry. {For a dually convex set,} the case (c) leads to a further comparison of the {new} divergence {with} the one introduced by Henmi and Kobayashi.
\end{abstract}

\pacs{Classical differential geometry (02.40.Hw), Riemannian geometries (02.40.Ky), Inverse problems (02.30.Zz).}

\maketitle

\newpage

\tableofcontents

\section[Introduction: Inverse problem and divergence functions in Information Geometry]{Introduction: {Inverse problem and divergence functions in Information Geometry}}\label{Intro}
The Inverse Problem within Information Geometry \cite{Amari00} concerns the search for a divergence function $\Div$ which {recovers a given dual structure $(\metric,\nabla,\nabla^*)$ of a smooth manifold $\Ma$.}

A {{\it dual structure} (or {\it dualistic structure}) on $\Ma$ is specified in terms of a metric tensor $\metric$ and two linear connections, $\nabla$ and $\nabla^*$, on the tangent bundle $\tangent\Ma$ such that 
\begin{equation}
\label{dualconnection}
X\ \metric\left(Y,Z\right)=\metric\left(Y,\nabla^*_X Z\right)+\metric\left(\nabla_X Y,Z\right),\ \forall\ X,Y,Z\in\Tau(\Ma),
\end{equation}
where $\Tau(\Ma)$ denotes the space of vector fields on $\Ma$, namely $C^{\infty}$ sections $X:\Ma\rightarrow\tangent\Ma$, $X_p\in\tangent_p\Ma$. The quadruple $(\Ma,\metric,\nabla,\nabla^*)$ is called {\it statistical manifold} whenever the dual connections $\nabla$ and $\nabla^*$ are both torsion free \cite{Ay17}.} The notion of a statistical manifold, introduced by Lauritzen \cite{Lauritzen87}, is usually referred to the triple $(\Ma,\metric,T)$, where $T(X,Y,Z)=\metric\left(\nabla^*_X Y-\nabla_X Y,Z\right)$ is a $3$-symmetric tensor. However, when $\nabla$ and $\nabla^*$ are both torsion free connections, the structures $(\Ma,\metric,\nabla,\nabla^*)$ and $(\Ma,\metric,T)$ are equivalent \cite{Ay17}.

A distance-like function $\Div:\Ma\times\Ma\rightarrow\RR $ {satisfying  
\begin{equation}
\label{distance-like}
\Div(p,q)\geq 0 \,\,\, \forall p,q\in\Ma\quad \mbox{and} \quad \Div(p,q)=0\ \mbox{iff}\ p=q
\end{equation}
is called a {\it divergence} or {\it contrast function}} on $\Ma$ if the matrix
\begin{equation}
\label{metricfromdiv}
\metric_{ij}(p)=-\left.\partial_i\partial_j^{\prime} \Div(\bxi_p,\bxi_q)\right|_{p=q}=\left.\partial^{\prime}_i\partial_j^{\prime} \Div(\bxi_p,\bxi_q)\right|_{p=q}
\end{equation}
is strictly positive definite everywhere on $\Ma$ \cite{eguchi1983}. Here, 
\begin{equation*}
\partial_i=\frac{\partial}{\partial \xi_p^i} \quad \mbox{and}\quad \partial^{\prime}_i=\frac{\partial}{\partial \xi_q^i}
\end{equation*} 
and $\{\bxi_p:=(\xi_p^1,\ldots,\xi_p^n)\}$ and $\{\bxi_q:=(\xi_q^1,\ldots,\xi_q^n)\}$ are local coordinate systems of $p$ and $q$, respectively. 
{Given a dual structure $(\metric,\nabla,\nabla^*)$ on $\Ma$, the divergence fucntion} \eqref{distance-like} is {said to be} compatible with $(\metric,\nabla,\nabla^*)$ if $\metric$ is obtained by \eqref{metricfromdiv} and furthermore the following holds \cite{eguchi1992}:
\begin{align}\label{dualfromDiv}
& \Gamma_{ijk}(p)=-\left.\partial_i\partial_j\partial_k^{\prime} \Div(\bxi_p,\bxi_q)\right|_{p=q}, \qquad {\Gamma}^*_{ijk}(p)=-\left.\partial^{\prime}_i\partial^{\prime}_j\partial_k \Div(\bxi_p,\bxi_q)\right|_{p=q}\ ,
\end{align} 
where $\Gamma_{ijk}=\metric\left(\nabla_{\partial_i}\partial_j,\partial_k\right)$, ${\Gamma}^*_{ijk}=\metric\left(\nabla^*_{\partial_i}\partial_j,\partial_k\right)$ are the symbols of the dual connections $\nabla$ and $\nabla^*$, respectively. In this article, we address our investigation to {the attempt of finding out} a divergence {which is {\it canonical} in a suitable way and} recovers the dual structure of a given statistical manifold $(\Ma,\metric,\nabla,\nabla^*)$. It is worth {noting} from Eqs. \eqref{metricfromdiv} and \eqref{dualfromDiv} that it is sufficient to define a divergence on a neighborhood of the diagonal in $\Ma\times\Ma$.

Matumoto \cite{matumoto1993} showed that, {given a torsion free dual structure $(\metric,\nabla,\nabla^*)$ on $\Ma$, there always exists a divergence function on $\Ma$ which induces that structure.} However, this is not unique and there are infinitely many divergences that give the same dual structure. {Hence, the search for a divergence which can be  considered as {\it the most natural}, { in some sense}, is of uppermost importance. To this end, Amari and Nagaoka introduced a Bregman type divergence on dually flat manifolds in terms of the Legendre transform between dual affine coordinates \cite{Amari00}. More precisely, since the statistical manifold $(\Ma,\metric,\nabla,\nabla^*)$ is dually flat, the curvature tensors $\RC(\nabla)$ and $\RC^*(\nabla^*)$ are zero (for more details see Appendix \ref{StatManifold}). Therefore, we are given on $\Ma$ mutually dual affine coordinates $(\{\theta^i\},\{\eta_i\})$ and their potentials $(\phi,\phi^*)$, namely two smooth functions on $\Ma$ such that $\frac{\partial \phi}{\partial \theta^i}=\eta_i$ and $\frac{\partial \phi^*}{\partial \eta_i}=\theta^i$. Given $p,q\in\Ma$, the Bregman type divergence is then defined by
\begin{equation}
\label{BregmanDiv}
D[p:q]:=\phi(p)+\phi^*(q)-\sum_i\theta^i(p)\,\eta_i(q)\,.
\end{equation} 
This divergence has relevant properties concerning the generalized Pythagorean theorem and the geodesic projection theorem \cite{Amari16} and it is referred to as {\it canonical divergence} and  commonly assessed as the natural solution of the inverse problem in Information Geometry for dually flat manifolds. Extensions of the canonical divergence within conformal geometry have been analysed by Kurose \cite{kurose1994} and Matsuzoe \cite{matsuzoe1998}. However, the need for a general canonical divergence, which applies to any dualistic structure, is a very crucial issue as pointed out in \cite{AyTusch}. In any case, such a divergence should recover the canonical divergence of Bregman type if applied to a dually flat structure.
In addition, in the self-dual case where $\nabla=\nabla^*$ coincides with the Levi-Civita connection of $\metric$, the supposed canonical divergence should be one half of the squared Riemannian distance \cite{Ay17}. { However, these natural properties that a canonical divergence is required to satisfy are not sufficient for a unique characterisation.} An attempt { towards the definition of} a general canonical divergence on a general statistical manifold has been put forth in \cite{Ay15}. Here, the authors have introduced a canonical divergence that satisfies all these requirements. Such a divergence is defined in terms of geodesic integration of the inverse exponential map. This one is interpreted as a {\it difference} vector that translates $q$ to $p$ for all $q,p$ suitably close in $\Ma$.}

\begin{figure}
\begin{center}
\includegraphics[scale=1]{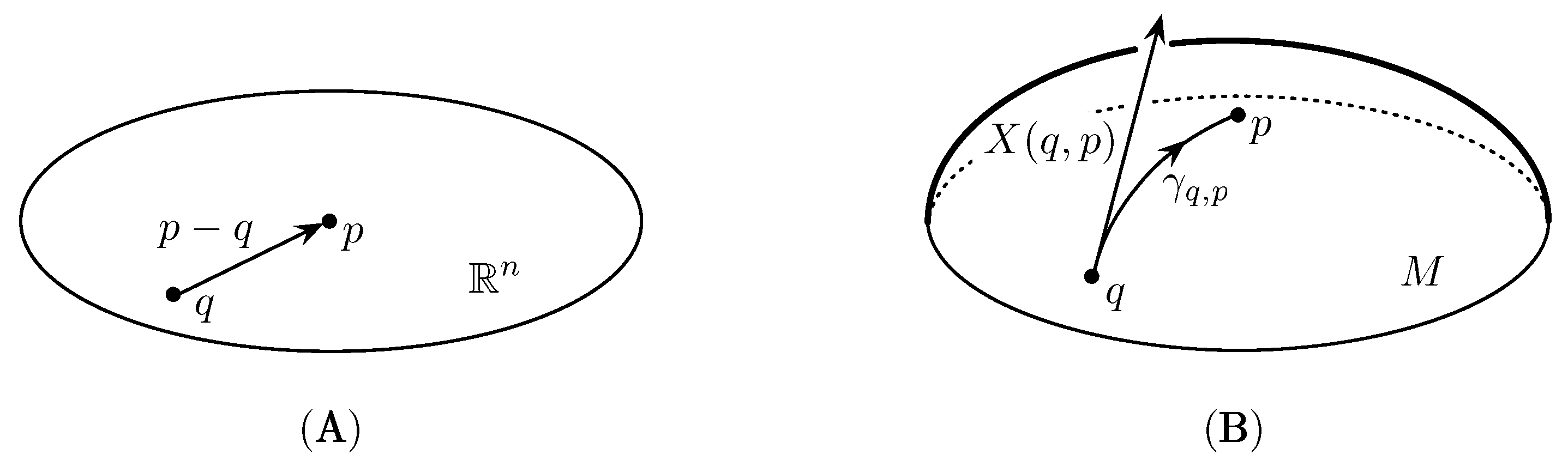}
\caption{On the left, (A) illustrates the difference vector $p-q$ in the linear vector space $\RR^n$; whereas, in (B) we can see the difference vector $\nihat(q,p)=\dot{\gamma}_{q,p}(0)$ in $\Ma$ as the inverse of the exponential map at $q$ (this Figure comes from \cite{Ay15}).}\label{differencevect}
\end{center}
\end{figure}

To be more precise, the inverse exponential map provides a generalization to $\Ma$ of the notion of difference vector in the linear vector space. In detail, let $p,q\in\RR^n$, the difference between $p$ and $q$ is given by the vector $p-q$ pointing to $p$ (see side (A) of Fig. \ref{differencevect}). Then, the difference between $p$ and $q$ in $\Ma$ is supplied by the exponential map of the connection $\nabla$ {(see Appendix \ref{StatManifold} for a more detailed discussion on linear connections and corresponding exponential maps)}. In particular, assuming that $p\in \mathrm{U}_q$ and $\mathrm{U}_q\subset \Ma$ is a $\nabla$-geodesic neighborhood of $q$, the difference vector from $q$ to $p$ is defined as (see (B) of Fig. \ref{differencevect})
\begin{equation}
\label{nihat}
\nihat_q(p):=\nihat(q,p):=\exp_q^{-1}(p)=\dot{\gamma}_{q,p}(0)\ ,
\end{equation}
where $\gamma_{q,p}$ is the $\nabla$-geodesic from $q$ to $p$ laying in $\mathrm{U}_q$. Clearly, by fixing $p\in\Ma$ and letting $q$ vary in $\Ma$, we obtain a vector field $\nihat(\cdot,p)$ whenever a $\nabla$-geodesic from $q$ to $p$ exists. From here on, we equally use both the notations, $\nihat(q,p)$ and $\nihat_q(p)$, for representing the difference vector from $q$ to $p$. Obviously, $\nihat_p(q):=\exp_p^{-1}(q)=:\nihat(p,q)$ denotes the difference vector from $p$ to $q$.

Therefore, the divergence proposed by Ay and Amari in \cite{Ay15} is defined as the path integral
\begin{equation}
\label{canonicaldiv}
D(p,q):=\int_0^1 \langle \nihat_t(p),\dot{\gamma}_{q,p}(t)\rangle_{\gamma_{q,p}(t)}\ dt\ ,
\end{equation}
where $\gamma_{q,p}$ is the $\nabla$-geodesic from $q$ to $p$ and $\langle\cdot,\cdot\rangle_{\gamma_{q,p}(t)}$ denotes the inner product with respect to $\metric$ evaluated at $\gamma_{q,p}(t)$. In Eq. \eqref{canonicaldiv}, $\nihat_t(p)$ is the vector field along $\gamma_{q,p}(t)$ given by Eq. \eqref{nihat} as follows,
\begin{equation}
\label{nihatvstime}
\nihat_t(p)=\nihat(\gamma_{q,p}(t),p)=\exp_{\gamma_{q,p}(t)}^{-1}(p)\ . 
\end{equation}
After elementary computations, Eq. \eqref{canonicaldiv} reduces to,
\begin{equation}
\label{Aydivergence}
D(p,q)=\int_0^1 \, t\, \|\dot{\gamma}_{p,q}(t)\|^2\ dt\ ,
\end{equation}
where $\gamma_{p,q}(t)$ is the $\nabla$-geodesic from $p$ to $q$ \cite{Ay15}.
If we consider {the} definition \eqref{canonicaldiv} for {a} general path $\gamma$ and denoting this as $D_{\gamma}(p,q)$, we will have dependence on $\gamma$. On the contrary, if the vector field $\nihat_t(p)$ is integrable, then $D_{\gamma}(p,q)=: D(p,q)$ turns out to be independent of the path from $q$ to $p$. A relevant issue closely related to the integrability of the vector field $\nihat_t(p)$ regards the geodesic projection property. Given a submanifold $\widetilde{\Ma}\subset\Ma$, we say that the divergence $D$ holds the geodesic projection property if { the following holds: for any ${p}\in\Ma\,\backslash\,\widetilde{\Ma}$ and any $\hat{p}\in \widetilde{\Ma}$ that locally minimizes $D(p,q)$ with respect to $q$, the geodesic from $p$ to $\hat{p}$ intersects $\widetilde{\Ma}$ orthogonally.} It turns out that a divergence holds the geodesic projection if and only if $\nihat(q,p)=c(p,q)\, \grad_q\,D(p,\cdot)$. Clearly, if $\nihat(q,p)$ is integrable, the latter relation is satisfied. In \cite{Ay15}, the authors have proved the integrability of $\nihat(q,p)$ in several contexts. However, the general problem is still open.

The dual divergence of $D$ has been defined in terms of the inverse exponential map with respect to the $\nabla^*$-connection. It turns out to be closely related to the divergence of the article \cite{Kobayashi00}. Here the authors applied the Hook's law to a ``$\nabla^*$-spring'' and defined the divergence as the physical work that is necessary to move a unit mass from $q$ to $p$ along the $\nabla$-geodesic $\gamma_{q,p}$ connecting them against the force field described by the inverse exponential map of $\nabla^*$. To be more precise, the force field along the $\nabla$-geodesic $\gamma_{q,p}(t)$ connecting $q$ with $p$ is defined in terms of $\nabla^*$-geodesics connecting $\gamma_{q,p}(t)$ with $q$:
\begin{equation}
\label{ForcefieldHenmi}
\mathrm{F}(\gamma_{q,p}(t)):=\left.\frac{\total}{\total s}\right|_{s=0}\sigma_t^*(s)=\stackrel{*}{\exp}^{-1}_{\gamma_{q,p}(t)}(q)=:\nihat^*_t(q)\,,
\end{equation}
where $\sigma_t^*(s)\,(0\leq s\leq 1)$ is the $\nabla^*$-geodesic such that $\sigma_t^*(0)=\gamma_{q,p}(t)$ and $\sigma^*_t(1)=q$, whereas $\stackrel{*}{\exp}$ denotes the exponential map of the $\nabla^*$-connection. Hence, the divergence of Henmi and Kobayashi is the work necessary to move a point of unit mass from $q$ to $p$ against the force field $\nihat_t^*(q)$ \cite{Kobayashi00}:
\begin{equation}
\label{HenmiDivergence}
W(p\| q):=-\int_0^1\left\langle\dot{\gamma}_{q,p}(t),\nihat_t^*(q)\right\rangle_{\gamma_{q,p}(t)}\,\total t\,.
\end{equation}
Under suitable conditions on the curvature tensors, $W(p\|q)$ turns out to be a potential function in the sense that it is independent of the particular path from $q$ to $p$ \cite{Kobayashi00}.

{ In this article, we further explore the idea of a canonical divergence as a path integral based on a vector field. This will be motivated and guided by the classical Gauss Lemma, which we first generalise to the context of a statistical manifold. Our analysis will naturally lead us to a canonical divergence which satisfies the outlined conditions and coincides with the previously proposed canonical divergences in some special but important cases. However, we should highlight that our approach, while providing a new perspective within the search for a general canonical divergence, does not conclude with a final answer. We propose a new candidate of a canonical divergence, which not only emphasises new geometric aspects but also sheds light on the previously proposed canonical divergences and thereby refines the search for the most natural one.}

\vspace{.2cm}

The layout of this article is as follows. In Section \ref{Overwiev} we outline the contents of the paper by stating the main results.  All the theory is then developed in the part entitled ``{\bf Technical details and proofs}" { which includes Sections \ref{GeodGeometry} to \ref{sec4}. In this part,} all the results discussed in Section \ref{Overwiev} are proved. In Section \ref{Conclusions} we draw some conclusions by outlining the results obtained in this work and discussing possible extensions. Useful tools related to statistical manifolds appear in Appendix \ref{StatManifold}.

\section{Contents and line of reasoning}\label{Overwiev}
In this manuscript we propose a new definition of a canonical divergence  on a general statistical manifold $(\Ma,\metric,\nabla,\nabla^*)$ through an extensive investigation of the geometry of geodesics. In this way and by combining  {the approaches carried out in \cite{Ay15} and in \cite{Kobayashi00} we succeed to supply an intrisic meaning to the new divergence. { Note} that the concept of a statistical manifold in Information Geometry generalizes the one of a Riemannian manifold. Indeed, when $\nabla=\nabla^*$, the statistical manifold $(\Ma,\metric,\nabla,\nabla^*)$ turns out to be a Riemannian manifold $(\Ma,\metric,\Lc)$ endowed with the Levi-Civita connection \cite{Amari16}. In this case, the geodesic structure of $(\Ma,\metric,\Lc)$ is fully understood thanks to the celebrated Gauss Lemma. This asserts that ``{\it the geodesic rays starting from $p\in\Ma$ are all orthogonal to the geodesic spheres centered at $p$}" \cite{Michor}.
\begin{theorem}[Gauss Lemma]\label{GaussLemma}
Let $(\Ma,\metric,\overline{\nabla}_{\mbox{\rm LC}})$ be a Riemannian manifold endowed with the Levi-Civita connection. For $p\in\Ma$ let $\varepsilon>0$ be so small that $\overline{\exp}_p:\overline{\mathcal{E}}_p(\varepsilon)\rightarrow\Ma$ is  a diffeomorphism on its image, where
$$
\overline{\mathcal{E}}_p(\varepsilon):=\left\{X\in\tangent_p\Ma\, |\, \sqrt{\langle X,X\rangle_p}<\varepsilon\right\}\,.
$$
Then in $\overline{\exp}_p(\overline{\mathcal{E}}_p(\varepsilon))$ the geodesic rays starting
from $p$ are all orthogonal to the ``geodesic spheres" 
$$
\overline{\mathrm{S}}_{\kappa}(p)=\left\{q\in \overline{\exp}_p\left(\overline{\mathcal{E}}_p(\varepsilon)\right)\ | \ \sqrt{\langle\overline{\exp}_p^{-1}(q),\overline{\exp}_p^{-1}(q)\rangle_p} = \kappa \right\}\,,
$$ 
for any $\kappa\in(0,\varepsilon)$\,.
\end{theorem} 
Here, $\overline{\exp}$ denotes the Levi-Civita exponential map (see Appendix\ref{StatManifold}). A relevant consequence of the Gauss Lemma is that the $\Lc$-geodesics generate the integral curves of the squared Riemannian distance gradient \cite{Lee97}.
\begin{theorem}
\label{thhgraddistance}
Let $\overline{\mathrm{S}}_{\kappa}(p)$ be a geodesic sphere centered at $p$ in a Riemannian manifold $(\Ma,\metric,\overline{\nabla}_{\mbox{\rm LC}})$. Consider the Riemannian distance $d_p:\overline{\mathrm{S}}_{\kappa}(p)\rightarrow\RR^+$ given by
\begin{equation}
\label{Riemanniandistance}
d(p,q):=\sqrt{\left\langle\overline{\exp}_p^{-1}(q),\overline{\exp}_p^{-1}(q)\right\rangle_p}\,,\quad \forall\,q\in \overline{\mathrm{S}}_{\kappa}(p)\,.
\end{equation}
Then, if $\overline{\sigma}(t)$ is the $\overline{\nabla}_{\mbox{\rm LC}}$-geodesic connecting $p$ and $q$,  we have that
\begin{equation}
\label{Energygradient}
\grad_q\, d(p,q)^2=2\,\dot{\overline{\sigma}}(1)\,.
\end{equation}
\end{theorem}  
In this article, we mimic the theory built around the Gauss Lemma and provide an extensive investigation of the geodesic geometry of a general statistical manifold $(\Ma,\metric,\nabla,\nabla^*)$. To this aim, we restrict our attention to a set $\U\subset\Ma$ such that for every $p,q\in\U$ there exist a unique $\nabla$-geodesic and a unique $\nabla^*$-geodesic connecting $p$ and $q$.
\begin{definition}
\label{dualconvexset}
Let $(\Ma,\metric,\nabla,\nabla^*)$ be a statistical manifold. A subset $\U\subset\Ma$ is said to be a {\bf dually convex set} if, given any pair of points $p,q\in\U$, we can find a unique $\nabla$-geodesic $\sigma:[0,1]\rightarrow\U$ and a unique $\nabla^*$-geodesic $\sigma^*:[0,1]\rightarrow\U$ such that $\sigma(0)=p=\sigma^*(0)$ and $\sigma(1)=q=\sigma^*(1)$.
\end{definition}

Notice that, given $\U\subset\Ma$ a dually convex set and $p\in\U$, we can find a neighborhood $\mathcal{E}_p$ of $0$ in $\tangent_p\Ma$ such that the $\nabla$ and $\nabla^*$ exponential maps $\exp_p:\mathcal{E}_p\rightarrow \exp_p(\mathcal{E}_p)$ and $\stackrel{*}{\exp}_p:\mathcal{E}_p\rightarrow \, \stackrel{*}{\exp}_p(\mathcal{E}_p)$ are diffeomorphisms.

Thus, given $\U$ a dually convex set, we introduce two vector fields, $\Pi$ and $\Pi^*$ herein, which generalize the concept of the geodesic tangent vector field as related to the gradient of the squared Riemannian distance according to Eq. \eqref{Energygradient}. Let $p,q\in\U$, we start by considering the difference vector from $p$ to $q$ with respect to the $\nabla$-connection:  
$$
\nihat_p(q):=\exp_p^{-1}(q)= \dot{\sigma}(0)\in\tangent_p\Ma\,.
$$
We then $\nabla$-parallel translate it along the $\nabla^*$-geodesic $\sigma^*$ from $p$ to $q$ (see Fig. \ref{perpendicular} side (\textbf{A})), and obtain} 
\begin{equation}
\label{P}
{\Pi_q(p):=}\paralleltransport_{\sigma^*}\nihat_p(q)\in\tangent_q\Ma\ .
\end{equation} 
At this point, by fixing $p$ and letting $q$ be varied, we obtain a vector field ${\Pi_q(p)}$ which corresponds to minus the difference vector defined in Eq. \eqref{nihat}. Analogously, we define the dual vector of $\Pi_q(p)$ as {the $\nabla^*$-parallel transport of $\dot{\sigma}^*(0)$ along the $\nabla$-geodesic $\sigma$ connecting $p$ with $q$ (see Fig. \ref{perpendicular} side (\textbf{B})),
\begin{equation}
\label{P*}
\Pi_q^*(p):=\paralleltransport^*_{\sigma}\nihat_p^*(q)\in\tangent_q\Ma,
\end{equation}
where
\begin{equation}
\label{nihat*}
\nihat_p^*(q):=\stackrel{*}{\exp}_p^{-1}(q) = \dot{\sigma}^*(0)\ .
\end{equation}

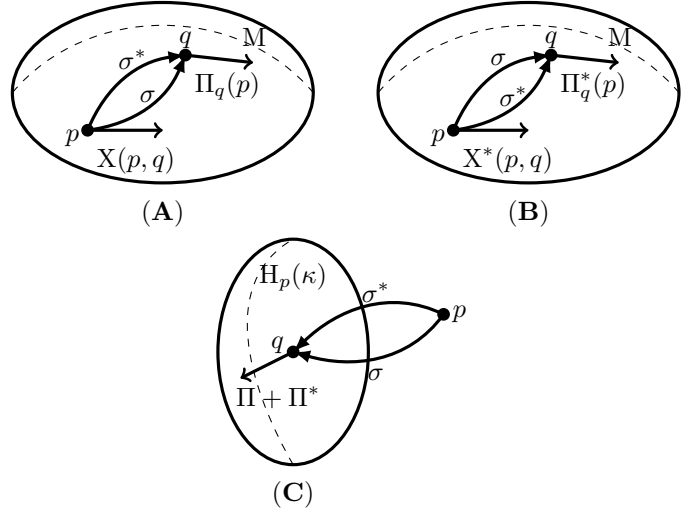
\begin{figure}[h!]
\centering
\begin{tikzpicture}
\draw [very thick] (0,0) ellipse (2 and 1.2);
\coordinate [label=above left:$\Ma$] (t) at (1.5,.5);
\draw [->, very thick] (-1,-.5) -- (0,-.5);
\draw [->, very thick] (0.3,.5) -- (1.2,.4);
\node (p) at (-1.2,-.6) {$p$};
\node (q) at (0.3,.7) {$q$};
\coordinate [label=below right:{\small $\nihat(p,q)$}] (X) at (-1,-.6);
\draw [-latex, bend left, very thick] (-1,-.5) to (0.3,.5) 
node at (-.4,.45) {$\sigma^*$};
\coordinate [label=below right:{\small $\Pi_q(p)$}] (X) at (.3,.4);
\draw [-latex, bend right, very thick] (-1,-.5) to (0.3,.5) 
node at (-.2,-.05) {$\sigma$};
\draw [black, thick] plot [mark=*, only marks]
coordinates {(-1.,-.5) (.3,.5)};
\draw [dashed] (-2,0) to [out=50, in=130] (2,0);
\node at (0,-1.6) {$(\textbf{A})$};
\end{tikzpicture}
\qquad
\begin{tikzpicture}
\draw [very thick] (0,0) ellipse (2 and 1.2);
\coordinate [label=above left:$\Ma$] (t) at (1.5,.5);
\draw [->, very thick] (-1,-.5) -- (0,-.5);
\draw [->, very thick] (0.3,.5) -- (1.2,.4);
\node (p) at (-1.2,-.6) {$p$};
\node (q) at (0.3,.7) {$q$};
\coordinate [label=below right:{\small $\nihat^*(p,q)$}] (X) at (-1,-.6);
\draw [-latex, bend left, very thick] (-1,-.5) to (0.3,.5) 
node at (-.4,.45) {$\sigma$};
\coordinate [label=below right:{\small $\Pi^*_q(p)$}] (X) at (.3,.4);
\draw [-latex, bend right, very thick] (-1,-.5) to (0.3,.5) 
node at (-.2,-.05) {$\sigma^*$};
\draw [black, thick] plot [mark=*, only marks]
coordinates {(-1.,-.5) (.3,.5)};
\draw [dashed] (-2,0) to [out=50, in=130] (2,0);
\node at (0,-1.6) {$(\textbf{B})$};
\end{tikzpicture}
\\
\begin{tikzpicture}
\draw [xscale=1,very thick](0,0) ellipse (1 and 1.5);
\draw [dashed] (0,-1.5) to [out=120, in=210] (0,1.5);
\coordinate [label=left:$q$] (p) at (0,0.1);
\coordinate [label=right:$p$] (q) at (2,.5);
\draw[-latex, bend left, very thick] (2,.5) to [out=-40,in=-150] (0,0)
       node at (1.1,.8) {${\sigma^*}$}
       node at (1.1,-.3) {$\sigma$}
         node at (0,1) {\small $\hyperm_p(\kappa)$};
\draw[-latex, bend right, very thick] (2,.5) to [out=40, in=150] (0,0) ;
\draw [black, thick] plot [mark=*, only marks]
coordinates { (0,0) (2,0.5)};
 \draw [->, very thick] (0,0) -- (-.7,-.35)
  node at (-.2,-.6) {\small $\Pi+\Pi^*$};
  \node at (0,-1.9) {$(\textbf{C})$};
\end{tikzpicture}
\caption{From the top to the bottom, (\textbf{A}) isllustrates the vector $\Pi$ that is the $\nabla$-parallel transport of $\nihat(p,q)=\dot{\sigma}(0)$ along the $\nabla^*$-geodesic $\sigma^*$ from $p$ to $q$; while (\textbf{B}) illustrates the vector $\Pi^*$ that is the $\nabla^*$-parallel transport of $\nihat^*(p,q)=\dot{\sigma}^*(0)$ along the $\nabla$-geodesic $\sigma$ from $p$ to $q$. Finally, (\textbf{C}) shows that the sum $\Pi+\Pi^*$ is orthogonal to the level-hypersurface $\hyperm_p(\kappa)$ of constant pseudo-squared-distance $r_p(q)$.}\label{perpendicular}
\end{figure}

{The relevance of the vectors $\Pi$ and $\Pi^*$ { for the definition of our} new canonical divergence and its dual function relies on the next two results. The first one extends  Theorem \ref{thhgraddistance} to the context of Information Geometry, where a smooth manifold $\Ma$ is endowed with a dualistic structure $(\metric,\nabla,\nabla^*)$.}
\begin{theorem}\label{thmgradient}
Let $(\Ma,\metric,\nabla,\nabla^*)$ be a statistical manifold and $\U\subset\Ma$ be a dually convex set. For every $p,q\in\U$,  {consider the pseudo-squared-distance defined by
\begin{equation}
\label{pseudonorm}
r(p,q):=\langle\exp_p^{-1}(q),\stackrel{*}{\exp}_p^{-1}(q)\rangle_p\,,
\end{equation}
where $\exp$ and $\stackrel{*}{\exp}$ denote the exponential maps of $\nabla$ and $\nabla^*$, respectively.} Then we have
\begin{equation}
\label{gradpseudo}
\grad_q r_p =\Pi_q(p)+\Pi^*_q(p)\,.
\end{equation}
\end{theorem}
Notice that the function $r(p,q)$ is not symmetric in its argument. Therefore, it is not a distance in the classical sense. This justifies its classification as pseudo-distance. The proof of Theorem \ref{thmgradient} is obtained through the extension of the Gauss Lemma to Information Geometry which claims that the sum $\Pi+\Pi^*$ is orthogonal to the hypersurfaces of constant pseudo-squared-distance $r_p(q)$.
\begin{theorem}\label{thmorthog}
{Let  $(\Ma,\metric,\nabla,\nabla^*)$ be a statistical manifold and $\U\subset\Ma$ be a dually convex set. For $p\in\U$ let $\varepsilon>0$ so small that $\exp_p,\stackrel{*}{\exp}_p:\mathcal{E}_p(\varepsilon)\rightarrow\U$ are diffeomorphisms on their images, $\mathcal{E}_p(\varepsilon)=\{X\in\tangent_p\Ma\, |\, \sqrt{\langle X,X\rangle_p}<\varepsilon\}$. For $\kappa\in(0,\varepsilon)$ consider the hypersurface $\hyperm_p(\kappa)\subset\U$ centered at $p$ and defined by
$
\hyperm_p(\kappa):=\{q\in\U\, |\, r_p(q)=\kappa\}\,.
$
Then, for every $q\in\hyperm_p(\kappa)$, the sum $\Pi_q(p)+\Pi_q^*(p)$ is orthogonal to $\hyperm_p(\kappa)$ at $q$.}
\end{theorem}
{\it The geometry of geodesics is developed in  Section \ref{GeodGeometry} and  Theorems \ref{thmgradient} and \ref{thmorthog} are proved therein.}

{ Note that the set $\hyperm_p(\kappa)$ is a hypersurface within $\Ma$ in the sense that it is a properly embedded submanifold with codimension equal to $1$ (see Remark \ref{hyperH} in Section \ref{GeodGeometry}).}

\vspace{.3cm}

{Clearly, when $\nabla=\nabla^*$ the exponential map of $\nabla$ coincides with the one of $\nabla^*$. Therefore, in the self-dual case the pseudo-squared-distance $r(p,q)$ becomes the squared Riemannian distance, 
$$
r_p(q)=\left\langle\overline{\exp}^{-1}_p(q),\overline{\exp}^{-1}_p(q)\right\rangle_{p}=d(p,q)^2\,,
$$
where again $\overline{\exp}$ denotes the Levi-Civita exponential map. Obviously, also the $\nabla$ and $\nabla^*$ geodesics coincide. Thus, when the statistical manifold $(\Ma,\metric,\nabla,\nabla^*)$ is self-dual, the vectors $\Pi_q(p)$ and $\Pi_q^*(p)$ can be written as follows:
$$
\Pi_q(p)=\paralleltransport_{\sigma^*}\dot{\sigma}(0)=\dot{\bar{\sigma}}(1)= \paralleltransport^*_{\sigma}\dot{\sigma}^*(0)=\Pi_q^*(p)\,,
$$
where $\bar{\sigma}:[0,1]\rightarrow\Ma$ is the $\Lc$-geodesic such that $\bar{\sigma}(0)=p$ and $\bar{\sigma}(1)=q$. Applying Theorem \ref{thmgradient} to this particular case, and recalling Theorem \ref{thhgraddistance}, if $\bar{\sigma}$ is the $\Lc$-geodesic from $p$ to $q$,  we then get
$$
 \grad_q r_p=\left(\Pi_q(p)+\Pi_q^*(p)\right)=2\,\dot{\bar{\sigma}}(1)=\grad_q\, d^2_p\,.
$$
This proves, on one side, the consistency of our approach with the classical Riemannian theory. On the other hand, since a canonical divergence has to be one half the squared Riemannian distance, it suggests that the pseudo-squared-distance $r_p(q)$ holds information of both, the canonical divergence and its dual function. Indeed, $r(p,q)$ is obtained by summing up the canonical divergence and its dual function.

A further support to this claim can be found in the dually flat case. In general, a divergence function is not symmetric in its argument and this asymmetry property plays an important role for providing a dualistic structure on a smooth manifold according to Eq. \eqref{metricfromdiv} and Eq. \eqref{dualfromDiv} \cite{eguchi1992}.  The canonical divergence \eqref{BregmanDiv} defined on dually flat statistical manifolds paves the way along this avenue. It holds the following symmetry property,
\begin{equation}
\label{duallyflatsymmetry}
D[q:p]=D^*[p:q]\,,
\end{equation}
where $D^*[p:q]$ is the dual function of \eqref{BregmanDiv} \cite{Amari16}. This nice property implies that the  pseudo-squared-distance $ r(p,q) $ is obtained by summing up the canonical divergence and its dual function \cite{Amari00},
\begin{equation}
\label{pseudoBregman}
D[p:q]+D^*[p:q] = \langle\exp_p^{-1}(q),\stackrel{*}{\exp}_p^{-1}(q)\rangle_p=r(p,q)\,.
\end{equation}

\vspace{.3cm}

Inspired by Eq. \eqref{gradpseudo} and Eq. \eqref{pseudoBregman}, we attempt to define a canonical divergence $\Div$ as a potential function of $\Pi$ and its dual $\Div^*$ as a potential function of $\Pi^*$,
\begin{equation}
\label{potentialfunction}
\grad\,\Div\,=\,\Pi\,,\quad\grad\,\Div^*\,=\,\Pi^*\,.
\end{equation} 
This is a very natural requirement as we know from Theorem \ref{thmgradient} that the pseudo-squared-distance $r(p,q)$ is a potential function of the sum $\Pi_q(p)+\Pi^*_q(p)$ on a dually convex set $\U$. To explicitly show the latter claim, we introduce  from Eqs. \eqref{P}, \eqref{P*} two vector fields $\Pi_t(p),\,\Pi^*_t(p)\in\Tau(\gamma)$ along any arbitrary path $\gamma:[0,1]\rightarrow\U$ such that $\gamma(0)=p$ and $\gamma(1)=q$.  Since $\gamma(t)\in\mathrm{U}$  for every $t\in[0,1]$, we can find a $\nabla$-geodesic $\sigma_t(s)\,(0\leq s\leq 1)$ and a $\nabla^*$-geodesic $\sigma_t^*(s)\,(0\leq s\leq 1)$ such that $\sigma_t(0)=p=\sigma_t^*(0)$ and $\sigma_t(1)=\gamma(t)=\sigma_t^*(1)$. Hence, according to Eq. \eqref{P} and Eq. \eqref{P*}, we can write
\begin{align}
\Pi_{t}(p)= \paralleltransport_{\sigma^*_t} \nihat_p(\gamma(t)),&\qquad \nihat_p(\gamma(t))=\exp_p^{-1}(\gamma(t))\label{vectorfieldPiI}\\
\Pi_t^*(p)=\paralleltransport^*_{\sigma_t}\nihat_p^*(\gamma(t)),&\qquad \nihat_p^*(\gamma(t))=\stackrel{*}{\exp}_p^{-1}(\gamma(t)) \label{vectorfieldPI*I}\,,
\end{align}
where $\paralleltransport_{\sigma^*_t}:\tangent_p\Ma\rightarrow\tangent_{\gamma(t)}\Ma$ is the $\nabla$-parallel transport along $\sigma^*_t(s)$ and $\paralleltransport^*_{\sigma_t}:\tangent_p\Ma\rightarrow\tangent_{\gamma(t)}\Ma$ is the $\nabla^*$-parallel transport along $\sigma_t(s)$. Then, after some computations (for more details see Section \ref{TheoremPi&Pi*orthogonal}) we obtain from Eq. \eqref{gradpseudo} that the sum
\begin{equation}
\label{indipendenceI}
\int_0^1\ \langle \Pi_t(p),\dot{\gamma}(t)\rangle_{\gamma(t)}\,\total t+\int_0^1\  \langle \Pi_t^*(p),\dot{\gamma}(t)\rangle_{\gamma(t)}\,\total t= r_p(q)
\end{equation}
is independent of the particular path from $p$ to $q$.

In view of the theories developed in  \cite{Ay15} and \cite{Kobayashi00}, the geometry of geodesics of a general statistical manifold carried out through the extension of the Gauss Lemma, i.e. Theorems \ref{thmgradient} and \ref{thmorthog}, suggests to introduce a new divergence and its dual function by the path integrals of the vector fields $\Pi$ and $\Pi^*$. In particular, we define the divergence $\Div(p,q)$ as the path integral of $\Pi_t(p)$ along a $\nabla$-geodesic.

\begin{definition}
\label{divergenceIdefi}
Let $(\Ma,\metric,\nabla,\nabla^*)$ be a statistical manifold and $\U$ be a dually convex set. For every $p,q\in\U$, we can consider a $\nabla$-geodesic $\sigma:[0,1]\rightarrow\U$ such that $\sigma(0)=p$ and $\sigma(1)=q$. Then from  Eq. \eqref{vectorfieldPiI}, we can set 
$$
\Pi_t(p)=\paralleltransport_{\sigma_t^*}\,\nihat_p(t),\quad \nihat_p(t)=\exp_p^{-1}(\sigma(t))\,,
$$
where $\paralleltransport_{\sigma^*_t}:\tangent_p\Ma\rightarrow\tangent_{\sigma(t)}\Ma$ is the $\nabla$-parallel transport along the $\nabla^*$-geodesic $\sigma_t^*$ such that $\sigma_t^*(0)=p$ and $\sigma_t^*(1)=\sigma(t)$. We define the  function $\Div:\mathrm{U}\times\U\rightarrow\RR $  by the path integration of the vector field $\Pi_t(p)$ along $\sigma$,
\begin{equation}
\label{divergenceI}
\Div(p,q):=\int_0^1\ \langle\Pi_t(p),\dot{\sigma}(t)\rangle_{\sigma(t)}\ \total t\ .
\end{equation}
We refer to $\Div(p,q)$ as {a} {canonical divergence} {on a dually convex set $\U$} on $\Ma$ from $p$ to $q$.
\end{definition}

Analogously, we define the dual divergence $\Div^*(p,q)$ as the path integral of $\Pi^*_t(p)$ along a $\nabla^*$-geodesic. Then, consider a $\nabla^*$-geodesic $\sigma^*:[0.1]\rightarrow\U$ such that $\sigma^*(0)=p$ and $\sigma^*(1)=q$. Thus, according to Eq. \eqref{vectorfieldPI*I}, we can set
$$
\Pi^*_t(p)=\paralleltransport^*_{\sigma_t}\,\nihat^*_p(t),\quad \nihat^*_p(t)=\stackrel{*}{\exp}_p^{-1}(\sigma^*(t))\,,
$$
where $\paralleltransport^*_{\sigma_t}:\tangent_p\Ma\rightarrow\tangent_{\sigma^*(t)}\Ma$ is the $\nabla^*$-parallel transport along the $\nabla$-geodesic $\sigma_t$ such that $\sigma_t(0)=p$ and $\sigma_t(1)=\sigma^*(t)$. Then, we define the dual function $\Div^*:\mathrm{U}\times \U\rightarrow\RR$ by the path integration of the vector field $\Pi_t^*(p)$ along the $\nabla^*$-geodesic $\sigma^*$,
\begin{equation}
\label{divergence*I}
\Div^*(p,q):=\int_0^1\ \langle\Pi^*_t(p),\dot{\sigma}^*(t)\rangle_{\sigma^*(t)}\ \total t\ .
\end{equation}
 We refer to $\Div^*(p,q)$ as {the} {dual divergence} of $\Div(p,q)$ on $\Ma$ from $p$ to $q$.

\begin{figure}[h!]\label{actionfunctional}
\begin{tikzpicture}
\draw[thick]  (0,0) to [out=60, in=120] (5,3);
\coordinate [label=below left:$p$]
(p) at (0,0);
\coordinate [label=right:$q$] (q) at (5,3);
\coordinate [label=above left:$\sigma(t)$] (t) at (1.8,2.2);
\draw[very thick] (0,0) to [out=-30,in=-60] (2,2.55)
       node at (2.2,.8) {$\sigma^*_t$}
node at (0.5,1.3) {$\sigma_t$} ;
\draw [black, thick] plot [mark=*, only marks]
coordinates {(0,0) (5,3) (2,2.55)};
\draw [->,ultra thick] (2,2.55) -- (1.5,3);
\coordinate [label=below left:{\small $\Pi_t(p)$}] (P) at (2,3.6); 
\draw [->,ultra thick] (0,0) -- (.5,.8);
\coordinate [label=left:{\small $\nihat_p(t)$}] (P) at (0.3,0.6); 
\end{tikzpicture}
\quad
\begin{tikzpicture}
\draw[thick]  (0,0) to [out=60, in=120] (5,3);
\coordinate [label=below left:$p$]
(p) at (0,0);
\coordinate [label=right:$q$] (q) at (5,3);
\coordinate [label=above left:$\sigma^*(t)$] (t) at (1.8,2.2);
\draw[very thick] (0,0) to [out=-30,in=-60] (2,2.55)
       node at (2.2,.8) {$\sigma_t$}
node at (0.5,1.3) {$\sigma^*_t$} ;
\draw [black, thick] plot [mark=*, only marks]
coordinates {(0,0) (5,3) (2,2.55)};
\draw [->,ultra thick] (2,2.55) -- (1.5,3);
\coordinate [label=below left:{\small $\Pi^*_t(p)$}] (P) at (2,3.6); 
\draw [->,ultra thick] (0,0) -- (.5,.8);
\coordinate [label=left:{\small $\nihat^*_p(t)$}] (P) at (0.3,0.6); 
\end{tikzpicture}
\caption{On the left side of the figure,  the $\nabla$-geodesic $\sigma$ connects $p$ with $q$. The vector field $\Pi_t(p)\in\Tau(\sigma)$  is obtained by {$\nabla$-}parallel translating the vector $\nihat_p(t)=\exp^{-1}_p(\sigma(t))$ along the $\nabla^*$-geodesic $\sigma^*_t$ that connects $p$ with $\sigma(t)$.  On the right side of the figure, the $\nabla^*$-geodesic $\sigma^*$ connects $p$ with $q$. The vector field $\Pi^*_t(p)\in\Tau(\sigma^*)$  is obtained by {$\nabla^*$}-parallel translating the vector $\nihat^*_p(t)=\stackrel{*}{\exp}_p^{-1}(\sigma^*(t))$ along the $\nabla$-geodesic $\sigma_t$ that connects $p$ with $\sigma^*(t)$.}
\end{figure}
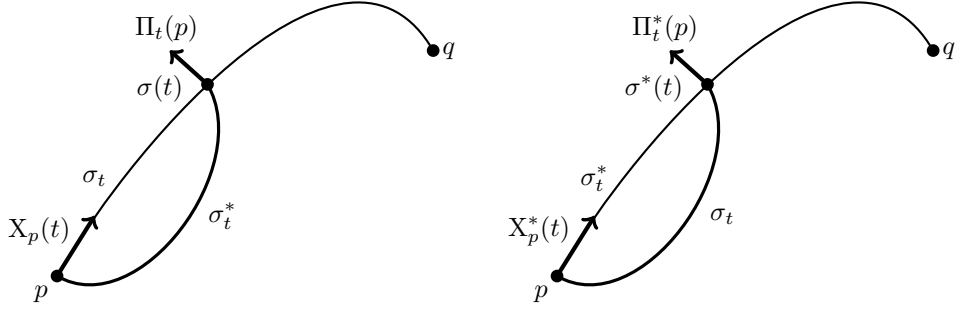

{ We emphasize that within Eqs. \eqref{divergenceI} and \eqref{divergence*I} we have made a particular choice. Specifically, we integrate $\Pi_t$ along a $\nabla$-geodesic and $\Pi^*_t$ along a $\nabla^*$-geodesic in order to get the canonical divergence and its dual function, respectively. If $\Pi_t$ and $\Pi^*_t$ are integrable, then the canonical divergence \eqref{divergenceI} and the dual divergence \eqref{divergence*I} would be independent of the particular path from $p$ to $q$. Unfortunately,} it turns out that, at least in general, the canonical divergence $\Div(p,q)$ is not a potential function of $\Pi_q(p)$, i.e.  $\grad_q\,\Div(p,\cdot)\neq\Pi_q(p)$, as well as $\Div^*(p,q)$ is not a potential function of $\Pi_q^*(p)$, i.e. $\grad_q\,\Div^*(p,\cdot)\neq\Pi^*_q(p)$. However, we succeed to supply an orthogonal decomposition of $\Pi$ in terms of the canonical divergence gradient  and a vector field that is orthogonal to $\nabla$-geodesics. On the other hand, we provide an orthogonal decomposition of $\Pi^*$ in terms of the dual divergence gradient  and a vector field that is orthogonal to $\nabla^*$-geodesics.

\begin{theorem}
\label{Pi&Pi*orthogonal}
Let $(\Ma,\metric,\nabla,\nabla^*)$ be a statistical manifold and $\U\subset\Ma$ be a dually convex set. For $p,q\in\U$ consider the canonical divergence $\Div(p,q)$ and its dual divergence $\Div^*(p,q)$. Let $\Pi_q(p)$ and $\Pi^*_q(p)$ given by Eq. \eqref{P} and Eq. \eqref{P*}, respectively. Then we have
\begin{align}
\label{Piorth}
& \Pi_q(p)=\grad_q\,\Div_p+X_q, \qquad X_q\in\tangent_q\Ma\quad \mbox{and}\quad \langle X_q,\dot{\sigma}(1)\rangle_q=0\ ,\\
\label{Pi*orth}
& \Pi^*_q(p)=\grad_q\,\Div^*_p+X^*_q, \qquad X^*_q\in\tangent_q\Ma\quad \mbox{and}\quad \langle X^*_q,\dot{\sigma}^*(1)\rangle_q=0\,.
\end{align}
Here, $\sigma,\sigma^*:[0,1]\rightarrow\U$ are the $\nabla$-geodesic and the $\nabla^*$-geodesic, respectively, such that $\sigma(0)=\sigma^*(0)=p$ and $\sigma(1)=\sigma^*(1)=q$. Furthermore, the decompositions \eqref{Piorth} and \eqref{Pi*orth} are orthogonal ones in the sense that $\langle\grad_q\,\Div_p,X_q\rangle_q=0$ and $\langle\grad_q\,\Div^*_p,X^*_q\rangle_q=0$ for all $q\in\U$.
\end{theorem}
\noindent{\it The proof of Theorem \ref{Pi&Pi*orthogonal} is presented in Section \ref{TheoremPi&Pi*orthogonal} and it is obtained therein within several steps.}

In order to prove Theorem \ref{Pi&Pi*orthogonal} we are now going to introduce two further functions ({\it Phi-functions} herein) which will turn out to be complementary, in some specific sense, to the canonical divergences $\Div$ and $\Div^*$.
Let $(\Ma,\metric,\nabla,\nabla^*)$ be a statistical manifold and $\U$ be a dually convex set. For $p,q\in\U$ we define the function $\varphi:\mathrm{U}\times\U\rightarrow\RR$ by the path integral of $\Pi_t(p)$ along the $\nabla^*$-geodesic $\sigma^*(t)$ from $p$ to $q$,
\begin{equation}
\label{phifunctionI}
\varphi(p,q):= \int_0^1 \langle \Pi_t(p),\dot{\sigma}^*(t)\rangle_{\sigma^*(t)}\ \total t\,,\quad \Pi_t(p)=\paralleltransport_{\sigma^*_t}\exp_p^{-1}(\sigma^*(t))\,,
\end{equation}
where $\sigma^*_t(s)\,(0\leq s\leq 1)$ is the $\nabla^*$-geodesic such that $\sigma^*_t(0)=p$ and $\sigma^*_t(1)=\sigma^*(t)$ and $\paralleltransport_{\sigma^*_t}:\tangent_p\Ma\rightarrow\tangent_{\sigma^*(t)}\Ma$ is the $\nabla$-parallel transport along $\sigma_t^*(s)$.
The dual function $\varphi^*:\U\times\U\rightarrow\RR$ is instead defined by the path integral of $\Pi_t^*(p)$ along the $\nabla$-geodesic $\sigma(t)$ from $p$ to $q$,
\begin{equation}
\label{phifunction*I}
\varphi^*(p,q):= \int_0^1 \langle \Pi^*_t(p),\dot{\sigma}(t)\rangle_{\sigma(t)}\ \total t\,,\quad \Pi^*_t(p)=\paralleltransport^*_{\sigma_t}\stackrel{*}{\exp}_p^{-1}(\sigma(t))\,,
\end{equation}
where $\sigma_t(s)\,(0\leq s\leq 1)$ is the $\nabla$-geodesic such that $\sigma_t(0)=p$ and $\sigma_t(1)=\sigma(t)$ and $\paralleltransport^*_{\sigma_t}:\tangent_p\Ma\rightarrow\tangent_{\sigma(t)}\Ma$ is the $\nabla^*$-parallel transport along $\sigma_t(s)$.
We refer to $\varphi(p,q)$ and $\varphi^*(p,q)$ as {\it \bf Phi-functions}.

The relevance of the Phi-functions emerges in this paper from the following statement:``the parallel transport of pairs of vectors with respect to a pair of dual connections is 'isometric' in the sense that inner product is preserved" \cite{Lauritzen87}. Indeed, thanks to this nice property held within a dualistic structure $(\metric,\nabla,\nabla^*)$, we get the following representation of the Phi-function $\varphi(p,q)$:
$$
\varphi(p,q)=\int_0^1\ \left\langle \stackrel{*}{\exp}_p^{-1}(\sigma^*(t)),\exp_p^{-1}(\sigma^*(t))\right\rangle_p\ \total t\,,
$$
for $p,q$ in a dually convex set $\U$ (see Lemma \ref{Lemma} of  Section \ref{Theorem 6} for more details). This representation allows us to naturally characterize $\varphi(p,q)$ by the decomposition of $\Pi$ in terms of a gradient vector field and another vector field which is orthogonal to $\nabla^*$-geodesics. Indeed, both Phi-functions, $\varphi$ and $\varphi^*$, are uniquely determined by the statement claimed in the following theorem.
\begin{theorem}
\label{localdecompositionPiI}
Given a statistical manifold $(\Ma,\metric,\nabla,\nabla^*)$ and a dually convex set $\U$. For $p,q\in\U$ consider $\sigma^*(t)\,(0\leq t\leq 1)$ and $\sigma(t)\,(0\leq t\leq 1)$ the $\nabla^*$-geodesic and the $\nabla$-geodesic, respectively, such that $\sigma^*(0)=p=\sigma(0)$ and $\sigma^*(1)=q=\sigma(1)$. Then, we have 
\begin{align}
\label{locdecPiI}
& \Pi_q(p)=\grad_q \varphi_p+V_q,\qquad V_q\in\tangent_q\Ma\quad \mbox{and}\quad \langle V_q,\dot{\sigma}^*(1)\rangle_q=0\ ,\\
\label{locdecPi*I}
& \Pi^*_q(p)=\grad_q {\varphi}^*_p+V^*_q,\qquad V^*_q\in\tangent_q\Ma\quad \mbox{and}\quad \langle V^*_q,\dot{\sigma}(1)\rangle_q=0\,.
\end{align}

In addition, decomposition \eqref{locdecPiI} and decomposition \eqref{locdecPi*I} are unique in $\mathrm{U}$.
\end{theorem}
{\it The proof of Theorem \ref{localdecompositionPiI} can be found in  Section \ref{Theorem 6}}.

\vspace{.2cm}

Both functions, $\varphi(p,q)$ and $\Div(p,q)$, are divergence (or contrast) functions according to the theory by Eguchi \cite{eguchi1983}.
\begin{theorem}\label{positivityI}
Consider a statistical manifold $(\Ma,\metric,\nabla,\nabla^*)$ and a dually convex set $\U\subset\Ma$. Then the canonical divergence $\Div(p,q)$ and the {\it Phi-function} $\varphi(p,q)$ are both non-negative in a neighborhood of the diagonal set $\Delta$ of $\U\times\U$ and vanish only on $\Delta$. Furthermore, they both induce the dual structure $(\metric,\nabla,\nabla^*)$ of $\Ma$ according to Eqs. \eqref{metricfromdiv} and \eqref{dualfromDiv}.
\end{theorem}
\noindent {\it The proof of Theorem \ref{positivityI} will be based on Taylor expansions of $\varphi$ and $\Div$,  presented in Propositions \ref{proptaylor} and \ref{TaylorD}.}

\vspace{.2cm}

Actually, $\Div(p,q)$ and $\varphi(p,q)$ differ from each other by the decompositions \eqref{Piorth} and \eqref{locdecPiI} as well as $\Div^*(p,q)$ and $\varphi^*(p,q)$ by the decompositions \eqref{Pi*orth} and \eqref{locdecPi*I}.  In fact, by the theory of minimum contrast geometry of Eguchi \cite{eguchi1992}, we have that $\grad_q\,\Div_p$ and $\grad_q\,\varphi_p$ are parallel to the  vector $\dot{\sigma}(1)$ whereas $\grad_q\,\Div^*_p$ and $\grad_q\,\varphi^*_p$ are parallel to the  vector $\dot{\sigma}^*(1)$. Therefore, it turns out that \eqref{Piorth} and \eqref{Pi*orth} are orthogonal decompositions while \eqref{locdecPiI} and \eqref{locdecPi*I} are not (see the end of Section \ref{ConclusionTh5} for more details).

In order to prove Theorem \ref{Pi&Pi*orthogonal}, we exploit the complementarity among the canonical divergences $\Div$ and $\Div^*$ and the Phi-functions $\varphi(p,q)$ and $\varphi^*(p,q)$ which appears in the decomposition of the pseudo-squared-distance $r(p,q)$.
Indeed, for $p,\,q$ in a dually convex set $\U\subset\Ma$, we may notice from Eq. \eqref{indipendenceI} and Eqs. \eqref{phifunction*I}, \eqref{divergenceI} that
$$
r_p(q)=\Div_p(q)+\varphi_p^*(q)\quad\mbox{and}\quad \grad_q\,r_p=\grad_q\,\Div_p +\grad_q\,\varphi_p^*\,.
$$ 
Therefore, by means of Eq. \eqref{gradpseudo} and Eq. \eqref{locdecPi*I} we get the decomposition \eqref{Piorth} by defining $X_q:=-V^*_q$.  Analogously, from Eq. \eqref{indipendenceI} and Eqs. \eqref{phifunctionI}, \eqref{divergence*I} we get
$$
r_p(q)=\Div^*_p(q)+\varphi_p(q)\quad\mbox{and}\quad \grad_q\,r_p=\grad_q\,\Div^*_p +\grad_q\,\varphi_p\,,
$$ 
for every $p,\,q\in\U$. Thus, by Eq. \eqref{gradpseudo} and Eq. \eqref{locdecPiI} we obtain the decomposition \eqref{Pi*orth} by defining $X^*_q:=-V_q$.

\vspace{.2cm}

Following a general theorem of Matumoto \cite{matumoto1993}, the Phi-function $\varphi(p,q)$ and the canonical divergence $\Div(p,q)$ coincide up to the third order in their Taylor expansion for $p$ and $q$ sufficiently close to each other. However, we prove in this paper that a further relation  between $\Div(p,q)$ and $\varphi(p,q)$ holds. Indeed, we show  that there exist  functions $\Upsilon,\Upsilon^*:[0,K]\rightarrow\RR^+\,(K>0)$ satisfying $\Upsilon(0)=0=\Upsilon^*(0)$ and $\Upsilon^{\prime}(0),\,\Upsilon^{\prime}(0)>0$ (here ``$\,\prime\,$" denotes the derivative operation) such that
$$
\varphi(p,q)=\Upsilon(\Div(p,q)),\qquad \varphi^*(p,q)=\Upsilon^*(\Div^*(p,q))
$$
for every $p,q$ in a dually convex set $\U$ (see Section \ref{SymmetryProperties} for more details). Although this proves a slight difference between $\Div(p,q)$ and $\varphi(p,q)$, we strengthen the reasons why we refer to $\Div(p,q)$ as the canonical divergence instead of the Phi-function $\varphi(p,q)$ by showing that $\Div(p,q)$ is one half the squared Riemannian distance on self-dual statistical manifolds and it reduces to the canonical divergence \eqref{BregmanDiv} of Bregman type on dually flat manifolds. Furthermore, beyond the characterization of $\Div(p,q)$ in terms of  the orthogonal decomposition \eqref{Piorth} we use  the divergence \eqref{Aydivergence} of Ay and Amari and the divergence \eqref{HenmiDivergence} of Henmi and Kobayashi as benchmark of our proposal. The relevance of the divergence \eqref{Aydivergence} in Information Geometry is above all based on the geodesic projection property. On the contrary, the divergence of Henmi and Kobayashi holds a nice potential theoretic property on the class of statistical manifolds satisfying the condition $(S)$:
\begin{equation}
\label{Scondition}
(S)\,=\,\left\{\begin{array}{ll}
(i) &\qquad  \RC(X,Y,Y,Y)=0\quad \forall\ X,Y\in\Tau(\Ma)\\
(ii) & \qquad  \nabla\ \RC =0,
\end{array}\right\}
\end{equation}
where $\RC$ denotes the Riemann curvature tensor of $\nabla$ (see Section \ref{SymmetricStatisticalmanifold} for more details). The dual condition $(S)^*$ is obtained by interchanging the role of $\nabla$ and $\RC$ with the $\nabla^*$-connection and the Riemann curvature tensor $\RC^*$ of $\nabla^*$. We refer to a statistical manifold $(\Ma,\metric,\nabla,\nabla^*)$ satisfying $(S)$ and $(S)^*$ as {\it symmetric statistical manifold}. {Note that, in particular, a dually flat manifold satisfies condition \eqref{Scondition}.  In this paper}, we succeed to prove that, on symmetric statistical manifolds, $\Div(p,q)$ coincides with the divergence $D(p,q)$ of Ay and Amari as well as with the dual divergence $W^*(q\|p)$ of Henmi and Kobayashi.
\begin{theorem}
\label{FeaturesDiv}
Let $(\Ma,\metric,\nabla,\nabla^*)$ be a statistical manifold and $\U\subset\Ma$ be a dually convex set. For $p,q\in\U$ consider the canonical divergence $\Div(p,q)$ defined by Eq. \eqref{divergenceI}. Then, the following properties hold true.
\begin{itemize}
\item If $\nabla=\nabla^*=\overline{\nabla}_{\mbox{\rm LC}}$, then $\Div(p,q)=\frac{1}{2}\, d(p,q)^2$, were $d(p,q)$ is the Riemannian distance \eqref{Riemanniandistance}.
\item If $\RC(\nabla)=0=\RC^*(\nabla^*)$, then $\Div(p,q)=D[p:q]$, where $D[p:q]$ is the divergence \eqref{BregmanDiv} of Bregman type.
\item If the condition \eqref{Scondition} holds, then $\Div(p,q)=D(p,q)$, where $D(p,q)$ is defined by Eq. \eqref{Aydivergence}. In addition, $\Div(p,q)=W^*(q\|p)$, where $W^*(q\|p)$ is the dual divergence of $W(p\|q)$ given in Eq. \eqref{HenmiDivergence}.
\end{itemize}
\end{theorem}
\noindent{\it The proof of Theorem \ref{FeaturesDiv} is going to be presented in Sections \ref{SelfDual&DuallyFlat} and \ref{SymmetricStatisticalmanifold}.}

{It turns out that the Phi-function $\varphi(p,q)$ and the canonical divergence $\Div(p,q)$ coincide in the dually flat case. However, we will show in Section \ref{SymmetricStatisticalmanifold} (Remark \ref{PhiinStatMan}) that they are not necessarily equal when a statistical manifold $(\Ma,\metric,\nabla,\nabla^*)$ holds the more general property \eqref{Scondition}. As a consequence, in comparison with the Phi-function, the function $\Div(p,q)$ appears to be a more natural choice towards the definition of a canonical divergence that is in line with} the divergences introduced in \cite{Ay15} and \cite{Kobayashi00}. 

\vspace{.2cm}

In this paper, we also address the issue of the symmetry property \eqref{duallyflatsymmetry} originally established by the canonical divergence of Bregamn type on dually flat manifolds. In general, only a weak version of it is true. However, the analysis carried out around the next theorem allows us to conjecture that this nice property is closely related to the decomposition of the pseudo-squared-distance $r(p,q)$ in terms of $\Div(p,q)$ and $\Div^*(p,q)$.
\begin{theorem}\label{ThmD*&D}
Let $(\Ma,\metric,\nabla,\nabla^*)$ be a statistical manifold and $\U$ be a dually convex set. For $p,q\in\U$ consider the canonical divergence $\Div(q,p)$ and the dual divergence $\Div^*(p,q)$ (resp. $\Div^*(q,p)$ and $\Div(p,q)$). Then, there exists a function $f$ (resp. $f^*$) satisfying the conditions $f(0)=0$ and $f^{\prime}(0)>0$ (resp. $f^*(0)=0$ and $f^{*\prime}(0)>0$) such that
\begin{equation}\label{D&D*}
\Div(q,p)=f\left(\Div^*(p,q)\right)\qquad \left(\mbox{resp.}\;\Div^*(q,p)=f^*\left(\Div(p,q)\right)\right) \ .
\end{equation}
\end{theorem}
\noindent{\it The proof of Theorem \ref{ThmD*&D} will be presented in Section \ref{SymmetryProperties}.}

\newpage

\begin{center}
\title{\bf \large Technical details and proofs}\label{Proofs}
\end{center}

\section{Geometry of geodesics in Information Geometry}\label{GeodGeometry}

Given a statistical manifold $(\Ma,\metric,\nabla,\nabla^*)$ we can recover the Levi-Civita connection by averaging the dual connections $\nabla$ and $\nabla^*$ \cite{Amari16},
\begin{equation}
\label{LeviCivita}
{\Lc}=\frac{1}{2}\left(\nabla+\nabla^*\right)\ .
\end{equation}

Thus, a statistical manifold can be understood as a generalization of a Riemannian manifold. Indeed, when $\nabla=\nabla^*$ the quadruple $(\Ma,\metric,\nabla,\nabla^*)$ reduces to the Riemannian manifold $(\Ma,\metric,\Lc)$ endowed with the Levi-Civita connection $\Lc$. In this case, the geodesic structure of $(\Ma,\metric,\Lc)$ is fully understood thanks to the celebrated Gauss Lemma, which is stated in  Theorem \ref{GaussLemma} of  Section \ref{Overwiev}. 

The proof of this classical result is crucially based on the claim that the function
\begin{equation}
\label{LCparallel}
t\mapsto\left\langle \frac{\total}{\total t}\overline{\exp}_p\left(t X_p\right),\frac{\total}{\total t}\overline{\exp}_p\left(t X_p\right)\right\rangle_{\overline{\exp}_p(tX_p)},
\end{equation}
is constant with respect to $t$ for all $X_p$ in a neighborhood $\overline{\mathcal{E}}_{p}(\varepsilon)\,(\varepsilon>0)$ of the null vector $O_p\in\tangent_p\Ma$ such that $\overline{\exp}_p:\overline{\mathcal{E}}_{p}(\varepsilon)\rightarrow\Ma$ is a diffeomorphism on to its image. To be more precise, in the energy of the $\Lc$-geodesic $\bar{\sigma}(t):=\overline{\exp}_p(t\,X_p)$,
$$
E(\bar{\sigma})=\frac{1}{2}\,\int_0^1\,\|\dot{\bar{\sigma}}(t)\|^2_{\bar{\sigma}(t)}\,\total t\,,
$$
the integrand is constant. Then, the Gauss Lemma is obtained by considering the first geodesic variation of $E(\bar{\sigma})$ \cite{Lee97}.
 
Actually, the Gauss Lemma tells us even more. Indeed, we can use it to prove that the $\Lc$-geodesics are related to the gradient of the squared Riemannian distance gradient as claimed in  Theorem \ref{thhgraddistance} of  Section \ref{Overwiev}. In particular, given for $q\in\overline{\exp}_p(\overline{\mathcal{E}}_p(\varepsilon))$ the $\Lc$-geodesic $\bar{\sigma}(t)\,(0\leq t\leq 1)$ connecting $p$ with $q$, we have
\begin{equation}
\label{LCgradient}
\grad_q \ d^2_p= 2\,\dot{\bar{\sigma}}(1),\quad d_p(q)=\sqrt{\langle\overline{\exp}_p^{-1}(q),\overline{\exp}_p^{-1}(q)\rangle_p}\,,
\end{equation}
where $\|\cdot\|_q$ denotes the norm induced by the metric tensor $\metric$. Notice that, since $\bar{\sigma}$ is a $\Lc$-geodesic, we can write
$$
\grad_q \ d^2_p= 2\,\overline{\paralleltransport}_{\bar{\sigma}}\overline{\nihat}_p(q)\,,\quad \mbox{with}\quad \overline{\nihat}_p(q)=\overline{\exp}_p^{-1}(q)\,,
$$
where $\overline{\paralleltransport}_{\bar{\sigma}}:\tangent_p\Ma\rightarrow\tangent_q\Ma$ denotes the $\Lc$-parallel transport along $\bar{\sigma}$.  As a consequence, for $q\in\overline{\mathrm{S}}_{\kappa}(p)$  in the geodesic sphere 
$$
\overline{\mathrm{S}}_{\kappa}(p):=\left\{q\in\overline{\exp}_p(\overline{\mathcal{E}}_p(\varepsilon))\,|\,\sqrt{\left\langle\overline{\exp}_p^{-1}(q),\overline{\exp}_p^{-1}(q)\right\rangle_p}=\kappa\right\}\,,\quad \kappa\in(0,\varepsilon)
$$
every vector ${X}_q\in \tangent_q\Ma$ can be decomposed in the following way,
\begin{equation}
\label{localdecomposition}
X_q=\lambda(q)\,\overline{\paralleltransport}_{\bar{\sigma}}\overline{\nihat}_p(q) + W_q=\widetilde{\lambda}(q)\,\grad_q\,d^2_p(q)+W_q,
\end{equation}
where $\lambda(q)$ is a coefficient depending on $q$, $\widetilde{\lambda}(q)=2\,\lambda(q)$ and $W_q$ is a tangent vector at $q$ to the geodesic sphere $\overline{\mathrm{S}}_{\kappa}(p)$ centered at $p$.

\begin{remark}
\label{StandardDivergence}
In \cite{Amari15} the authors proposed the function
\begin{equation}
D[p:q]:=\langle \exp^{-1}_p(q),\exp^{-1}_p(q)\rangle_p
\end{equation}
as the \textit{Standard Divergence} of the statistical manifold $(\Ma,\metric,\nabla,\nabla^*)$. Here, $\exp_p$ denotes the exponential map with respect to the $\nabla$-connection. In contrast to the Levi-Civita connection, the function in Eq. \eqref{LCparallel}, now computed by the $\nabla$-exponential map, is not constant with respect to $t$. However, when $\nabla=\nabla^*$ the standard divergence becomes the square of the Riemannian distance,
$$
D[p:q]=\langle \overline{\exp}^{-1}_p(q),\overline{\exp}^{-1}_p(q)\rangle_p=\left(d_p(q)\right)^2\,.
$$
Unfortunately, it turns out that this divergence is unable, at least in general, to recover the dual structure of $\Ma$.
\end{remark}

\vspace{.3cm}

\noindent In order to formulate an extension to the framework of Information Geometry, we basically hold two messages from the nice classic theory around the Gauss Lemma:
\begin{enumerate}
\item the $\Lc$-geodesics are orthogonal to the geodesic spheres;
\item the integral curves of the Riemannian distance gradient are $\Lc$-geodesics.
\end{enumerate}

Thus, to mimic the claim established about the Eq. \eqref{LCparallel} and develop our theory around the extension of the Gauss Lemma to Information Geometry, we need to select an appropriate function which turns out to be constant along a $\nabla$-geodesic or a $\nabla^*$-geodesic. In this way, we can define pseudo-energies of $\nabla$-geodesics and $\nabla^*$-geodesics and investigate their first variation.
Recall that in Information Geometry the following statement holds true: ``the parallel transport of pairs of vectors with respect to a pair of dual connections is `isometric' in the sense that inner product is preserved" \cite{Lauritzen87}. To prove this, consider $p,q\in\Ma$ and $\gamma(t):[0,1]\rightarrow\Ma$ an arbitrary path such that $\gamma(0)=p$ and $\gamma(1)=q$. Let $X(t)\in\Tau(\gamma)$ be a $\nabla$-parallel section and $Y(t)\in\Tau(\gamma)$ be a $\nabla^*$-parallel section along $\gamma$. Then, according to Eq. \eqref{dualconnection} we have that
$$
\dot{\gamma}\left\langle X(t),Y(t)\right\rangle_{\gamma(t)}=\left\langle\nabla_{\dot{\gamma}} X(t),Y(t)\right\rangle_{\gamma(t)}+\left\langle X(t),\nabla^*_{\dot{\gamma}}Y(t)\right\rangle_{\gamma(t)}=0\,,
$$
where $\dot{\gamma}\left\langle X(t),Y(t)\right\rangle_{\gamma(t)}=\frac{\total}{\total t}\left\langle X(t),Y(t)\right\rangle_{\gamma(t)}$. Therefore, for $p,q\in\U$ in a dually convex set we define the functional  over the set of paths connecting $p$ and $q$ as,
\begin{equation}
\label{pseudoenergy}
\Lag(\gamma):=\int_0^1 \ \langle \dot{\gamma}(t),\paralleltransport_{t} \nihat(p,q)\rangle_{\gamma(t)}\,\total t\,, \quad \nihat(p,q)=\exp_p^{-1}(q)
\end{equation}
where $\paralleltransport_{t}:\tangent_p\Ma\rightarrow\tangent_{\gamma(t)}\Ma$ is the $\nabla$-parallel transport along $\gamma$. Note that $\Lag(\gamma)$ is not positive in general. For this reason, we refer to $\Lag(\gamma)$ as the {\it $\nabla$-pseudo-energy} of $\gamma$. When $\gamma$ is a $\nabla^*$-geodesic we have that the integrand of $\Lag$ is constant and $\Lag$ assumes a very useful form.

\begin{pro}\label{EnergyProp}
Let $p,q\in\U$ in a dually convex set. Consider the $\nabla^*$-geodesic $\sigma^*:[0,1]\rightarrow\U$  such that $\sigma^*(0)=p$ and $\sigma^*(1)=q$. Then
\begin{equation}
\label{Lagconst}
\Lag(\sigma^*)=\langle \nihat^*(p,q),\nihat(p,q)\rangle_p=r(p,q)\,,
\end{equation}
where $\nihat^*(p,q)={\stackrel{*}{\exp}_p}^{-1}(q)$ and $r(p,q)$ is the pseudo-squared-distance defined by Eq. \eqref{pseudonorm}.
\end{pro}

\noindent {\bf Proof}. Consider the map 
$$t\mapsto\langle \dot{\sigma}^*(t),\paralleltransport_{t}\nihat(p,q)\rangle_{\sigma^*(t)},$$
where $\paralleltransport_{t}:\tangent_p\Ma\rightarrow\tangent_{\sigma^*(t)}\Ma$ denotes the $\nabla$-parallel transport along $\sigma^*(t)$.
Then, by taking the derivative with respect to $t$ it follows from Eq. \eqref{dualconnection} that
$$
\frac{\total}{\total t}\langle \dot{\sigma}^*(t),\paralleltransport_{t}\nihat(p,q)\rangle_{\sigma^*(t)}=\langle \nabla_t  \paralleltransport_{t}\nihat(p,q),\dot{\sigma}^*(t)\rangle_{\sigma^*(t)}+\langle \nabla^*_t \dot{\sigma}^*(t), \paralleltransport_{t}\nihat(p,q)\rangle_{\sigma^*(t)}\ ,
$$ 
where $\nabla_t=\nabla_{\dot{\sigma}^*(t)}$ and $\nabla^*_t=\nabla^*_{\dot{\sigma}^*(t)}$ are the covariant derivatives of $\nabla$ and $\nabla^*$, respectively, on $\sigma^*$. By recalling that $\textbf\paralleltransport_{t}\nihat(p,q)$ is the $\nabla$-parallel transport along $\sigma^*$, we have $\nabla_t\paralleltransport_{t}\nihat(p,q)\equiv 0$. Analogously, we have $\nabla^*_t \dot{\sigma}^*(t)\equiv 0$ because $\sigma^*$ is a $\nabla^*$-geodesic. Therefore, we obtain 
$$
\frac{\total}{\total t}\langle \dot{\sigma}^*(t),\paralleltransport_{t}\nihat(p,q)\rangle_{\sigma^*(t)}=0,
$$
and  finally, we arrive at
$$
\langle\dot{\sigma}^*(t),\paralleltransport_{t}\nihat(p,q)\rangle_{\sigma^*(t)}=\langle\dot{\sigma}^*(0),\nihat(p,q)\rangle_{p}\ .
$$
Hence, we can conclude by noticing that $\dot{\sigma}^*(0)=\stackrel{*}{\exp}_p^{-1}(q)=\nihat^*(p,q)$.
 \hfill $\square$

\begin{remark}
The functional $\Lag$ can be also computed on a $\nabla$-geodesic $\sigma$ from $p$ to $q$. In this case, it assumes the following expression
\begin{equation}
\label{nablaLag}
\Lag(\sigma)=\int_0^1\|\dot{\sigma}(t)\|_{\sigma(t)}^2\ dt,
\end{equation}
where the integrand is now not constant with respect to $t$.
\end{remark}

\vspace{0.3cm}

Instead of geodesic spheres, we consider hypersurfaces of constant pseudo-squared-distance $r$ in a dually convex set $\U\subset\Ma$. Thus, let $p\in\U$,  we  define the set $\hyperm_p(\kappa)$  as follows
\begin{equation}
\label{hyperthm}
\hyperm_p(\kappa):=\{q\in\U \ |\ r(p,q)=\langle\nihat_p(q),\nihat^*_p(q)\rangle_p=\kappa\}\ ,
\end{equation}
where 
$$
\kappa>0,\quad \nihat_p(q)=\exp_p^{-1}(q), \qquad \nihat_p^*(q)={\stackrel{*}{\exp}}^{-1}_p(q)\ .
$$
\begin{remark}
\label{hyperH}
{ Let $p\in\mathrm{U}$ in a dually convex set. For $q\in\mathrm{U}-\{p\}$ there exists a vector $X_p(q)\in\mathcal{E}_p$ such that $q=\exp_p(X_p(q))$. In addition, there exists also a vector $X^*_p(q)\in\mathcal{E}_p$ such that $q=\stackrel{*}{\exp}_p(X^*_p(q))$. Since, both the exponential maps, the $\nabla$ and the $\nabla^*$ ones, are diffeomorphism from $\mathcal{E}_p$ onto their images, we can define $X_p(q):=\exp_p^{-1}(q)$ and $X_p^*(q):=\stackrel{*}{\exp}_p^{-1}(q)$. Clearly, the map $r_p:\Ma\rightarrow\RR,\,r_p(q):=\left\langle X_p(q),X_p^*(q)\right\rangle_p$ is a smooth one because it is a composition of smooth functions, namely $\metric(\cdot,\cdot)$, $\exp_p(\cdot)$ and $\stackrel{*}{\exp}_p(\cdot)$. Therefore, the gradient of $r_p(\cdot)$ is a vector field on $\mathrm{U}$ and the differential of $r_p$ is defined by the following relation:
$$
\left(\total r_p\right)_q(Y)=\left\langle\grad_q\,r_p(\cdot),Y\right\rangle_q\,.
$$
Obviously, this differential is surjective except at $q\neq p$. This implies that $q\neq p$ is a {\it regular} point of $r_p$ \cite{Lee}. Furthermore, we call $\kappa\in\RR$ a {\it regular value} of $r_p$ if $r_p^{-1}(\kappa)$ is a regular point. Now, a classical result in Riemannian geometry states that {\it every regular level set of a smooth map on a manifold $\Ma$ is a properly embedded submanifold whose codimension is the dimension of the codomain} \cite{Lee}.

We can then conclude that the set $\hyperm_p(\kappa)\subset\mathrm{U}$ is a submanifold of $\Ma$ with codimension $1$.}
\end{remark}
We refer to $\hyperm_p(\kappa)$ as {\it the hypersurface of constant pseudo-squared-distance centered at $p$}. We will shortly prove that the combination of vectors $\Pi_q(p)$ and $\Pi^*_q(p)$ of the Eqs. \eqref{P} and \eqref{P*} defines the rays of the hypersurface $\hyperm_p(\kappa)$. In particular, for $q\in\hyperm_p(\kappa)$ we will show that $\Pi_q(p)+\Pi^*_q(p)$ is orthogonal to $\hyperm_p(\kappa)$ at $q$.

\vspace{0.3cm}

Owing to the duality structure of the statistical manifold $(\Ma,\metric,\nabla,\nabla^*)$, we can introduce in $\tangent_p\Ma$ two notions of pseudo-spheres. 

\begin{definition}\label{defispheres}
{Let $\U\subset\Ma$ be a dually convex set, and let $p\in\U$. Consider the set $\mathcal{E}_p\subset\tangent_p\Ma$ such that $\exp_p,\,\stackrel{*}{\exp}_p:\mathcal{E}_p\rightarrow\Ma$ are diffeomorphisms on their images. For $X_p\in\mathcal{E}_p$ we can find  a $\nabla$-geodesic $\sigma(t)\,(0\leq t\leq 1)$ and a $\nabla^*$-geodesic $\sigma^*(t)\,(0\leq t\leq 1)$ such that
$$
\sigma(0)=p, \; \dot{\sigma}(0)=X_p \qquad \mbox{and}\qquad \sigma^*(0)=p, \; \dot{\sigma}^*(0)=X_p\ .
$$
Then we define 
\begin{equation}
\label{hyper}
\hyper_p(\kappa):=\left\{X_p\in\mathcal{E}_p\ |\  \langle {\exp}^{-1}_p(\stackrel{*}{\exp}_p(X_p)),X_p\rangle_p=\kappa \right\}
\end{equation}
and
\begin{equation}
\label{hyper*}
\hyper^*_p(\kappa):=\left\{X_p\in\mathcal{E}_p\ |\  \langle \stackrel{*}{\exp}_p^{-1}\left({\exp}_p(X_p)\right),X_p\rangle_p=\kappa \right\}.
\end{equation}}
\end{definition}

\begin{remark}\label{hyperimage}
The image of $\hyper_p(\kappa)$ through the $\nabla^*$-exponential map is given by
$$
\stackrel{*}{\exp}_p\left(\hyper_p(\kappa)\right)=\hyperm_p(\kappa)\,,
$$
where $\hyperm_p(\kappa)$ is defined in Eq. \eqref{hyperthm}.
 
Indeed, if we set $q=\stackrel{*}{\exp}_p(X_p)$ for some $X_p\in\hyper_p(\kappa)$, by definition we have $X_p=\stackrel{*}{\exp}_p^{-1}(q)$. Hence, from Eq. \eqref{hyper} we get $\langle\exp_p^{-1}(q),\stackrel{*}{\exp}_p^{-1}(q)\rangle_p=\kappa$, which proves that $q\in\hyperm_p(\kappa)$. Vice versa, if $q\in\hyperm_p(\kappa)$ we have $\langle\exp_p^{-1}(q),\stackrel{*}{\exp}_p^{-1}(q)\rangle_p=\kappa$. Since $\hyperm_p(\kappa)\subset\U$, we can find $X_p\in\mathcal{E}_p$ such that $q=\stackrel{*}{\exp}_p^{-1}(X_p)$, because $\hyperm_{p}(\kappa)\subset\mathrm{U}_p$. Therefore, we get $\langle \exp^{-1}_p(\stackrel{*}{\exp}_p(X_p)),X_p\rangle_p=\kappa$ which proves that $X_p\in\hyper_p(\kappa)$.

Analogously, we can prove that the action of the $\nabla$-exponential map on $\hyper^*_p(\kappa)$ gives
$$
{\exp}_p\left(\hyper^*_p(\kappa)\right)=\hyperm_p(\kappa),
$$  
where $\hyperm_p(\kappa)$ is defined by Eq. \eqref{hyperthm}.

{ Since, for $p,q\in\mathrm{U}$ in a dually convex set, bot the maps, $\exp_p$ and $\stackrel{*}{\exp}_p$, are diffeomorphisms onto their images, we can employ the same arguments as the ones in Remark \ref{hyperH} and conclude that both sets, $\hyper_p(\kappa)$ and $\hyper_p^*(\kappa)$, are hypersurfaces within $\tangent_p\Ma$.}
\end{remark}

The spheres of Definition \ref{defispheres} are not but almost the same object. Indeed, consider the map
\begin{align}
\label{Iso}
& I_p:\hyper_p(\kappa)\rightarrow\hyper^*_p(\kappa)\,,\nonumber\\
& I_p(X_p):=\exp_p^{-1}\left(\stackrel{*}{\exp}_p(X_p)\right), \; \forall \ X_p\in\hyper_p(\kappa)\ . 
\end{align}
Then we have,
\begin{pro}
{The map $I_p:\hyper_p(\kappa)\rightarrow\hyper^*_p(\kappa)$ defined by Eq. \eqref{Iso} is a  diffeomorphism}.
In addition, the following diagram 
\begin{center}
\begin{tikzpicture}
    \node (S) at (0,2) {$\hyper_p(\kappa)$};
    \node (H) at (2,2) {$\hyperm_p(\kappa)$};
    \node (Sd) at (0,0) {$\hyper^*_p(\kappa)$};
    \draw[->] (S)--(H) node [midway,above] {$\stackrel{*}{\exp}_p$};
    \draw[->] (S)--(Sd) node [midway,left] {$I_p$};
    \draw[->] (Sd)--(H) node [midway,below,right] {$\exp_p$};
 \end{tikzpicture}
\end{center}
is commutative.
\end{pro}
{\bf Proof}. Consider $X\in\hyper_p(\kappa)$. Firstly, we have that $I_p(X)\in\hyper^*_p(\kappa)$. Indeed, 
\begin{eqnarray*}
\langle \stackrel{*}{\exp}_p^{-1}\left(\exp_p(I_p(X))\right),I_p(X)\rangle_p&=& \langle X,\exp_p^{-1}\left(\stackrel{*}{\exp}_p(X)\right)\rangle_p\\
&=&\kappa \ .
\end{eqnarray*}
Consider now the map $\tilde{I}_p:\hyper^*_p(\kappa)\rightarrow \hyper_p(\kappa)$ defined by 
$$
\tilde{I}_p(X)=\stackrel{*}{\exp}_p^{-1}\left(\exp_p(X)\right)\ .
$$
Then, we can trivially see that
\begin{eqnarray*}
&&I_p\circ \tilde{I}_p(X)=X\ ,\quad \forall X\in\hyper_p(\kappa)\\
&&\tilde{I}_p\circ I_p(X)=X\ ,\quad \forall X\in\hyper_p^*(\kappa)\ .
\end{eqnarray*}
Therefore, we can conclude that $\tilde{I}_p=I^{-1}_p$.  In order to prove that the diagram is commutative, let us consider $q=\stackrel{*}{\exp}_p(X)$ by some $X\in\hyper_p(\kappa)$. From Remark \ref{hyperimage} we know that $q\in\hyperm_p(\kappa)$. In addition, by the definition \eqref{Iso} we also have that $q=\exp_p(I_p(X))$. This proves that $\stackrel{*}{\exp}_p\equiv\exp_p\circ I_p$. Finally, since both $\hyper_p(\kappa)$ and $\hyper^*_p(\kappa)$ are in the set $\mathcal{E}_p$ of the Definition \ref{defispheres}, we can conclude that $I_p$ is a diffeomorphism, as well. \hfill $\square$ 

\vspace{.5cm}

We now proceed to investigating the first variation of the $\nabla$-pseudo-energy  $\Lag$. In order to pursue this goal, let us firstly introduce the notion of path variation.
Given an arbitrary path $\gamma:[0,1]\rightarrow\Ma$ such that $\gamma(0)=p$ and $\gamma(1)=q$, {a continuous map $\Sigma:(-\varepsilon,\varepsilon)\times [0,1]\rightarrow\Ma$ is called a {\it variation} of $\gamma$ if $\Sigma(0,t)\equiv\gamma(t)$. In addition, we require that for any $s\in(-\varepsilon,\varepsilon)$ the {\it main curve} $\Sigma_s\equiv \Sigma(s,\cdot)$ is a smooth curve. Moreover, also the {\it transverse curve} $\Sigma^{(t)}\equiv\Sigma(\cdot,t)$ to the variation $\Sigma$ is a smooth curve for any $t\in[0,1]$. Finally, if $\gamma$ is a $\nabla$-geodesic, a variation $\Sigma$ is said to be a $\nabla$-{\it geodesic variation} of $\gamma$ if all the main curves $\Sigma_s(\cdot)$ are $\nabla$-geodesics. The same applies to a $\nabla^*$-geodesic and its $\nabla^*$-geodesic variation.}  

A  vector field along $\Sigma$ is a smooth map $\Xi:(-\varepsilon,\varepsilon)\times [0,1]\rightarrow\tangent\Ma$ such that $\Xi(s,t)\in\tangent_{\Sigma(s,t)}\Ma$ for each $(s,t)\in(-\varepsilon,\varepsilon)\times [0,1]$. Two very special vector fields {on $\Sigma$} are defined as follows
\begin{equation}
\label{SandT}
T(s,t):=\frac{\total}{\total t}\Sigma_s(t)=\partial_t\Sigma(s,t)\ , \quad S(s,t):=\frac{\total}{\total s}{\Sigma}^{(t)}(s)=\partial_s\Sigma(s,t)\ .
\end{equation}
{Clearly, $T(s,t)$ is the velocity vector field of the main curve $\Sigma_s(t)$ whereas $S(s,t)$ is the velocity vector field of the transverse curve $\Sigma^{(t)}(s)$.} Finally, $V(t)=\partial_s\Sigma(0,t)\in\Tau(\gamma)$ is called the {\it variation field}  of {the variation} $\Sigma$.

\vspace{.3cm}

Let $\Sigma(s,t)$ be a variation   of an arbitrary path $\gamma(t)\subset\mathrm{U}$ in a dually convex set $\U$ such that $\gamma(0)=p$ and $\gamma(1)=q$. For every $s\in(-\varepsilon,\varepsilon)$ we can consider the  vector,
\begin{equation}
\label{*}
\nihat_p(s)\equiv\nihat(p,\Sigma_s(1)):= \exp_p^{-1}(\Sigma_s(1))\ ,
\end{equation}
which is the velocity vector at $p$ of the $\nabla$-geodesic connecting $p$ and $\Sigma_s(1)$.
The first variation of $\Lag$ is then provided by the following Proposition.
\begin{pro}\label{FirstProp}
Let $p,q\in\U$ in a dually convex set $\U\subset\Ma$. Consider an arbitrary path  $\gamma:[0,1]\rightarrow\U$ such that $\gamma(0)=p$ and $\gamma(1)=q$ and a variation of $\gamma$, $\Sigma:(-\varepsilon,\varepsilon)\times [0,1]\rightarrow\Ma$. Let $V\in\Tau(\gamma)$ be the variation vector field of $\Sigma$. Finally, define the functional $\Lag(s):=\Lag(\Sigma_s)$. Then we have
\begin{equation}
\label{first}
\frac{\total\Lag}{\total s}(0)=\left. \langle V(t),\paralleltransport_{t}\nihat_p(q)\rangle_{\gamma(t)}\right|_0^1+\int_0^1\langle\dot{\gamma}(t),\nabla_V\paralleltransport_{t} \nihat_p(q)\rangle_{\gamma(t)} dt,
\end{equation}
where {$\paralleltransport_{t}:\tangent_p\Ma\rightarrow\tangent_{\gamma(t)}\Ma$} denotes the $\nabla$-parallel transport {along} the curve $\gamma(t)$.
\end{pro}
\noindent{\bf Proof.} Let us first see the definition of $\Lag$ evaluated at $\Sigma_s$. Recalling the definition of the vector field $T(s,t)$ given in Eq. \eqref{SandT} we have
\begin{eqnarray*}
\Lag(s)\equiv \Lag(\Sigma_s)&=&\int_0^1\left\langle\frac{\total\Sigma_s}{\total t}(t),\paralleltransport_{s,t}\nihat_p(s)\right\rangle_{\Sigma_s(t)}\,\total t\\
&=&\int_0^1\left\langle T(s,t),\paralleltransport_{s,t}\nihat_p(s)\right\rangle_{\Sigma_s(t)}\,\total t\, ,
\end{eqnarray*}
where {$\paralleltransport_{s,t}:\tangent_p\Ma\rightarrow\tangent_{\Sigma(s,t)}\Ma$} is the $\nabla$-parallel transport along the curve $\Sigma_s(t)$.

Therefore, by taking the derivative {with respect to $s$} and  exploiting the Eq. \eqref{dualconnection} we obtain
\begin{eqnarray}
\frac{\total\Lag(s)}{\total s}&=&\int_0^1\ \left(\langle \nabla^*_s T(s,t),\paralleltransport_{s,t}\nihat_p(s)\rangle_{\Sigma_s(t)}+ \langle  T(s,t),\nabla_s\paralleltransport_{s,t}\nihat_p(s)\rangle_{\Sigma_s(t)}\right)\ dt,\nonumber\\
\end{eqnarray}
where $\nabla_s=\nabla_{\dot{\Sigma}{^{(t)}}(s)}$ and $\nabla^*_s=\nabla^*_{\dot{\Sigma}{^{(t)}}(s)}$ are the covariant derivatives along $\Sigma^{(t)}(s)$ with respect to $\nabla$ and $\nabla^*$, respectively.  Since the connection $\nabla^*$ is torsion-free, we have that $\nabla^*_s \left(\partial_t \Sigma(s,t)\right)=\nabla^*_t \left(\partial_s\Sigma(s,t)\right)$, {which is equivalent to writing $\nabla^*_s T=\nabla^*_t S$, where $\nabla^*_t=\nabla^*_{\dot{\Sigma}_s(t)}$ is the $\nabla^*$-covariant derivative along the curve $\Sigma_s(t)$ \cite{Lee97}}. Therefore, {we can perform} the following computations:
\begin{eqnarray}
\frac{\total\Lag}{\total s}(s)&=&\int_0^1\ \left(\langle  \nabla^*_t S(s,t),\paralleltransport_{s,t}\nihat_p(s)\rangle_{\Sigma_s(t)}+ \langle  T(s,t),\nabla_s\paralleltransport_{s,t}\nihat_p(s)\rangle_{\Sigma_s(t)}\right)\,\total t\nonumber\\
&=& \int_0^1\ \left(\frac{\total}{\total t}\langle  S(s,t),\paralleltransport_{s,t}\nihat_p(s)\rangle_{\Sigma_s(t)}-\langle  S(s,t),\nabla_t\paralleltransport_{s,t}\nihat_p(s))\rangle_{\Sigma_s(t)}\right)\,\total t\nonumber\\
&&+ \int_0^1\ \left(\langle  T(s,t),\nabla_s\paralleltransport_{s,t}\nihat_p(s)\rangle_{\Sigma_s(t)}\right)\,\total t\nonumber\\
&=& \int_0^1\ \left(\frac{\total}{\total t}\langle  S(s,t),\paralleltransport_{s,t}\nihat_p(s)\rangle_{\Sigma_s(t)}+ \langle  T(s,t),\nabla_s\paralleltransport_{s,t}\nihat_p(s)\rangle_{\Sigma_s(t)}\right)\,\total t\ ,\nonumber
\end{eqnarray}
{where we exploited the property \eqref{dualconnection} and we used  $\nabla_t\paralleltransport_{s,t}\nihat_p(s)=0$. Here $\nabla_t=\nabla_{\dot{\Sigma}_s(t)}$ is the $\nabla$-covariant derivative along the curve $\Sigma_s(t)$.} Hence, we arrive at
\begin{equation}
\frac{\total\Lag}{\total s}(s)= \left.{\langle  S(s,t),\paralleltransport_{s,t}\nihat_p(s)\rangle_{\Sigma_s(t)}}\right|_{0}^1+\int_0^1 \langle  T(s,t),\nabla_s\paralleltransport_{s,t}\nihat_p(s)\rangle_{\Sigma_s(t)}\, \total t\ .
\label{firstvar}
\end{equation} 

Finally, setting $s=0$ and recalling that $S(0,t)=V(t)$, $T(0,t)=\dot{\gamma}(t)$, $\Sigma(0,t)=\gamma(t)$ and $\left.\nihat_p(s)\right|_{s=0}=\exp_p^{-1}(\Sigma(0,1))=\exp_p^{-1}(\gamma(1))=\nihat_p(q)$ we obtain that
\begin{equation}
\label{firstvariational}
\frac{d\Lag}{ds}(0)=\left. \langle V(t),\paralleltransport_{t}\nihat_p(q)\rangle_{\gamma(t)}\right|_0^1+\int_0^1\langle\dot{\gamma}(t),\nabla_V\paralleltransport_{t} \nihat_p(q)\rangle_{\gamma(t)}\, \total t.
\end{equation}

\hfill $\square$

\vspace{.2cm}

\vspace{.3cm}

We are now in the position to prove the extension of the Gauss Lemma to Information Geometry which was previously stated in Section \ref{Overwiev}. 

\noindent\textbf{\large Proof of Theorem \ref{thmorthog}.}\,  Let us consider a curve within the pseudo-sphere $\hyper_p(\kappa)\subset\tangent_p\Ma$, namely $\tau:(-\varepsilon,\varepsilon)\rightarrow\hyper_p(\kappa)$, such that  $\tau(0)=\nihat^*_p(q)=\stackrel{*}{\exp}_p^{-1}(q)$.
The map $\Sigma^*(s,t):=\stackrel{*}{\exp}_p(t\ \tau(s))$ is a $\nabla^*$-geodesic variation of the $\nabla^*$-geodesic\newline
$\sigma^*(t):=\stackrel{*}{\exp}_p(t\ \tau(0))$. {Notice that $\sigma^*(1)=\stackrel{*}{\exp}_p(\tau(0))=\stackrel{*}{\exp}_p(\nihat_p^*(q))=q$. In addition, for every $s\in(-\varepsilon,\varepsilon)$, we have $\Sigma^*(s,0)=\stackrel{*}{\exp}_p(O_p)=p$, where $O_P$ is the null vector in $\tangent_p\Ma$. Recall that the map $I_p:\hyper_p(\kappa)\rightarrow\hyper^*_p(\kappa)$ defined by
$$
I_p(X_p)=\exp_p^{-1}\left(\stackrel{*}{\exp}_p(X_p)\right),\quad \forall\, X_p\in\hyper_p(\kappa)\,,
$$
is a diffeomorphism. Hence, we can define a map $\Sigma:(-\varepsilon,\varepsilon)\times[0,1]\rightarrow\Ma$ by
\begin{equation}
\label{nabla-geod-var}
\Sigma(s,t):=\exp_p\left(t\,I_p(\tau(s))\right)\,,
\end{equation}
which turns out to be a $\nabla$-geodesic variation of the $\nabla$-geodesic connecting $p$ and $q$. In fact, recalling that $\tau(0)=\nihat_p^*(q)$, $\stackrel{*}{\exp}_p(\nihat_p^*(q))=q$ and $I_p(\nihat_p^*(q))=\nihat_p(q)=\exp_p^{-1}(q)$, we have that
\begin{eqnarray*}
\Sigma(0,t)&=&\exp_p\left(t\,I_p(\tau(0))\right)\\
&=&\exp_p\left(t\,I_p(\nihat_p^*(q))\right)\\
&=&\exp_p\left(t\,\exp_p^{-1}\left(\stackrel{*}{\exp}_p(\nihat_p^*(q))\right)\right)\\
&=&\exp_p\left(t\,\exp_p^{-1}(q)\right)\,.
\end{eqnarray*}
This proves that $\Sigma(0,t)=:\sigma(t)$ is the $\nabla$-geodesic connecting $p$ and $q$ and then, by definition, $\Sigma(s,t)$ is a $\nabla$-geodesic variation of $\sigma(t)$. In addition, by the definition of the map $I_p$, we also have 
\begin{equation}
\label{Variation}
\Sigma(s,1)=\exp_p(I_p(\tau(s)))=\exp_p\left(\exp_p^{-1}(\stackrel{*}{\exp}_p(\tau(s)))\right)=\Sigma^*(s,1)\, 
\end{equation}
for every $s\in(-\varepsilon,\varepsilon)$ (For a reference to the aforementioned construction, see Fig. \ref{Variations}).}

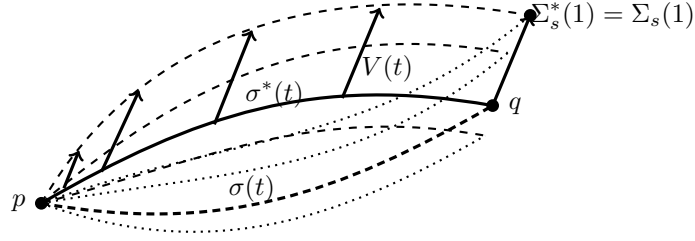
\begin{figure}
\centering
\begin{tikzpicture}
\node (p) at (-3.3,-1.5) {$p$};
\node (q) at (3.3,-.2) {$q$};
\node (S) at (4.6,1) {$\Sigma_s^*(1)=\Sigma_s(1)$};
\node (V) at (1.6,.3) {$V(t)$};
\draw[very thick]  (-3,-1.5) to [out=30, in=170] (3,-.2) 
node at (-.2,-1.3) {$\sigma(t)$};
\draw [very thick, densely dashed] (-3,-1.5) to [out=-10, in=-150] (3,-.2) 
node at (.1,-.05) {$\sigma^*(t)$};
\draw [black, thick] plot [mark=*, only marks]
coordinates {(-3,-1.5) (3,-.2) (3.5,1)};
\draw[thick,dashed]  (-3,-1.5) to [out=40, in=170] (3.2,.5);
\draw[thick,dotted]  (-3,-1.5) to [out=10, in=-140] (3.2,.5);
\draw[thick,dashed]  (-3,-1.5) to [out=60, in=170] (3.5,1);
\draw[thick,dotted]  (-3,-1.5) to [out=20, in=-140] (3.5,1);
\draw[thick,dashed]  (-3,-1.5) to [out=15, in=170] (2.9,-.6);
\draw[thick,dotted]  (-3,-1.5) to [out=-20, in=-150] (2.9,-.6);
\draw[->,very thick] (-2.7,-1.3) to (-2.5,-.8);
\draw[->,very thick] (-2.2,-1.08) to (-1.7,0);
\draw[->,very thick] (-.7,-.45) to (-0.2,.8);
\draw[->,very thick] (1,-.1) to (1.5,1.1);
\draw[->,very thick] (3,-.2) to (3.5,1);
\end{tikzpicture}
\caption{Given $p,\,q$ in a dually convex set $\U\subset\Ma$, we have a $\nabla^*$-geodesic $\sigma^*(t)$ (in solid line) and a $\nabla$-geodesic $\sigma(t)$ (in densely dashed line) such that $\sigma^*(0)=p=\sigma(0)$ and $\sigma^*(1)=q=\sigma(1)$. A $\nabla^*$-geodesic variation $\Sigma^*(s,t)$ (in dashed lines) of $\sigma^*$ with variation field $V(t)$ induces a $\nabla$-geodesic variation $\Sigma(s,t)$ (in dotted lines) of the $\nabla$-geodesic $\sigma$.}\label{Variations}
\end{figure}

\vspace{.2cm}

\noindent {At this point, we can compute the $\nabla$-pseudo-energy $\Lag$ at $\Sigma^*_s(t)$.} We then obtain 
$$
\Lag(s)=\int_0^1\langle\dot{\Sigma}^*_s(t),\paralleltransport_{s,t}\nihat_p(s)\rangle_{\Sigma^*_s(t)}\ \total t\ ,
$$
where $\paralleltransport_{s,t}:\tangent_p\Ma\rightarrow\tangent_{\Sigma^*(s,t)}\Ma$ is the $\nabla$-parallel transport along the curve $\Sigma^*_s(t)$ and $\nihat_p(s)=\exp_p^{-1}(\Sigma^*_s(1))$. Since $\Sigma^*_s(t)$ is the $\nabla^*$-geodesic from $p$ to $\Sigma_s^*(1)$,  we can write $\dot{\Sigma}_s^*(t)=\paralleltransport^*_{s,t}\dot{\Sigma}_s^*(0)$, where $\paralleltransport^*_{s,t}:\tangent_p\Ma\rightarrow\tangent_{\Sigma^*(s,t)}\Ma$ is the $\nabla^*$-parallel transport along the curve $\Sigma^*_s(t)$. Therefore, by recalling that the inner product is preserved under the parallel transport of pairs of vectors with respect to
a pair of dual connections, we get
\begin{eqnarray*}
\Lag(s)&=&\int_0^1 \langle \dot{\Sigma}^*_s(t),\paralleltransport_{s,t}\nihat_p(s)\rangle_{\Sigma^*_s(t)}\ \total t= \int_0^1 \langle \dot{\Sigma}^*_s(0),\nihat_p(s)\rangle_{p}\ \total t\\
&=& \langle\dot{\Sigma}^*_s(0),\nihat_p(s)\rangle_{p}  \ .
\end{eqnarray*}

Now, we may observe that $\dot{\Sigma}^*_s(0)=\tau(s)$. Moreover, since we know that $\nihat_p(s)=\exp^{-1}_p(\Sigma^*_s(1))$ and $\Sigma^*_s(1)=\stackrel{*}{\exp}_p(\tau(s))$, we can write
$$
\nihat_p(s)=\exp^{-1}_p\left(\stackrel{*}{\exp}_p(\tau(s))\right)\,.
$$
Hence, we obtain 
\begin{eqnarray}\label{Lag*}
\Lag(s)&=&\langle \dot{\Sigma}^*_s(0),\nihat_p(s)\rangle_p\nonumber\\
&=& \left\langle \tau(s),\exp_p^{-1}\left(\stackrel{^*}{\exp}_p(\tau(s))\right)\right\rangle_p = \kappa\,,
\end{eqnarray}
because $\tau(s)\in\hyper_p(\kappa)$. Finally,  we trivially have that
\begin{equation}\label{zero}
\left.\frac{\total}{\total s}\Lag(s)\right|_{s=0}=0\ .
\end{equation}

\vspace{.2cm}

\noindent Consider now the first variation of $\Lag(s)$,
$$
\frac{\total\Lag}{\total s}(s)= \left.{\langle  \partial_s\Sigma^*(s,t),\paralleltransport_{s,t}\nihat_p(s)\rangle_{\Sigma_s^*(t)}}\right|_{0}^1+\int_0^1 \langle  \dot{\Sigma}^*_s(t),\nabla_s\paralleltransport_{s,t}\nihat_p(s)\rangle_{\Sigma_s^*(t)}\,\total t\,.
$$
{Let us focus our investigation on the integrand $\langle  \dot{\Sigma}^*_s(t),\nabla_s\paralleltransport_{s,t}\nihat_p(s)\rangle_{\Sigma_s^*(t)}$. Define $\Xi(s,t):=\paralleltransport_{s,t}\nihat_p(s)$. This is a section on the variation $\Sigma^*(s,t)$. By exploiting the property that the inner product is invariant under the parallel transport with respect to dual connections, we then obtain
$$
\langle  \dot{\Sigma}^*_s(t),\nabla_s\paralleltransport_{s,t}\nihat_p(s)\rangle_{\Sigma_s^*(t)}=\langle \paralleltransport_{s,t}^{*-1} \dot{\Sigma}^*_s(t),\paralleltransport_{s,t}^{-1}\nabla_s\Xi(s,t)\rangle_{p}\,,
$$
where $\paralleltransport_{s,t}^{*-1},\,\paralleltransport_{s,t}^{-1}:\tangent_{\Sigma^*(s,t)}\Ma\rightarrow\tangent_p\Ma$ are the $\nabla^*$ and $\nabla$ parallel transports, respectively, along the curve $\Sigma^*_s(t)$. Recalling that $\Sigma_s^*(t)$ is a $\nabla^*$-geodesic, we can write $\dot{\Sigma}^*_s(t)=\paralleltransport^*_{s,t} \dot{\Sigma}^*_s(0)$ and then
\begin{equation}
\label{Equtile}
\langle  \dot{\Sigma}^*_s(t),\nabla_s\paralleltransport_{s,t}\nihat_p(s)\rangle_{\Sigma_s^*(t)}=\langle \dot{\Sigma}^*_s(0),\paralleltransport_{s,t}^{-1}\nabla_s\Xi(s,t)\rangle_{p}\,.
\end{equation}
Choose a basis $\{e_i(s,0)\}\subset\tangent_p\Ma$. Then by the definition of $\Xi(s,t)$ and writing $\nihat_p(s)=\sum_i p^i(s,0)\,e_i(s,0)$, we get $\Xi(s,t)=\sum_i p^i(s,0)\,\paralleltransport_{s,t}\,e_i(s,0)$. Now, for a given $t\in[0,1]$, we can choose a $\nabla$-parallel frame $\{e_i(s,t)\}$ along the curve $\Sigma^{*(t)}(s)$ such that $e_i(s,t)=\paralleltransport_{s,t}\,e_i(s,0)$ for all $i$. By applying the Leibniz rule for the connection $\nabla$, we obtain
$$
\nabla_s\Xi(s,t)=\sum_i(\partial_s p^i)(s,0)\,\paralleltransport_{s,t} e_i(s,0)\,,
$$
because the sections $e_i(s,t)$ are $\nabla$-parallel along $\Sigma^{*(t)}(s)$. Therefore, we get
$$
\paralleltransport_{s,t}^{-1}\nabla_s\Xi(s,t)=\sum_i(\partial_s p^i)(s,0)\,e_i(s,0)\,.
$$

For $\tilde{t}$ ranging in $[0,1]$, consider the $\nabla$-parallel transport along the curve $\Sigma_s(\tilde{t})$, namely $\widetilde{\paralleltransport}_{s,\tilde{t}}:\tangent_p\Ma\rightarrow\tangent_{\Sigma(s,\tilde{t})}\Ma$, where $\Sigma(s,t)$ is the $\nabla$-geodesic variation of $\sigma$ defined in Eq. \eqref{nabla-geod-var}. Hence, we can write
$$
\widetilde{\paralleltransport}_{s,\tilde{t}}\paralleltransport_{s,t}^{-1}\nabla_s\Xi(s,t)=\sum_i(\partial_s p^i)(s,0)\,\widetilde{\paralleltransport}_{s,\tilde{t}}\,e_i(s,0)\,.
$$
For a given $\tilde{t}\in[0,1]$, let us choose a $\nabla$-parallel frame $\{\tilde{e}_i(s,\tilde{t})\}$ along the curve $\Sigma^{(\tilde{t})}(s)$ such that $\tilde{e}_i(s,\tilde{t})=\widetilde{\paralleltransport}_{s,\tilde{t}}\,e_i(s,0)$. Hence, we have that
\begin{eqnarray*}
\nabla_s \widetilde{\paralleltransport}_{s,\tilde{t}}\nihat_p(s)&=&\nabla_s\left(\sum_i p^i(s,0)\,\widetilde{\paralleltransport}_{s,\tilde{t}}\,e_i(s,0)\right)\\
&=& \sum_i (\partial_s p^i)(s,0)\, \widetilde{\paralleltransport}_{s,\tilde{t}}\,e_i(s,0) \\
&=&\widetilde{\paralleltransport}_{s,\tilde{t}}\paralleltransport_{s,t}^{-1}\nabla_s\Xi(s,t)\,,
\end{eqnarray*}
where we applied the Leibniz rule for the $\nabla$-connection and used $\nabla_s\tilde{e}_i(s,\tilde{t})=0$, with $\nabla_s$ denoting the $\nabla$-covariant derivative along the curve $\Sigma^{(\tilde{t})}(s)$. By exploiting again the property that the inner product is invariant under the parallel transport with respect to dual connections, we then obtain 
\begin{eqnarray}
\langle \dot{\Sigma}^*_s(0),\paralleltransport_{s,t}^{-1}\nabla_s\Xi(s,t)\rangle_{p}&=&\left\langle \widetilde{\paralleltransport}^*_{s,\tilde{t}}\dot{\Sigma}^*_s(0),\widetilde{\paralleltransport}_{s,\tilde{t}}\paralleltransport_{s,t}^{-1}\nabla_s\Xi(s,t)\right\rangle_{\Sigma_s(\tilde{t})}\nonumber\\
\label{Equtile1}
&=& \left\langle \widetilde{\paralleltransport}^*_{s,\tilde{t}}\dot{\Sigma}^*_s(0),\nabla_s \widetilde{\paralleltransport}_{s,\tilde{t}}\nihat_p(s)\right\rangle_{\Sigma_s(\tilde{t})}\,,
\end{eqnarray} 
where $\widetilde{\paralleltransport}^*_{s,\tilde{t}}:\tangent_p\Ma\rightarrow\tangent_{\Sigma(s,\tilde{t})}\Ma$ is the $\nabla^*$-parallel transport along the curve $\Sigma_s(\tilde{t})$. Finally, by comparing Eq. \eqref{Equtile} and  Eq. \eqref{Equtile1} we get
\begin{equation}
\label{EqImportante}
\langle  \dot{\Sigma}^*_s(t),\nabla_s\paralleltransport_{s,t}\nihat_p(s)\rangle_{\Sigma_s^*(t)}=\left\langle \widetilde{\paralleltransport}^*_{s,\tilde{t}}\dot{\Sigma}^*_s(0),\nabla_s \widetilde{\paralleltransport}_{s,\tilde{t}}\nihat_p(s)\right\rangle_{\Sigma_s(\tilde{t})}\,.
\end{equation}
Now observe that $\dot{\Sigma}_s^*(0)$ is the velocity vector at $p$ of the $\nabla^*$-geodesic connecting $p$ and $\Sigma_s^*(1)=\Sigma_s(1)$. Then, in agreement with Eq. \eqref{*} we use the following notation,
 \begin{equation}
 \label{**}
 \dot{\Sigma}_s^*(0)= \stackrel{*}{\exp}_p^{-1}(\Sigma_s(1))=:\nihat^*_p(s)\equiv\nihat^*(p,\Sigma_s(1))\,.
 \end{equation}
Moreover, since $\nihat_p(s)=\dot{\Sigma}_s(0)$ and $\Sigma_s(\tilde{t})$ is a $\nabla$-geodesic, we can write $\widetilde{\paralleltransport}_{s,\tilde{t}}\nihat_p(s)=\dot{\Sigma}_s(\tilde{t})$. Therefore, Eq. \eqref{EqImportante} becomes
\begin{equation}
\label{Eqmain}
\langle  \dot{\Sigma}^*_s(t),\nabla_s\paralleltransport_{s,t}\nihat_p(s)\rangle_{\Sigma_s^*(t)}=\left\langle \widetilde{\paralleltransport}^*_{s,\tilde{t}}\nihat_p^*(s),\nabla_s\,\dot{\Sigma}_s(\tilde{t})\right\rangle_{\Sigma_s(\tilde{t})}\,.
\end{equation}
Then, we can use this expression to compute the integral in the first variation of the $\nabla$-pseudo-energy $\Lag$. In particular, by recalling that the $\nabla$-connection is torsion-free we have $\nabla_s(\partial_t\Sigma(s,t))=\nabla_t(\partial_s\Sigma(s,t))$. This allows us to perform the following computation
\begin{eqnarray}
\int_0^1 \left\langle  \dot{\Sigma}^*_s(t),\nabla_s\paralleltransport_{s,t}\nihat_p(s)\right\rangle_{\Sigma_s^*(t)}\total t &=& \int_0^1 \left\langle \widetilde{\paralleltransport}^*_{s,{t}}\nihat_p^*(s),\nabla_s\,\dot{\Sigma}_s({t})\right\rangle_{\Sigma_s({t})}\total t \nonumber\\
&=& \int_0^1 \left\langle \widetilde{\paralleltransport}^*_{s,{t}}\nihat_p^*(s),\nabla_t\,\partial_s\Sigma(s,t)\right\rangle_{\Sigma_s({t})}\total t\nonumber\\
&=&\int_0^1\frac{\total}{\total t}\left\langle \widetilde{\paralleltransport}^*_{s,{t}}\nihat_p^*(s),\partial_s\Sigma(s,t)\right\rangle_{\Sigma_s({t})}\total t \nonumber\\
&&-\int_0^1 \left\langle \nabla^*_t \widetilde{\paralleltransport}^*_{s,{t}}\nihat_p^*(s),\partial_s\Sigma(s,t)\right\rangle_{\Sigma_s({t})}\total t\nonumber\\
&=&\left.\left\langle \widetilde{\paralleltransport}^*_{s,{t}}\nihat_p^*(s),\partial_s\Sigma(s,t)\right\rangle_{\Sigma_s({t})}\right|_0^1\,,
\end{eqnarray}
where we exploited the property \eqref{dualconnection} and used $\nabla^*_t\widetilde{\paralleltransport}^*_{s,{t}}\nihat_p^*(s)=0$. By plugging the latter expression into the first variation of $\Lag(s)$ we then obtain
\begin{equation*}
 \frac{\total\Lag}{\total s}(s)=\left.\langle\partial_s\Sigma^*(s,t), \paralleltransport_{s,t}\nihat_p(s)\rangle_{\Sigma_s^*(t)}\right|_0^1+\left.\langle\partial_s\Sigma(s,t), \paralleltransport^*_{s,t}\nihat^*_p(s)\rangle_{\Sigma_s(t)}\right|_0^1\, .
 \end{equation*}
By setting $s=0$ we have that
\begin{align*}
&\left.\nihat_p(s)\right|_{s=0}=\exp_p^{-1}(q)=\nihat_p(q), \qquad \left.\nihat^*_p(s)\right|_{s=0}=\stackrel{*}{\exp}_p^{-1}(q)=\nihat^*_p(q)\\
&\left.\paralleltransport_{s,t}\nihat_p(s)\right|_{s=0}=\paralleltransport_{\sigma^*}\nihat_p(q)=\Pi_q(p),\quad \left.\widetilde{\paralleltransport}^*_{s,t}\nihat^*_p(s)\right|_{s=0}=\widetilde{\paralleltransport}^*_{\sigma}\nihat^*_p(q)=\Pi^*_q(p)\,,
\end{align*}
since $\Sigma_s^*(0,t)=\sigma^*(t)$ and $\Sigma_s(0,t)=\sigma(t)$. Hence, we get
\begin{equation}
\label{Sum}
 \frac{\total\Lag}{\total s}(0)=\langle\partial_s \Sigma^*(0,1),\Pi_q(p)\rangle_q
+\langle\partial_s\Sigma(0,1),\Pi_q^*(p)\rangle_q\,,
\end{equation}
because $\Sigma(s,0)=p=\Sigma^*(s,0)\,\forall\,s\in(-\varepsilon,\varepsilon)$ implies that 
$
\partial_s\Sigma^*(0,0)=O_p=\partial_s\Sigma(0,0)\ ,
$
where $O_p$ is the null vector in $\tangent_p\Ma$.
}

\vspace{.2cm}

In order to conclude the proof of  Theorem \ref{thmorthog}, let us {recall that $\Sigma(s,1)=\exp_p(I_p(\tau(s)))$ and $\tau(0)=\nihat_p^*(q)$. Then, we can carry out the following computation:} 
\begin{eqnarray}
\partial_s \Sigma(0,1)&=&\left(\total \exp_p\right)_{I_p(\nihat_p^*(q))}\left[\left.\partial_s I_p(\tau(s))\right|_{s=0}\right]\nonumber\\
&=&\left(\total \exp_p\right)_{I_p(\nihat_p^*(q))}\left[\left(\total I_p\right)_{\nihat_p^*(q)}(\tau^{\prime}(0))\right]\nonumber\\
&=& \left(\total \exp_p\right)_{I_p(\nihat_p^*(q))}\left[\left(\total \exp_p\right)^{-1}_{I_p(\nihat_p^*(q))}\left(\total \stackrel{*}{\exp_p}\right)_{\nihat^*_p(q)}(\tau^{\prime}(0))\right]\nonumber\\
&=&\left(\total \exp_p\right)_{I_p(\nihat_p^*(q))}\left[\left(\total \exp_p\right)^{-1}_{I_p(\nihat_p^*(q))}\left(\partial_s \Sigma^*(0,1)\right)\right]\nonumber\\
\label{samevariation}
&=&\partial_s \Sigma^*(0,1)\,,
\end{eqnarray}
{where we used $\total(\stackrel{*}{\exp}_p)_{\nihat_p^*(q)}\left(\tau^{\prime}(0)\right)=\partial_s\Sigma^*(0,1)$ which directly follows from\newline $\Sigma^*(s,1) = \stackrel{*}{\exp}_p(\tau(s))$.}

Finally, from Eq. \eqref{zero} {and plugging $\partial_s\Sigma(0,1)=\partial_s\Sigma^*(0,1)$ in Eq. \eqref{Sum},}  we obtain
\begin{equation}
0=\langle\partial_s\Sigma^*(0,1),\Pi_q(p)+\Pi_q^*(p)\rangle_q
\end{equation}
for an arbitrary tangent vector $\partial_s\Sigma^*(0,1)$ at $q$ to $\stackrel{*}{\exp}_p(\hyper_p(\kappa))=\hyperm_p(\kappa)$. This proves that the sum $\Pi_q(p)+\Pi^*_q(p)$ is orthogonal at $q$ to the hypersurface $\hyperm_p(\kappa)$ of constant pseudo-squared-distance $r_p(q)$. \hfill $\square$

\vspace{.5cm}

We are now ready to prove Theorem \ref{thmgradient} of Section \ref{Overwiev}. This claims that for $p,q\in\U$ in a dually convex set $\U\subset\Ma$, the pseudo-squared-distance $r(p,q)$ is the potential function of the sum $\Pi_q(p)+\Pi^*_q(p)$ of the vectors $\Pi_q(p)$ and $\Pi_q(p)$ of Eqs. \eqref{P} and \eqref{P*}.

\noindent{\bf \large Proof of  Theorem \ref{thmgradient}}. \, Consider $p\in\Ma$ and $\mathrm{U}_p\subset \Ma$ a dually convex neighborhood of $p$.  {Recall that} the pseudo-squared-distance $r_p(q)\equiv r(p,q)$ is defined for all $q\in\mathrm{U}_p$ by
$$
r_p(q)=\langle\exp_p^{-1}(q),\stackrel{*}{\exp}_p^{-1}(q)\rangle_p\ .
$$
In order to prove that
$$
\grad_q r_p=\Pi_q(p)+\Pi^*_q(p),
$$
{consider a variation of the end point $q$. This is given by $\Sigma^*(s,1)=\stackrel{*}{\exp}_p(\tau(s)),\,\forall\,s\in(-\varepsilon,\varepsilon)$. Here $\tau:(-\varepsilon,\varepsilon)\rightarrow\mathcal{E}_p$ is a curve within the neighborhood $\mathcal{E}_p\ni O_p\in\tangent_p\Ma$ such that $\exp_p,\,\stackrel{*}{\exp}_p:\mathcal{E}_p\rightarrow\Ma$ are diffeomorpsims on their images. Moreover, we require that $\tau(0)=\nihat_p^*(q)=\stackrel{*}{\exp}_p^{-1}(q)$. We may observe that a variation of the end point $q$ is $\Sigma(s,1)=\exp_p(I_p(\tau(s)))$, as well. Indeed, by recalling the definition
$$
I_p:\mathcal{E}_p\rightarrow\mathcal{E}_p,\quad X\mapsto \exp_p^{-1}\left(\stackrel{*}{\exp}_p(X)\right),\,\, \forall\,X\in\mathcal{E}_p\,,
$$
we have that $\Sigma^*(s,1)=\Sigma(s,1)$ for every $s\in(-\varepsilon,\varepsilon)$.} 

\vspace{.2cm}

From Eq. \eqref{Lagconst} we know  that the pseudo-squared-distance $r_p(q)$ is obtained by computing $\Lag$ at the $\nabla^*$-geodesic $\sigma^*$, namely
$
\Lag(\sigma^*)=r_p(q).
$
{Consider the end point variation $\Sigma^*(s,1)\equiv \Sigma_s^*(1)$. Then for every $s\in(-\varepsilon,\varepsilon)$ we have that
\begin{align*}
&r_p(\Sigma_s^*(1))=\left\langle\nihat_p(s),\nihat^*_p(s)\right\rangle_p\,,\quad \mbox{where} \\
&\nihat_p(s)=\exp_p^{-1}(\Sigma_s^*(1)),\quad \nihat^*_p(s)=\stackrel{*}{\exp}_p^{-1}(\Sigma_s^*(1))\,.
\end{align*}
By recalling the definition of the pseudo-energy $\Lag$, we can evaluate it on the $\nabla^*$-geodesic $\Sigma_s^*(t)$:
\begin{eqnarray*}
\Lag(s)\equiv\Lag(\Sigma_s^*)&=&\int_0^1\left\langle\dot{\Sigma}^*_s(t),\paralleltransport_{s,t}\nihat_p(s)\right\rangle_{\Sigma_s^*(s)}\total t\\
&=& \int_0^1\left\langle\paralleltransport^*_{s,t}\dot{\Sigma}^*_s(0),\paralleltransport_{s,t}\nihat_p(s)\right\rangle_{\Sigma_s^*(s)}\total t\\
&=&\int_0^1 \left\langle \dot{\Sigma}^*_s(0),\nihat_p(s)\right\rangle_{p}\total t\\
&=&\left\langle\nihat_p^*(s),\nihat_p(s)\right\rangle_p\,=\,r_p\left(\Sigma^*_s(1)\right)\,,
\end{eqnarray*}
where $\paralleltransport_{s,t},\paralleltransport^*_{s,t}:\tangent_p\Ma\rightarrow\tangent_{\Sigma^*(s,t)}\Ma$ are the $\nabla$ and $\nabla^*$ parallel transports along $\Sigma_s^*(t)$ and we used the property that the inner product is invariant under the parallel transport with respect to dual connections. Therefore, we can write
\begin{equation}
\label{diffLag}
\left.\frac{\total r_p(\Sigma_s^*(1))}{\total s}\right|_{s=0}=\left.\frac{\total \Lag(s)}{\total s}\right|_{s=0}\,.
\end{equation}
}
On the other hand, we can write 
$$
\left.\frac{\total r_p(\Sigma_s^*(1))}{\total s}\right|_{s=0}=\left(\total r_p\right)_q\left(\partial_s\Sigma^*(0,1)\right)=\left\langle\grad_q\,r_p,\partial_s\Sigma_s^*(0,1)\right\rangle_q\,.
$$
From Eq. \eqref{Sum} and Eq. \eqref{samevariation} we know that
$$
\left.\frac{\total \Lag(s)}{\total s}\right|_{s=0}=\left\langle\Pi_q(p)+\Pi^*_q(p),\partial_s\Sigma^*(0,1)\right\rangle_q\,.
$$
Hence, we obtain
\begin{equation}
\label{gradr}
\left\langle\grad_q\,r_p,\partial_s\Sigma_s^*(0,1)\right\rangle_q=\left\langle\Pi_q(p)+\Pi^*_q(p),\partial_s\Sigma^*(0,1)\right\rangle_q\,.
\end{equation}
In addition, from Theorem \ref{thmorthog} we know that $\Pi_q(p)+\Pi^*_q(p)$ is orthogonal at $q$ to hypersurfaces of constant $r_p(q)$. Therefore, Eq. \eqref{gradr} implies that
$$
\grad_q r_p=\Pi_q(p)+\Pi^*_q(p)\,.
$$
\hfill $\square $

\vspace{.5cm}

Owing to the dual structure of a statistical manifold $(\Ma,\metric,\nabla,\nabla^*)$, the theory  developed so far around $\Lag$, can be accomplished also in terms of a different pseudo-energy. Given $p,q\in\U$ in a dually convex set, we define the $\nabla^*$-pseudo-energy $\Lag^*$ by 
\begin{equation}
\label{pseudoenergy*}
\Lag^*(\gamma):=\int_0^1 \left\langle \dot{\gamma}(t),\paralleltransport^*_{t} \nihat^*_p(q)\right\rangle_{\gamma(t)}\ \total t,\quad \nihat_p^*(q)=\stackrel{*}{\exp}_p^{-1}(q)\,,
\end{equation}
where $\gamma$ is an arbitrary path such that $\gamma(0)=p$ and $\gamma(1)=q$, and $\paralleltransport^*_{t}:\tangent_p\Ma\rightarrow\tangent_{\gamma(t)}\Ma$ is the $\nabla^*$-parallel transport along $\gamma(t)$. In this case the pseudo-squared-distance $r_p(q)$ is obtained by computing $\Lag^*$ at the $\nabla$-geodesic $\sigma(t)$ connecting $p$ and $q$ {as claimed by the following Proposition.} 
\begin{pro}\label{Energy*Prop}
Let $p,q\in\U$ in a dually convex set and let $\sigma:[0,1]\rightarrow\U$ be a $\nabla$-geodesic such that $\sigma(0)=p$ and $\sigma(1)=q$. Then
\begin{equation}
\label{Lagconst*}
\Lag^*(\sigma)=\langle \exp_p^{-1}(q),\stackrel{*}{\exp}_p^{-1}(q)\rangle_q\equiv r_p(q)\, .
\end{equation}
\end{pro}
\noindent{\bf Proof.}\ The proof is identical to the proof of Proposition \ref{EnergyProp}, just replacing $\nabla$ by $\nabla^*$.

\vspace{.2cm}

The first variation of $\Lag^*$ is stated in the following Proposition and {can be proved in the same way as the Proposition \ref{FirstProp}.}

\begin{pro}
For $p,q$ in a dually convex set $\U\subset\Ma$, consider an arbitrary path  $\gamma:[0,1]\rightarrow\U$ {such that $\gamma(0)=p$ and $\gamma(1)=q$} and a variation of $\gamma$, $\Sigma:(-\varepsilon,\varepsilon)\times [0,1]\rightarrow\Ma$. Let $V\in\Tau(\gamma)$ be the variation vector field of $\Sigma$. Finally, define the functional $\Lag^*(s):=\Lag^*(\Sigma_s)$. Then we have
\begin{equation}
\label{first*}
\frac{\total\Lag^*}{\total s}(0)=\left. \langle V(t),\paralleltransport^*_{t}\nihat^*_p(q)\rangle_{\gamma(t)}\right|_0^1+\int_0^1\langle\dot{\gamma}(t),\nabla^*_V\paralleltransport^*_{t} \nihat^*_p(q)\rangle_{\gamma(t)} dt,
\end{equation}
where $\paralleltransport^*_{t}:\tangent_p\Ma\rightarrow\tangent_{\gamma(t)}\Ma$ denotes the $\nabla^*$-parallel transport along $\gamma(t)$.
\end{pro}

\vspace{.2cm}

{Hence, by means of the first $\nabla$-geodesic variation of the $\nabla^*$-pseudo-energy $\Lag^*$, we can prove Theorem \ref{thmorthog} by following the same methods carried out around the $\nabla$-pseudo-energy $\Lag$.}

\section{Proof of Theorem 5}\label{TheoremPi&Pi*orthogonal}

In this section we aim to prove Theorem \ref{Pi&Pi*orthogonal} of Section \ref{Overwiev}. We will carry out its proof in several steps. First of all, we  show the complementarity of the canonical divergences $\Div$ and $\Div^*$ with the Phi-functions $\varphi$ and $\varphi^*$. This is obtained by developing the potential theoretic-property of the pseudo-squared-distance $r(p,q)$ claimed by Theorem \ref{thmgradient}. 

Thus, for  $p,q$ in a dually convex set $\U\subset\Ma$ consider an arbitrary path $\gamma:[0,1]\rightarrow\mathrm{U}$ such that $\gamma(0)=p$ and $\gamma(1)=q$. From the definition of the vectors $\Pi$ and $\Pi^*$ given by Eq. \eqref{P} and Eq. \eqref{P*}, respectively, we can supply two vector fields along $\gamma$ in the following way. Since $\gamma(t)$ is in $\mathrm{U}$  for every $t\in[0,1]$, we can find a $\nabla$-geodesic $\sigma_t(s)\,(0\leq s\leq 1)$ and a $\nabla^*$-geodesic $\sigma_t^*(s)\,(0\leq s\leq 1)$ such that $\sigma_t(0)=p=\sigma_t^*(0)$ and $\sigma_t(1)=\gamma(t)=\sigma_t^*(1)$. Hence, according to Eq. \eqref{P} and Eq. \eqref{P*}, we can write
\begin{align}
\Pi_{t}(p)= \paralleltransport_{\sigma^*_t} \nihat_p(\gamma(t)),&\qquad \nihat_p(\gamma(t))=\exp_p^{-1}(\gamma(t))\label{vectorfieldPi} \\
\Pi_t^*(p)=\paralleltransport^*_{\sigma_t}\nihat_p^*(\gamma(t)),&\qquad \nihat_p^*(\gamma(t))=\stackrel{*}{\exp}_p^{-1}(\gamma(t)) \label{vectorfieldPI*}\,,
\end{align}
where $\paralleltransport_{\sigma^*_t}:\tangent_p\Ma\rightarrow\tangent_{\gamma(t)}\Ma$ is the $\nabla$-parallel transport along the $\nabla^*$-geodesic $\sigma^*_t(s)$ and $\paralleltransport^*_{\sigma_t}:\tangent_p\Ma\rightarrow\tangent_{\gamma(t)}\Ma$ is the $\nabla^*$-parallel transport along the $\nabla$-geodesic $\sigma_t(s)$. Then, consider the pseudo-squared-distance between $p$ and $\gamma(t)$, namely 
$$
r_p(\gamma(t))=\langle\nihat_p(\gamma(t)),\nihat^*_p(\gamma(t))\rangle_p\,,
$$
where $\nihat_p(\gamma(t)),\,\nihat^*_p(\gamma(t))$ are given in Eq. \eqref{vectorfieldPi} and   Eq. \eqref{vectorfieldPI*}, respectively. Hence, by applying the claim stated in Theorem \ref{thmgradient} we can write
\begin{equation}
\label{gradrt}
\grad_{\gamma(t)}\,r_p=\Pi_t(p)+\Pi_t^*(p)\,.
\end{equation}
Furthermore, we have that
\begin{equation}
\label{differential}
(\total r_p)_{\gamma(t)}(X_t) =\langle \grad_{\gamma(t)} r_p,X_t\rangle_{\gamma(t)},\quad \forall\ X_t\in\tangent_{\gamma(t)}\Ma\ .
\end{equation}
Finally, we can compose the inner product of the curve velocity $\dot{\gamma}(t)$ with the vector field $\Pi_t(p)+\Pi^*_t(p)$ and perform the following computation by means of Eq. \eqref{gradrt} and Eq. \eqref{differential}:
\begin{align}
\int_0^1\ \langle \Pi_t(p)+\Pi_t^*(p),\dot{\gamma}(t)\rangle_{\gamma(t)}\ dt& = \int_0^1\ \langle \grad_{\gamma(t)}\ r_p,\dot{\gamma}(t)\rangle_{\gamma(t)}\ \total t\nonumber\\
&=   \int_0^1\, \left(\total\, r_p\right)_{\gamma(t)}\left(\dot{\gamma}(t)\right)\ \total t\nonumber\\
&=  \int_0^1\ \frac{\total \ r_p\circ \gamma}{\total t}(t)\,\total t\nonumber\\
&= r_p(\gamma(1))-  r_p(\gamma(0))\nonumber\\
&= r_p(q),\label{gammapath}
\end{align}
where, obviously, we have $r_p(p)\equiv 0$. {This proves that} the sum
\begin{equation}
\label{indipendence}
\int_0^1\ \langle \Pi_t(p),\dot{\gamma}(t)\rangle_{\gamma(t)}\,\total t+\int_0^1\  \langle \Pi_t^*(p),\dot{\gamma}(t)\rangle_{\gamma(t)}\,\total t= r_p(q)
\end{equation}
is independent of the particular path from $p$ to $q$.

\vspace{.2cm}

From Eq. \eqref{indipendence}, the definition of the canonical divergence $\Div(p,q)$ has been obtained by considering the $\nabla$-geodesic $\sigma$, such that $\sigma(0)=p$ and $\sigma(1)=q$, in the first integral of the left hand side instead of the arbitrary path $\gamma$:
\begin{equation}
\label{divergence}
\Div(p,q):= \int_0^1 \langle \Pi_t(p),\dot{\sigma}(t)\rangle_{\sigma(t)}\ \total t\,,\quad \Pi_t(p)=\paralleltransport_{\sigma^*_t}\,{\exp}_p^{-1}(\sigma(t))\,,
\end{equation}
where $\paralleltransport_{\sigma^*_t}:\tangent_p\Ma\rightarrow\tangent_{\sigma(t)}\Ma$ denotes the $\nabla$-parallel transport along the $\nabla^*$-geodesic $\sigma^*_t(s)\,(0\leq s\leq 1)$ such that $\sigma_t^*(0)=p$ and $\sigma^*_t(1)=\sigma(t)$. On the contrary, the Phi-function $\varphi^*(p,q)$ has been defined by considering  the $\nabla$-geodesic $\sigma$ in the second integral of the left hand side of \eqref{indipendence},
\begin{equation}
\label{phifunction*}
\varphi^*(p,q)=\int_0^1\ \langle\Pi^*_t(p),\dot{\sigma}(t)\rangle_{\sigma(t)}\ \total t\,,\quad \Pi^*_t(p)=\paralleltransport^*_{\sigma_t}\,\stackrel{*}{\exp}_p^{-1}(\sigma(t))\,,
\end{equation} 
where $\paralleltransport^*_{\sigma_t}:\tangent_p\Ma\rightarrow\tangent_{\sigma(t)}\Ma$ is the $\nabla^*$-parallel transport along the $\nabla$-geodesic $\sigma_t(s)\,(0\leq s\leq 1)$ such that $\sigma_t(0)=p$ and $\sigma_t(1)=\sigma(t)$.
Then, since the value of $r(p,q)$ is independent of the particular path from $p$ to $q$, we can see that these two functions are complementary in the sense that the following relation holds:
\begin{equation}
\label{r&Div&Phi*}
r(p,q)\,=\,\Div(p,q)+\varphi^*(p,q)\,.
\end{equation}

On the contrary, if we plug the $\nabla^*$-geodesic $\sigma^*$, such that $\sigma^*(0)=p$ and $\sigma^*(1)=q$, into the left hand side of Eq. \eqref{indipendence}, we get from the second integral the dual canonical divergence $\Div^*(p,q)$,
\begin{equation}
\label{divergence*}
\Div^*(p,q):= \int_0^1 \langle \Pi^*_t(p),\dot{\sigma}^*(t)\rangle_{\sigma^*(t)}\ \total t\,, \quad
\Pi^*_t(p)=\paralleltransport^*_{\sigma_t}\,\stackrel{*}{\exp}_p^{-1}(\sigma^*(t))\,,
\end{equation}
where $\paralleltransport^*_{\sigma_t}:\tangent_p\Ma\rightarrow\tangent_{\sigma^*(t)}\Ma$ is the $\nabla^*$-parallel transport along the $\nabla$-geodesic $\sigma_t(s)\,(0\leq s\leq 1)$ such that $\sigma_t(0)=p$ and $\sigma_t(1)=\sigma^*(t)$. Moreover, from the first integral of the left hand side of \eqref{indipendence} we obtain the Phi-function $\varphi(p,q)$:
\begin{equation}
\label{phifunction}
\varphi(p,q)=\int_0^1\ \langle\Pi_t(p),\dot{\sigma}^*(t)\rangle_{\sigma^*(t)}\ \total t\,,\quad \Pi_t(p)=\paralleltransport_{\sigma_t^*}\,\exp_p^{-1}(\sigma^*(t))\,,
\end{equation}
where $\paralleltransport_{\sigma_t^*}:\tangent_p\Ma\rightarrow\tangent_{\sigma^*(t)}\Ma$ is the $\nabla$-parallel transport along the $\nabla^*$-geodesic $\sigma_t^*(s)\,(0\leq s\leq 1)$ such that $\sigma_t^*(0)=p$ and $\sigma_t^*(1)=\sigma^*(t)$. Also in this case the functions $\Div^*(p,q)$ and $\varphi(p,q)$ play a complementary role in the sense that the pseudo-squared-distance is given by
\begin{equation}
\label{r&Div*&Phi}
r_p(q)\,=\, \Div^*(p,q)+\varphi(p,q)\,.
\end{equation}

\vspace{.2cm}

In order to employ  Eqs. \eqref{r&Div&Phi*} and \eqref{r&Div*&Phi} for carrying out the proof of Theorem \ref{Pi&Pi*orthogonal}, we need two more ingredients. Firstly, we will prove Theorem \ref{localdecompositionPiI} which naturally characterizes $\varphi$ and $\varphi^*$ by the unique decompositions of $\Pi$ and $\Pi^*$ in terms of gradient vector fields and vector fields that are orthogonal to $\nabla^*$-geodesics and $\nabla$-geodesics. Later, we will show the consistency of  canonical divergences and Phi-functions with the dualistic structure of a statistical manifold $(\Ma,\metric,\nabla,\nabla^*)$. This latter part is relevant for applying the theory of minimum contrast geometry by Eguchi \cite{eguchi1992}, which in turn will allow us to complete the proof of Theorem \ref{Pi&Pi*orthogonal} by getting that $\grad\,\Div$ and $\grad\,\varphi$ are parallel to $\nabla$-geodesics whereas $\grad\,\Div^*$ and $\grad\,\varphi^*$ are parallel to $\nabla^*$-geodesics.

\subsection{Proof of Theorem 6}\label{Theorem 6}

Without loss of generality we shall focus on $\varphi(p,q)$ as we can obtain the theory for $\varphi^*$ in the same way just by interchanging the role of $\nabla$ and $\nabla^*$. The next result supplies a nice representation of $\varphi(p,q)$ which will be very useful for proving Theorem \ref{localdecompositionPiI} of  Section \ref{Overwiev} as well as the consistency of $\varphi(p,q)$ with the dual structure of a statistical manifold $(\Ma,\metric,\nabla,\nabla^*)$.

\begin{lemma}\label{Lemma}
Let $p,q\in\U$ in a dually convex set  and $\sigma^*(t)\,(0\leq t\leq 1)$ be the $\nabla^*$-geodesic such that $\sigma^*(0)=p$ and $\sigma^*(1)=q$. Then, we have
\begin{equation}
\label{phi}
\varphi(p,q)=\int_0^1\ \langle I_p(t\nihat^*_p(q)),\nihat_p^*(q)\rangle_p\ \total t,
\end{equation}
where $I_p(t \nihat_p^*(q))=\exp_p^{-1}\left(\stackrel{*}{\exp}_p(t\nihat^*_p(q))\right)$ and $\nihat_p^*(q)=\stackrel{*}{\exp}^{-1}_p(q)=\dot{\sigma}^*(0)$.
\end{lemma}
\noindent {\bf Proof}.\ Consider the $\nabla^*$-geodesic $\sigma^*_t(s)$ that connects $p$ {with} $\sigma^*(t)$, namely $\sigma_t^*(0)=p$ and $\sigma_t^*(1)=\sigma^*(t)$. Then, {for any $t\in[0,1]$ we can write} $\sigma_t^*(s)=\sigma^*(s t)\ (s\in[0,1])$. A classical result in Riemannian geometry tells us that \cite{Lee97},
$$
\dot{\sigma}_t^*(1)=t\ \dot{\sigma}^*(t)\ .
$$
By substituting this expression into Eq. \eqref{phifunction} we immediately obtain that
\begin{equation}
\label{divgrad}
\varphi(p,q)=\int_0^1 \frac{1}{t}\langle \Pi_t(p), \dot{\sigma}_t^*(1)\rangle_{\sigma^*(t)}\ \total t\ .
\end{equation}
Recall that
$$
\Pi_t(p)=\paralleltransport_{\sigma^*_t}\nihat_p(\sigma^*(t)),\quad \nihat_p(\sigma^*(t))=\exp_p^{-1}(\sigma^*(t)), \qquad \dot{\sigma}_t^*(1)=\paralleltransport^*_{\sigma^*_t}\dot{\sigma}_t^*(0)\ .
$$
Hence, thanks to the invariance of the inner product under the combined action of $\paralleltransport$ and $\paralleltransport^*$ we get
$$
\varphi(p,q)=\int_0^1 \frac{1}{t}\langle \nihat_p(\sigma^*(t)), \dot{\sigma}_t^*(0)\rangle_{p}\ \total t\ .
$$
Now, since $\sigma^*$ is the $\nabla^*$-geodesic from $p$ to $q$, we can write\newline
$\sigma^*(t)=\stackrel{*}{\exp}_p(t\,\nihat_p^*(q))$. Then we may observe that
$$
\nihat_p(\sigma^*(t))=\exp_p^{-1}\left(\stackrel{*}{\exp}_p(t\ \nihat_p^*(q))\right),\qquad \dot{\sigma}_t^*(0)=t\ \nihat_p^*(q)
$$
because $\sigma_t^*$ is a re-parametrization of $\sigma^*$. Therefore, we obtain
\begin{eqnarray*}
\varphi(p,q)&=&\int_0^1 \frac{1}{t}\left\langle \exp_p^{-1}\left(\stackrel{*}{\exp}_p(t\ \nihat_p^*(q))\right), t\ \nihat_p^*(q)\right\rangle_{p}\ \total t\\
&=& \int_0^1 \left\langle I_p(t\ \nihat^*_p(t)), \nihat_p^*(q)\right\rangle_{p}\ \total t,
\end{eqnarray*}
by recalling definition \eqref{Iso} of the map $I_p$. This proves claim \eqref{phi}. \hfill $\square$

\vspace{.2cm}
We are now in the position to prove Theorem \ref{localdecompositionPiI} of Section \ref{Overwiev}, which claims that the Phi-functions $\varphi(p,q)$ and $\varphi^*(p,q)$ are naturally characterized by the local decomposition of the vectors $\Pi$ and $\Pi^*$, respectively, in terms of gradient vectors. Without loss of generality, we shall prove only Eq. \eqref{locdecPiI} as the decomposition \eqref{locdecPi*I} can be obtained in the same way.

\noindent{\bf \large Proof of Theorem \ref{localdecompositionPiI}.}
By the definition \eqref{Iso} of the map $I_p$, we can rewrite the representation \eqref{phi} of $\varphi(p,q)$ in the following way,
\begin{equation}
\label{representationPhi}
\varphi_p(q)\equiv\left(\varphi_p\circ\stackrel{*}{\exp}_p\right)(\nihat_p^*(q))=\int_0^1 \frac{1}{t}\left\langle\exp_p^{-1}\left(\stackrel{*}{\exp}_p(t \nihat_p^*(q))\right),t \nihat_p^*(q)\right\rangle_p\ \total t\ ,
\end{equation}
where we have written $q=\stackrel{*}{\exp}_p(\nihat_p^*(q))$. 

For $p\in\Ma$ we can consider normal coordinates $\{\xi^i\}$ with respect to the affine connection $\nabla^*$ \cite{Katanaev2018} in a dually convex neighborhood $\mathrm{U}_p$ of $p$. Therefore, the component $\metric_{ij}$ of the metric tensor $\metric$ at $p$ can be written as
$$
\metric_{ij}(p)=\delta_i^j\,,
$$
where $\delta_i^j$ denotes the Kronecker's delta function, i.e. $\delta_i^j=1$ if $i=j$ and $\delta_i^j=0$ otherwise.

By identifying $\tangent_{\nihat_p^*(q)}\left(\tangent_p\Ma\right)$ with $\tangent_p\Ma$ in the canonical way, we can write 
$
\nihat_p^*(q)=\sum_i \xi^i\partial_i\ .
$
Consider now the following function, 
$$
\Xi(t\ \nihat_p^*(q)):=\left\langle \exp_p^{-1}\left(\stackrel{*}{\exp}_p(t\ \nihat_p^*(q))\right),t\ \nihat_p^*(q)\right\rangle_p\ .
$$
By Eq. \eqref{representationPhi}, we have that 
$$
\left(\varphi_p\circ\stackrel{*}{\exp}_p\right)\left(\nihat_p^*(q)\right)=\int_0^1\frac{1}{t}\ \Xi(t\ \nihat_p^*(q))\ \total t\ .
$$
Moreover, by means of normal coordinates $\{\xi^i\}$ we obtain
$$
\left\langle\nihat_p^*(q),\grad_{\nihat_p^*(q)}\left(\varphi\circ\stackrel{*}{\exp}_p\right)\right\rangle_p=\left\langle\xi^i\partial_i,\frac{\partial\left(\varphi\circ\stackrel{*}{\exp}_p\right)}{\partial\xi^j}\partial_j\right\rangle_p=\sum_i\xi^i\frac{\partial\left(\varphi\circ\stackrel{*}{\exp}_p\right)}{\partial\xi^i}\ .
$$
This yields the following computation,{
\begin{eqnarray}
\left\langle \nihat_p^*(q),\grad_{\nihat_p^*(q)} \left(\varphi_p\circ\stackrel{*}{\exp_p}\right)\right\rangle_p&=& \sum_i \xi^i \frac{\partial}{\partial\xi^i}\ \int_0^1 \frac{1}{t}\ \Xi(t\, \nihat_p^*(q))\ \total t \nonumber\\
&=&\sum_i \xi^i\int_0^1\,\frac{1}{t}\,\left(\frac{\partial}{\partial\xi^i}\Xi(t\, \nihat_p^*(q))\right)\, \total t \nonumber \\
&=& \sum_i \xi^i \int_0^1\,\frac{1}{t} t \left(\frac{\partial}{\partial \xi^i}\nihat_p^*(q)\right) \frac{\partial\Xi}{\partial \xi^i}(t\, \nihat_p^*(q))\ \total t \nonumber \\
&=&\int_0^1\,\sum_i (\xi^i\partial_i)\frac{\partial\Xi}{\partial \xi^i}(t\, \nihat_p^*(q))\ \total t \nonumber \\
&=& \int_0^1 \frac{\total}{\total t}\Xi(t\ \nihat_p^*(q))\ \total t=\Xi(1\ \nihat_p^*(q)),\label{relationThm3}
\end{eqnarray}
where we used $\nihat_p^*(q)=\sum_i\xi^i\partial_i$ and  $\Xi(0)=\langle \exp_p(p),0\rangle_p=0$. Recalling that
\begin{equation}
\label{Xi}
\Xi(\nihat_p^*(q)) = \langle \nihat_p(q),\nihat_p^*(q)\rangle_p\ ,
\end{equation}
we then get}
\begin{equation}
\label{gradatp}
\grad_{\nihat_p^*(q)} \left(\varphi_p\circ\stackrel{*}{\exp_p}\right)=\nihat_p(q)+V_p,
\end{equation}
where $\langle V_p,\nihat_p^*(q)\rangle_p=0$. 
{Exploiting the property that the inner product is invariant under the parallel transport with respect to dual connections, we can write
$$
\langle\Pi_q(p),\dot{\sigma}^*(1)\rangle_q = \langle\paralleltransport_{\sigma^*}\nihat_p(q),\paralleltransport^*_{\sigma^*}\nihat^*_p(q)\rangle_q = \langle \nihat_p(q),\nihat_p^*(q)\rangle_p\,,
$$
which follows from the definition \eqref{P} of $\Pi_q(p)$ and because $\sigma^*$ is a $\nabla^*$-geodesic. In addition, from Eq. \eqref{relationThm3}, Eq. \eqref{Xi}  we have that
\begin{eqnarray*}
\langle\Pi_q(p),\dot{\sigma}^*(1)\rangle_q &=&\langle \nihat_p(q),\nihat_p^*(q)\rangle_p\\
&=& \left\langle \nihat_p^*(q),\grad_{\nihat_p^*(q)} \left(\varphi\circ\stackrel{*}{\exp_p}\right)\right\rangle_p\\
&=&\total( \varphi_p\circ\stackrel{*}{\exp}_p)_{\nihat_p^*(q)}(\nihat_p^*(q))\,.
\end{eqnarray*}
}
At this point, we can observe that $(\total\varphi_p)_q=\total\left(\varphi_p\circ \stackrel{*}{\exp}_p\right)_{\nihat_p^*(q)}$. {Therefore, by the chain rule of the differential and from the well-known relation \newline
 $(\total \stackrel{*}{\exp}_p)_{\nihat_p^*(q)}(\nihat_p^*(q))=\dot{\sigma}^*(1)$, we obtain} 
\begin{eqnarray*}
\langle\Pi_q(p),\dot{\sigma}^*(1)\rangle_q &=& \total( \varphi_p\circ\stackrel{*}{\exp}_p)_{\nihat_p^*(q)}(\nihat_p^*(q))\\ 
&=&  (\total \varphi_p)_q\left((\total \stackrel{*}{\exp}_p)_{\nihat_p^*(q)}(\nihat_p^*(q))\right)\\
 &=& (\total \varphi_p)_q\left(\dot{\sigma}^*(1)\right) = \langle \grad_q \varphi_p,\dot{\sigma}^*(1)\rangle_q\ .
\end{eqnarray*}
{This proves that}
$$
\Pi_q(p)= \grad_q \varphi_p +V_q
$$
with $V_q\in\tangent_q\Ma$ uniquely defined by $\Pi_q(p)-\grad_q \varphi_p$ and $\langle V_q,\dot{\sigma}^*(1)\rangle_q=0$.

\vspace{.2cm}

In order to prove that the decomposition \eqref{locdecPiI} is unique, suppose that there exists another decomposition of $\Pi_q(p)$ satisfying the conditions of Theorem \ref{localdecompositionPiI}, i.e. $\Pi_q(p)=\grad_q\tilde{\varphi}_p(q)+\tilde{V}_q$ with $\langle\tilde{V}_q,\dot{\sigma}^*(1)\rangle_q=0$. In addition, let us assume that $\tilde{\varphi}_p(p)=0$. We have then,
$$
0=\grad_q\left(\varphi_p-\tilde{\varphi}_p\right)+V_q-\tilde{V}_q \quad \mbox{and}\quad \langle \dot{\sigma}^*(1),V_q-\tilde{V}_q\rangle_q=0 .
$$
It is evident that 
\begin{equation}
\label{Homogeneous}
\left\langle\dot{\sigma}^*(1),\grad_q\left(\varphi_p-\tilde{\varphi}_p\right)\right\rangle_q=0\ .
\end{equation}

For $q$ in a $\nabla^*$-normal neighborhood of $p$ we have the diffeomorphism
$$
\stackrel{*}{\exp}_p\left(\xi_1(q)\frac{\partial}{\partial \xi^1}+\ldots+\xi_n(q)\frac{\partial}{\partial \xi^n}\right)\rightarrow(\xi^1(q),\ldots,\xi^n(q))\,,
$$
where $(\xi^1(q),\ldots,\xi^n(q))\in\RR^n$. Recalling that $\grad_q\,f=\metric^{ij}(q)\frac{\partial\,f}{\partial \xi^i}\partial_i$ for any smooth function $f$, we may notice that 
 \eqref{Homogeneous} is a homogeneous first-order linear equation in partial derivatives which has only the constant solution $\varphi_p-\tilde{\varphi}_p=\mbox{const.}$,
$$
\sum_i \frac{\partial\,(\varphi_p-\widetilde{\varphi}_p)(\xi^i)}{\partial\xi^i}\,\total\xi^i=0
$$
Hence,  due to the assumption $\varphi_p(\xi(p))=0=\tilde{\varphi}_p(\xi(p))$ we obtain that $\varphi_p(q)=\tilde{\varphi}_p(q)$ and $V_q=\tilde{V}_q$. \hfill $\square$

\begin{remark}
\label{presnov}
Methods in the proof of Theorem \ref{localdecompositionPiI} are inspired by \cite{PRESNOV}, where the author presented a decomposition for vector fields on a Riemannian manifold of non-positive curvature with
application to non-linear mechanics and irreversible thermodynamics.
\end{remark}

\vspace{.3cm}

\subsection{Consistency Theorem}\label{ConsistencySec}

In this section, we aim to prove Theorem \ref{positivityI}.

In order to show that $\varphi(p,q)$ {is a divergence function on a statistical manifold $\Ma$,} we have to prove its consistency with the dual structure $(\metric,\nabla,\nabla^*)$ of $\Ma$. This means that in a neighborhood of the diagonal set of $\Ma\times\Ma$ we need to verify that Eqs. \eqref{metricfromdiv} and \eqref{dualfromDiv} are satisfied.

Let us now assume that $p$ and $q$ are close to each other, that is
\begin{equation}
\label{small}
z^i=\xi^i_q-\xi^i_p
\end{equation}
is small. Here $\{\bxi_p\}$ and $\{\bxi_q\}$ are local coordinates at $p$ and $q$, respectively. Then, {we start by providing a Taylor expansion of $\varphi(p,q)$ up to $O\left(\|\bz\|^4\right)$ .}
\begin{pro}\label{proptaylor}
Consider $\|\mathbf{z}\|=\|\bxi_q-\bxi_p\|$ small enough. Then, the function $\varphi(p,q)$ is expanded up to $O\left(\|\mathbf{z}\|^4\right)$ as follows
\begin{equation}
\label{Taylor}
\varphi(p,q)=\frac{1}{2}\ \metric_{ij}(p)\ z^i z^j+\frac{1}{6}\Lambda_{ijk}(p)\ z^i z^j z^k+ O\left(\|\mathbf{z}\|^4\right),
\end{equation}
where
\begin{equation}
\label{Tensor}
\Lambda_{ijk}(p)=2\Gamma_{ijk}^*(p)+\Gamma_{ijk}(p)\ .
\end{equation}
\end{pro}
\noindent {\bf Proof}.\ Let us consider the representation \eqref{phi} of the function $\varphi(p,q)$. Then, recall that
$I_p(t\ \nihat_p^*(t))$ is nothing but the velocity vector at $p$ of the $\nabla$-geodesic $\sigma_t(s)\,(0\leq s\leq 1)$ that connects $p$ with $\sigma^*(t)$. On the other hand, $\nihat^*_p(q)$ is the velocity vector at $p$ of the $\nabla^*$-geodesic $\sigma^*$ from $p$ to $q$. Therefore, we need to Taylor expand up to $O\left(\|\bz\|^4\right)$ with respect to the local coordinate $\{\bxi\}$ the following expression
\begin{equation}
\label{Taylor1}
\int_0^1\metric_{ij}(p)\
\dot{\sigma}_t(0)^i\
\dot{\sigma}^*(0)^j\ \total t\,,
\end{equation}
where the Einstein notation is adopted.

The local coordinates $\bxi(t)$ of the $\nabla$-geodesic $\sigma^*(t)$ in Taylor series are given by
\begin{equation}
\label{geod}
\xi^j(t)=\xi^j_p+t z^j+ \frac{t}{2}(1-t){\Gamma^*}_{\mu\nu}^j(p)z^{\mu}z^{\nu}+ O\left(\|\mathbf{z}\|^3\right),
\end{equation}
where the summation over $\mu$ and $\nu$ is understood. Then we obtain,
\begin{equation}
\label{geodderiv}
\frac{\total}{\total t}\sigma^*(0)^j=z^j+\frac{1}{2}{\Gamma^*}_{\mu\nu}^j(p)z^{\mu}z^{\nu}+O\left(\|\mathbf{z}\|^3\right).
\end{equation}

In addition we have that
\begin{equation*}
\frac{\total }{\total s}\sigma_t(0)^i=\sigma^*(t)^i-\xi_p^i+\frac{1}{2}\Gamma_{\mu\nu}^i(p)(\sigma^*(t)^{\mu}-\xi_p^{\mu})(\sigma^*(t)^{\nu}-\xi_p^{\nu})+O\left(\|\sigma^*(t)-\xi_p\|^3\right)\ .
\end{equation*}
Now, as 
$$
\sigma^*(t)^j-\xi_p^j=t z^j+\frac{t}{2}{\Gamma^*}_{\mu\nu}^j(p) z^{\mu}z^{\nu}-\frac{t^2}{2}{\Gamma^*}_{\mu\nu}^j(p)
z^{\mu}z^{\nu}+O\left(\|\mathbf{z}\|^3\right)
$$
we arrive at
\begin{equation}\label{geodt}
\frac{\total }{\total s}\sigma_t(0)^i=t z^i+\frac{t}{2}{\Gamma^*}^i_{\mu\nu}z^{\mu}z^{\nu}+\frac{t^2}{2}(\Gamma_{\mu\nu}^i(p)-{\Gamma^*}^i_{\mu\nu}(p))z^{\mu}z^{\nu}+O\left(\|\mathbf{z}\|^3\right)\ .
\end{equation}

At this point, we can plug Eq. \eqref{geodderiv} and Eq. \eqref{geodt} into Eq. \eqref{Taylor1}
\begin{eqnarray*}
&&\int_0^1\metric_{ij}(p)\
\dot{\sigma}_t(0)^i\
\dot{\sigma}^*(0)^j\ \total t =\\
&&=\int_0^1 g_{ij}(p)\left[t z^i+\frac{t}{2}{\Gamma^*}^i_{\mu\nu}(p)z^{\mu}z^{\nu}+\frac{t^2}{2}(\Gamma_{\mu\nu}^i(p)-{\Gamma^*}^i_{\mu\nu}(p))z^{\mu}z^{\nu}\right]\\
&&\times\left[z^j+\frac{1}{2}{\Gamma^*}_{\mu\nu}^j(p)z^{\mu}z^{\nu}\right]\ \total t\\
&&= \frac{1}{2}\metric_{ij}(p)z^i z^j+\frac{1}{2}\metric_{ij}(p){\Gamma^*}_{\mu\nu}^i(p) z^{\mu}z^{\nu}z^j+\frac{1}{6}\metric_{ij}(p)z^i \left(\Gamma_{\mu\nu}^j(p)-{\Gamma^*}_{\mu\nu}^j(p)\right)z^{\mu}z^{\nu}\,,
\end{eqnarray*}
{where we dropped out the terms of $O(\|\bz\|^4)$.}

Finally, by symmetrizing the indices because of the multiplication $z^i z^j z^k$, we obtain
$$
\varphi(p,q)=\frac{1}{2}\metric_{ij}(p)z^i z^j+\frac{1}{6}\Lambda_{ijk}(p)z^i z^j z^k,
$$
where 
$
\Lambda_{ijk}(p)=2{\Gamma^*}_{ijk}(p)+\Gamma_{ijk}(p)
$
is obtained by recalling that $g_{il}\Gamma^l_{jk}=\Gamma_{ijk}$.

\hfill $\square$

\vspace{.2cm}

In order to show that $\Div$ is consistent to the dual structure $(\metric,\nabla,\nabla^*)$, as well, we supply a Taylor expansion of $\Div(p,q)$ about $\bz=\bxi_q-\bxi_p=0$, where $\{\bxi_p\}$ and $\{\bxi_q\}$ are local coordinates at $p$ and $q$, respectively.

\begin{pro}\label{TaylorD}
When $\|\mathbf{z}\|=\|\bxi_q-\bxi_p\|$ is small, the canonical divergence $\Div$ is expanded as
\begin{equation}
\label{DTaylor}
\Div(p,q)= \frac{1}{2}\metric_{ij}(p) z^i z^j+\frac{1}{6}\Lambda_{ijk}(p) z^i z^j z^k+O(\|\bz\|^4)\ ,
\end{equation}
where $\Lambda_{ijk}(p)=2 \stackrel{*}{\Gamma}_{ijk}+\Gamma_{ijk}$.
\end{pro}
\noindent {\bf Proof}. \ By looking at  Eq. \eqref{divergence} we need to Taylor expand with respect to the local coordinate $\{\bxi\}$ the following factors
$$
\metric_{ij}(\sigma(t)),\quad
\frac{\total}{\total t}\sigma^i(t),\quad
\paralleltransport^j_{\sigma^*_t}\left(\nihat_p(t)\right),
$$
where $\nihat_p(t)=\frac{\total }{\total s}\sigma_t(0)$ {is the velocity vector at $p$ of the $\nabla$-geodesic $\sigma_t(s)\,(0\leq s\leq1)$ connecting $p$ with $\sigma(t)$, whereas $\paralleltransport_{\sigma^*_t}:\tangent_p\Ma\rightarrow\tangent_{\sigma(t)}\Ma$ is the $\nabla$-parallel transport along the $\nabla^*$-geodesic $\sigma_t^*(s)\,(0\leq s\leq 1)$ connecting $p$ with $\sigma(t)$. Here,} $\paralleltransport^j_{\sigma^*_t}$ denotes the $j$th component of the parallel transport with respect to $\nabla$-connection. The Taylor expansion of the metric tensor is given by
\begin{equation}
\label{metricT}
\metric_{ij}(\sigma(t))=\metric_{ij}(p)+t\partial_k \metric_{ij}(p) z^k+O(\|\bz\|^2),
\end{equation}
where $\partial_k=\frac{\partial}{\partial\xi_p^k}$. Consider now the local coordinates $\bxi(t)$ of the geodesic $\sigma(t)$. By Taylor expanding it, we obtain
\begin{equation}
\label{geodD}
\xi^i(t)=\xi^i_p+t z^i+ \frac{t}{2}(1-t)\Gamma_{\mu\nu}^i(p)z^{\mu}z^{\nu}+ O\left(\|\bz\|^3\right),
\end{equation}
where the summation over $\mu$ and $\nu$ is {understood}. Then we have,
\begin{equation}
\label{geodderivD}
\frac{\total}{\total t}\sigma^i(t)=z^i+\frac{1}{2}(1-2t)\Gamma^i_{\mu\nu}(p)z^\mu z^\nu+O\left(\|\bz\|^3\right).
\end{equation}

Consider now the $\nabla$-geodesic $\sigma_t(s)$. From Eq. \eqref{geodderivD} we obtain the following expression for $\nihat^j_p(t)$,
$$
\nihat^j_p(t)=\frac{\total }{\total s}\sigma^j_t(0)= \xi^j(t)-\xi_p^j+\frac{1}{2}\Gamma^j_{\mu\nu}(\sigma(t))(\xi^{\mu}(t)-\xi^{\mu}_p)(\xi^{\nu}(t)-\xi^{\nu}_p).
$$
In addition, we have that
$$
\xi^j(t)-\xi_p^j=t z^j+\frac{t}{2}(1-t)\Gamma_{\mu\nu}^j(p) z^{\mu}z^{\nu}.
$$
Then, we arrive at
\begin{equation}
\label{x}
\nihat^j_p(t)=t z^j+\frac{t}{2}\Gamma_{\mu\nu}^j(p)z^{\mu}z^{\nu}.
\end{equation}
In the end, recalling that {the} $\nabla^*$-geodesic $\sigma_t^*$ connects $p$ with $\sigma(t)$, we use the following Taylor expansion of the $\nabla$-parallel transport along $\sigma^*_t$ \cite{Tapei},
$$
\paralleltransport^j_{\sigma_t^*}(\nihat_p(t))=\nihat^j_p(t)-\Gamma_{\mu\nu}^j(p)(\nihat^{\mu}_p(t))(\sigma^{\nu}(t)-\xi_p^{\nu})
$$
and from Eq. \eqref{x} we obtain
\begin{equation}
\label{partrans}
\paralleltransport^j_{\sigma_t^*}(\nihat_p(t))=t z^j+\frac{t}{2}(1-2t)\Gamma_{\mu\nu}^j(p)z^{\mu}z^{\nu}+O\left(\|\bz\|^3\right).
\end{equation}

We are now ready to provide the Taylor series of the path integral {given in} Eq. \eqref{divergence}. By collecting Eqs. \eqref{metricT}, \eqref{geodderivD} and \eqref{partrans} we obtain the following expression for $\Div(p,q)$,
\begin{align}
\Div(p,q)=& \int_0^1 dt \left[\metric_{ij}(p)+t\partial_k \metric_{ij}(p) z^k+O(\|\bz\|^2)\right]\\
&\times \left[z^i+\frac{1}{2}(1-2t)\Gamma^i_{\mu\nu}(p)z^\mu z^\nu+O\left(\|\bz\|^3\right) \right]\nonumber\\
&\times \left[t z^j+\frac{t}{2}(1-2t)\Gamma_{\mu\nu}^j(p)z^{\mu}z^{\nu}+O\left(\|\bz\|^3\right)\right]\ \total t . \nonumber
\end{align}
Finally, by computing this integral up to $O\left(\|\bz\|^4\right)$ and recalling the relation $\Gamma_{ijk}=g_{li}\Gamma^l_{jk}$ we arrive at
\begin{eqnarray*}
\Div(p,q)&=&\frac{1}{2}\metric_{ij}(p)z^i z^j+\frac{1}{3}\partial_k\metric_{ij}(p)z^iz^j z^k-\frac{1}{6}\Gamma_{ijk}(p)z^{i}z^{j} z^k\\
&=&\frac{1}{2}\metric_{ij}(p)z^i z^j+\frac{1}{6}\left(2\partial_k\metric_{ij}(p)z^i z^j z^k-\Gamma_{ijk}(p)z^{i}z^{j} z^k\right) \ .
\end{eqnarray*}

Now, from the relation $\partial_k\metric_{ij}=\Gamma_{ijk}+\Gamma^*_{ijk}$, we obtain
\begin{equation}\label{divpq}
\Div(p,q)=\frac{1}{2}\metric_{ij}(q) z^i z^j +\frac{1}{6}\left(2 \Gamma^*_{ijk}+\Gamma_{ijk}\right)z^iz^jz^k\ .
\end{equation}
Eq. \eqref{divpq} can be reduced to Eq. \eqref{DTaylor}  by using Eq. \eqref{Tensor}. \hfill $\square$

Eguchi introduced in \cite{eguchi1983} the concept of the {\it contrast function} in order to construct statistical structures on a given manifold $\Ma$. A contrast function $\rho(p,q)$ is defined everywhere on  $\Ma\times\Ma$. For a function $\rho(p,q)$ to be a contrast function, it is required that
$$
\rho(p,q)\geq 0, \qquad \rho(p,q)=0\ \Longleftrightarrow\ p=q\ ,
$$
and
$$
\left.\partial_i\partial_j\ \rho(p,q)\right|_{q=p}=\left.-\partial_i\partial_j^{\prime}\ \rho(p,q)\right|_{q=p}=\metric_{ij}(p)
$$
is strictly positive definite on $\Ma$. If $\rho(p,q)$ is a contrast function,
$$
\Gamma_{ijk}^{\rho}(p):=\left.-\partial_i\partial_j\partial_k^{\prime}\rho(p,q)\right|_{q=p},\qquad \Gamma_{ijk}^{\rho^*}(p):=\left.-\partial^{\prime}_i\partial^{\prime}_j\partial_k\rho(p,q)\right|_{q=p}
$$
define torsion free affine dual connections with respect to the Riemannian metric $\metric$. Now, the purpose of the present article is to recover a given dual structure $(\metric,\nabla,\nabla^*)$ on a manifold $\Ma$ by means of the divergence function since our investigation has been addressed from the very beginning to the inverse problem. In order to pursue this aim, it is enough to consider a contrast function to be defined in a neighborhood of the diagonal set of $\Ma\times\Ma$. We are then in the position to prove Theorem \ref{positivityI} of Section \ref{Overwiev}.

\noindent {\bf \large Proof of Theorem \ref{positivityI}.}\ By means of Eq. \eqref{DTaylor} we obtain that, if $p$ and $q$ are sufficiently close to each other, then
$$
\Div(p,q)\geq 0, \qquad \Div(p,q)=0\ \Longleftrightarrow\ p=q\ .
$$

Analogously, by Proposition \ref{proptaylor} and, in particularly, from Eq. \eqref{Taylor} we have that
$$
\varphi(p,q)\geq 0, \qquad \varphi(p,q)=0\ \Longleftrightarrow\ p=q\ ,
$$
when $p$ and $q$ are sufficiently close to each  other,
as well.

In order to prove that both functions, $\Div(p,q)$ and $\varphi(p,q)$, generate the dual structure, we consider the Taylor series \eqref{Taylor} and \eqref{DTaylor}. Since these are equal, we will show the consistency to the dual structure $(\metric,\nabla,\nabla^*)$ only for $\varphi(p,q)$.  Hence, by differentiating the Taylor series \eqref{Taylor} with respect to $\bxi_q$ we obtain,
\begin{align}
\partial^{\prime}_i \varphi(p,q)=&\metric_{ij}(p)z^j+\frac{1}{2}\Lambda_{ijk}(p)z^j z^k\\
\partial^{\prime}_j\partial^{\prime}_i \varphi(p,q)=& \metric_{ij}(p)+\Lambda_{ijk}(p) z^k\ . \label{metricdiv}
\end{align}
By evaluating $\partial^{\prime}_j\partial^{\prime}_i\varphi(p,q)$ at $\bxi_q=\bxi_p$, i.e. $\bz=0$, we obtain
\begin{equation}
\left.\partial^{\prime}_j\partial^{\prime}_i\varphi(p,q)\right|_{p=q}=\left.\partial^{\prime}_j\partial^{\prime}_i\Div(p,q)\right|_{p=q}=\metric_{ij}(p) \ .
\end{equation}
In addition, we differentiate Eq. \eqref{metricdiv} with respect to $\bxi_p$ and evaluate it at $\bz=0$. This computation leads to
\begin{align}
\left.\partial_k\partial^{\prime}_j\partial^{\prime}_i \varphi(p,q)\right|_{p=q}=&\partial_k\metric_{ij}(p)-\Lambda_{ijk}(p)+\left.\partial_k\Lambda_{ijk}z^k\right|_{z=0}\nonumber\\
=& \partial_k\metric_{ij}(p)-\Lambda_{ijk}(p)\nonumber\\
=&\Gamma_{ijk}+{\Gamma}^*_{ijk}-2{\Gamma}^*_{ijk}-\Gamma_{ijk}=-{\Gamma}^*_{ijk},
\end{align}
where we used Eq. \eqref{Tensor} and the relation $\partial_k\metric_{ij}=\Gamma_{ijk}+\Gamma^*_{ijk}$. Finally, we can conclude that
\begin{equation}
\left.\partial_k\partial^{\prime}_j\partial^{\prime}_i\varphi(p,q)\right|_{p=q}=\left.\partial_k\partial^{\prime}_j\partial^{\prime}_i\Div(p,q)\right|_{p=q}=-{\Gamma}^*_{ijk}(p)\,.
\end{equation}

\hfill $\square$

According to Theorem \ref{positivityI}, the non-negativity of $\Div$ and $\varphi$ holds only in a neighborhood of the diagonal $\Delta$ in $\U\times\U$. However, we will provide sufficient conditions for $\Div$ being non-negative on $\U\times\U$ (see Proposition \ref{positivityplus}).

\vspace{.2cm}

\begin{remark}\label{dualPhi}
Consider the dual functions of ${\varphi}(p,q)$ and $\Div(p,q)$. By interchanging the role of the $\nabla$-connection with the $\nabla^*$-connection, we obtain from Proposition \ref{proptaylor} and Proposition \ref{TaylorD} the following Taylor expansions for ${\varphi}^*$ and $\Div(p,q)$,
\begin{align}
\label{Taylor*}
{\varphi}^*(p,q)=\frac{1}{2}\ \metric_{ij}(p)\ z^i z^j+\frac{1}{6}\Lambda^*_{ijk}(p)\ z^i z^j z^k+O\left(\|\bz\|^4\right),\\
\label{TaylorD*}
{\Div}^*(p,q)=\frac{1}{2}\ \metric_{ij}(p)\ z^i z^j+\frac{1}{6}\Lambda^*_{ijk}(p)\ z^i z^j z^k+O\left(\|\bz\|^4\right),
\end{align}
where
\begin{equation}
\label{Tensor*}
\Lambda^*_{ijk}(p)=2{\Gamma}_{ijk}(p)+\stackrel{*}{\Gamma}_{ijk}(p)\ .
\end{equation}
Then, by repeating the same arguments as in the proof of Theorem \ref{positivityI} we get
\begin{align}
\label{dualconsistency}
 \metric_{ij}(p)&=\left.\partial^{\prime}_i\partial^{\prime}_j {\varphi}^*(\bxi_p,\bxi_q)\right|_{p=q}=\left.\partial^{\prime}_i\partial^{\prime}_j {\Div}^*(\bxi_p,\bxi_q)\right|_{p=q},\\
{\Gamma}_{ijk}(p)&=-\left.\partial^{\prime}_i\partial^{\prime}_j\partial_k {\varphi}^*(\bxi_p,\bxi_q)\right|_{p=q}=-\left.\partial^{\prime}_i\partial^{\prime}_j\partial_k {\Div}^*(\bxi_p,\bxi_q)\right|_{p=q}\ ,
\end{align}
which proves that $\varphi^*$ and $\Div^*$ succeed to recovering the dual structure $(\metric,\nabla,\nabla^*)$, as well.
\end{remark}

\subsection{Conclusion of the proof of Theorem 5}\label{ConclusionTh5}

In order to complete the proof of Theorem \ref{Pi&Pi*orthogonal}, we need to show that $\grad_q\,\Div_p$ is parallel to the $\nabla$-velocity vector $\dot{\sigma}(1)$, where $\sigma(t)\,(0\leq t\leq 1)$ is the $\nabla$-geodesic connecting $p$ with $q$. To this aim we exploit the theory of minimum contrast geometry by Eguchi \cite{eguchi1992}. Recall that the canonical divergence $\Div:\U\times\U\rightarrow\RR$ recovers the dual structure $(\metric,\nabla,\nabla^*)$ of a statistical manifold $\Ma$. In particular, we have that the symbols $\Gamma_{ijk}=\metric\left(\nabla_{\partial_i}\partial_j,\partial_k\right)$ of the $\nabla$-connection are given by
$$
\Gamma_{ijk}(p)=-\left.\partial_i\partial_j\partial_k^{\prime} \Div(\bxi_p,\bxi_q)\right|_{p=q}\,,
$$
where $\{\bxi_p\}$ and $\{\bxi_p\}$ are local coordinates at $p$ and $q$, as usual, and $\partial_i=\partial/\partial\xi_p^i$, $\partial_i^{\prime}=\partial/\partial\xi^i_q$. Consider the hypersurface
$
{\hyperm}_{\Div}:=\{q\in\Ma\,|\,\Div(p,q)=\kappa\}
$
of constant divergence $\Div$ centered at $p\in\Ma$, {where $\kappa$ is a suitably chosen positive constant.}
At each $q\in {\hyperm}_{\Div}$ we can define the {\it minimum contrast leaf} $L_q^{\Div}$ at $q$:
\begin{align}
\label{mincontrastleaf}
& L^{\Div}_q:=\left\{\widetilde{p}\in\U\ |\ \Div(\widetilde{p},q)=\min_{q^{\prime}\in {\hyperm}_{\Div}}\Div(\widetilde{p},q^{\prime})\right\}\ ,
\end{align}
{where $\U$ is a dually convex set as it is defined in Definition \ref{dualconvexset}. In particular, we have that the exponential maps of $\nabla$ and $\nabla^*$ are both diffeomorphisms. Hence, we can find a suitable constant $\varepsilon>0$ such that the $\nabla$-geodesic ball of radius $\varepsilon$ is contained in $\U$ and contains $\hyperm_{\Div}$. This implies, together with the constant rank theorem  that $\hyperm_{\Div}$ is a compact manifold \cite{Lee}. Since the divergence $\Div$ is smooth, we can conclude that it always attains the minimum on $\hyperm_{\Div}$. {Now, since $\hyperm_{\Div}$ is a $n-1$ submanifold of $\Ma$, there exists a normal tubular neighborhood of $\hyperm_{\Div}$ in $\Ma$ \cite{Hirsch}. More precisely, a tubular neighborhood of $\hyperm_{\Div}$ is a pair $(f,\xi)$, where $\xi=(p,\mathrm{E},\hyperm_{\Div})$ is a vector bundle over $\hyperm_{\Div}$ and $f:\mathrm{E}\rightarrow\Ma$ is an embedding such that: i) $\left.f\right|_{\hyperm_{\Div}}=1_{\hyperm_{\Div}}$ where $\hyperm_{\Div}$ is identified with the zero section of $\mathrm{E}$; ii) $f(\mathrm{E})$ is an open neighborhood of $\hyperm_{\Div}$ in $\Ma$. Since $\U$ is a dually convex set, we can achieve a tubular neighborhood by the $\nabla$ exponential map and the tangent bundle $\tangent\hyperm_{\Div}$. In particular, we can find an open set $\mathrm{W}$ of the tangent bundle $\tangent\Ma$ such that $\exp(W)$ is a tubular neighborhood of $\hyperm_{\Div}$. By choosing normal vector fields, we then obtain a normal tubular neighborhood. Recalling that $\exp$ is a diffeomorphism onto its image, we can find a sufficiently small open set $\widetilde{\U}\subset\mathrm{W}$ such that $\exp(\widetilde{\U}\cap N_q)$ is the set of points of $\exp(\widetilde{\U})$ whose nearest point of $\hyperm_{\Div}$ is $q$ \cite{Hirsch}. Here, $N_q$ denotes the normal vector at $q$ in $\tangent_q\hyperm_{\Div}$. Therefore, since $\Div$ is a smooth function, for any point $\widetilde{p}\in\widetilde{\U}$ there exists a unique point $q\in\hyperm_{\Div}$ such that $q$ minimizes $\Div(\widetilde{p},q^{\prime})$ in $\hyperm_{\Div}$  \cite{eguchi1992}.} Let then $q\in\hyperm_{\Div}$ be fixed and consider the set $L_q^{\Div}$. The map $i:L_q^{\Div}\rightarrow\RR$ such that $i(\widetilde{p})=\min_{q^{\prime}\in\hyperm_{\Div}}\Div(\widetilde{p},q^{\prime})$ is injective because of the above assumption. Moreover the differential at any $\widetilde{p}$ of $i$ has rank $1$ because the kernel of $\left.(\total\,i)\right|_{\widetilde{p}}$ is the null set. Hence, the constant-rank theorem implies that $i$ is an immersion and the $L_q^{\Div}$ is a $1$-dimensional submanifold of $\Ma$ \cite{Lee}.} Moreover, $\Ma$ is decomposed (at least locally) into a foliation $\Ma=\cup\{L_q^{\Div}\,|\,q\in\hyperm_{\Div}\}$. Therefore, we can decompose the tangent space of $\mathrm{M}$ at $q$  as follows:
$$
\tangent_q\Ma=\tangent_q\hyperm_{\Div}\oplus\tangent_q L_q^{\Div}\,.
$$ 
Let $q\in \hyperm_{\Div}$ be fixed. By the above assumption we have that the derivative at $q$ of $\Div(p,q)$ along any direction $U$ that is tangent at $q$ to ${\hyperm}_{\Div}$ is zero, i.e. $U_q\Div(p,q)=0$ for all $U\in\Tau({\hyperm}_{\Div})$ and for all $p\in L^{\Div}_q$. Therefore, from Eq. \eqref{metricfromdiv} we have that
$$
\left\langle X,U\right\rangle_q=X_{(p)} U_{(q)} \left.\Div(p,q)\right|_{p=q}=0,\,\quad \forall\, X\in\Tau(L_q^{\Div})\,\,\mbox{and}\,\, \forall\, U\in\Tau({\hyperm}_{\Div})\,,
$$
where $X_{(p)}$ is the derivative at $p$ of $\Div(p,q)$ by $X$ and $U_{(q)}$ is the derivative at $q$ of $\Div(p,q)$ by $U$.
This proves that the tangent space of $L_q^{\Div}$ at $q$ coincides with the normal space of $\hyperm_{\Div}$ at $q$ which means that the minimum contrast leaf $L_q^{\Div}$ is orthogonal at $q$ to the hypersurface $\hyperm_{\Div}$ of constant divergence. In addition, we straightforwardly obtain that
\begin{equation}
\label{zeronormal}
X_{1(p)}\ldots X_{n(p)}U_{(q)}\left.\Div(p,q)\right|_{p=q}=0,\quad \forall\, X_1,\ldots,X_n\in\Tau(L_q^{\Div})\,\,\mbox{and}\,\, U\in\Tau(\hyperm_{\Div})\,.
\end{equation}
Define the map $II:\Tau(L_q^{\Div})\times\Tau(L_q^{\Div})\rightarrow\Tau(\hyperm_{\Div})$ by
\begin{equation}
\label{II}
\metric\left(II(X,Y),U\right)=-\left.X_{(p)}Y_{(p)}U_{(q)}\Div(p,q)\right|_{p=q},\quad \forall\,U\in\Tau(\hyperm_{\Div})\,,
\end{equation}
where $X_{(p)},Y_{(p)}\in \Tau(L_q^{\Div})$ denote the derivatives at $p$ in the  directions $X,Y$ that are tangent to $L_q^{\Div}$. On the contrary, $U_{(q)}$ is the derivative at $q$ in the normal direction to $L_q^{\Div}$. Then $II$ is
the {\it second fundamental tensor} with respect to the $\nabla$-connection. Indeed, it is bilinear and it is decomposed by
$$
\nabla_X Y=\widetilde{\nabla}_X Y+II(X,Y)\,,
$$
where $\widetilde{\nabla}$ is the connection of the $1$-dimensional leaf $L_q^{\Div}$ \cite{eguchi1992}.  According to Eq. \eqref{zeronormal}, we can see that the $\nabla$-second fundamental tensor of $L_q^{\Div}$ vanishes at $q$. This implies that the family of all curves which are orthogonal to the hypersurface $\hyperm_{\Div}$ of constant divergence are all $\nabla$-geodesics starting from $p$ (with a suitable choice of the parameter).

Now, since the gradient at $q$ of the divergence $\grad_q\,\Div_p$ is orthogonal to the hypersurface of constant divergence centered at $p$, we obtain that $\grad_q\,\Div_p$ is parallel to $\dot{\sigma}(1)$, in symbols
$$
\grad_q\,\Div_p\,\parallel\, \dot{\sigma}(1)\,.
$$
Clearly, by interchanging the role of $\Div(p,q)$ with $\Div^*(p,q)$, we get that the gradient at $q$ of the dual divergence $\grad_q\,\Div^*_p$  is parallel to the $\nabla^*$-velocity $\dot{\sigma}^*(1)$, where $\sigma^*(t)\,(0\leq t\leq 1)$ is the $\nabla^*$-geodesic from $p$ to $q$.

We may notice that, since $\Gamma_{ijk}(p)=-\left.\partial_i\partial_j\partial_k^{\prime} \varphi(\bxi_p,\bxi_q)\right|_{p=q}$, we also get $\grad_q\,\varphi\parallel\dot{\sigma}(1)$. This proves that the decomposition \eqref{locdecPiI} is not necessarily an orthogonal one. The same applies to $\varphi^*(p,q)$ and decomposition \eqref{locdecPi*I}.

Recalling the Eq. \eqref{r&Div&Phi*}, i.e. $r_p(q)= \Div_p(q)+\varphi_p^*(q)$, from  Eqs. \eqref{gradpseudo} and \eqref{locdecPiI} we obtain 
$$
\Pi_q(p)+\Pi_q^*(p)=\grad_q\,\Div_p+\Pi^*_q(p)+V_q^*\,,
$$
where $\langle V_q^*,\dot{\sigma}(1)\rangle_q$. Analogously, recalling \eqref{r&Div*&Phi}, i.e. $r_p(q)= \Div^*_p(q)+\varphi_p(q)$, from Eqs. \eqref{gradpseudo} and \eqref{locdecPi*I} we have that
$$
 \Pi_q(p)+\Pi_q^*(p)=\grad_q\,\Div^*_p+\Pi_q(p)+V_q\,,
$$
where $\langle V_q,\dot{\sigma}^*(1)\rangle_q=0$. Finally, by defining $X_q=-V^*_q$ and $X^*_q=-V_q$ we obtain the decompositions \eqref{Piorth} and \eqref{Pi*orth}, i.e.
$$
\Pi_q= \grad_q\ \Div_p + X_q,\quad \Pi^*_q= \grad_q\ \Div^*_p + X^*_q\,.
$$
These are orthogonal decompositions as we know from above that $\grad_q\Div_p$ is parallel to $\dot{\sigma}(1)$ and $\grad_q\Div^*_p$ is parallel to $\dot{\sigma}^*(1)$. In addition, $\langle X_q,\dot{\sigma}(1)\rangle_p=0=\langle X^*_q,\dot{\sigma}^*(1)\rangle_p$.

\section{Symmetry properties of canonical divergence}\label{SymmetryProperties}

In general the canonical divergence $\Div(p,q)$ defined on a general statistical manifold $(\Ma,\metric,\nabla,\nabla^*)$ is not symmetric: $\Div(p,q)\neq\Div(q,p)$. However, it is natural to ask what the relation of $\Div(p,q)$ and $\Div(q,p)$ is. In this section, we address the issue of the symmetry inspired by the property \eqref{duallyflatsymmetry} which is held by the Bregman canonical divergence on dually flat manifolds and prove Theorem \ref{ThmD*&D} of Section \ref{Overwiev}.

\noindent{\bf \large Proof of Theorem \ref{ThmD*&D}}.\  Consider the Taylor expansion \eqref{DTaylor} claimed in Proposition \ref{TaylorD} and then interchange the role of $p$ and $q$, 
\begin{align*}
&\Div(q,p)= \frac{1}{2}\metric_{ij}(q) z^i z^j+\frac{1}{6}\Lambda_{ijk}(q) z^i z^j z^k+O(\|\bz\|^3)\ ,
\\ &\Lambda_{ijk}(q)=2 \stackrel{*}{\Gamma}_{ijk}+\Gamma_{ijk},\quad \mathbf{z}=\bxi_p-\bxi_q\,
\end{align*} 
where $\{\bxi_p\}$ and $\{\bxi_q\}$ are local coordinates at $p$ and $q$, as usual.

We can prove that
\begin{equation}
\label{Dnabla*}
\left.\partial_i\partial_j\partial_k^{\prime}\Div_q(p)\right|_{p=q}=-\Gamma^*_{ijk}(p)\ ,
\end{equation}
where $\partial_i=\frac{\partial}{\partial\xi^i_p}$ and $\partial^{\prime}_i=\frac{\partial}{\partial\xi^i_q}$. Consider now the dual divergence $\Div^*(p,q)$. From Remark \ref{dualPhi} we also know that $\left.\partial_i\partial_j\partial_k^{\prime}\Div^*(p,q)\right|_{p=q}=-\Gamma^*_{ijk}(p)$.  Hence, $\Div(q,p)$ and $\Div^*(p,q)$ generate the dual structure $(\metric,\nabla,\nabla^*)$ in the same way.

Let  $\widetilde{\hyperm}_{\Div}=\{q\in\Ma\ |\ \Div(q,p)=\kappa\}$ be the level hypersuface of $\Div(q,p)$. For a point  $q\in \widetilde{\hyperm}_{\Div}$ we define {the {minimum contrast leaf at $q$}:}
\begin{align*}
& L_q:=\left\{p\in\U\ |\ \Div(q,p)=\min_{q^{\prime}\in \widetilde{\hyperm}_{\Div}}\Div_q(q^{\prime})\right\}\ ,,
\end{align*}
{where $\U$ is a dually convex set.}

According to the theory of {\it minimum contrast geometry} by Eguchi (See \cite{eguchi1983}, \cite{eguchi1992}), we know that $\{L_q\}_{q\in \widetilde{\hyperm}_{\Div}}$ is (locally) a  foliation of $\Ma$ with $1$-dimensional leaves such that
\begin{itemize}
\item[(i)]\ each leaf $L_q$ is orthogonal to $\widetilde{\hyperm}_{\Div}$ at $q$,
\item[(ii)]\ the second fundamental form with respect to {$\nabla^*$} of $L_q$ is zero at $q$.
\end{itemize}
Therefore, the family of all curves orthogonally intersecting $\widetilde{\hyperm}_{\Div}$ are all {$\nabla^*$-geodesics}, with a suitable choice of the parameter, from $p$. 

This implies that the gradient of $\Div(q,p)$ at $q$ is parallel to $\dot{\sigma}^*(1)$, where $\sigma^*$ is the $\nabla^*$-geodesic from $p$ to $q$. Therefore, we can find a constant $c(q)$ depending on $q$ such that
$$
\Pi_q^*=c(q)\ \grad_q\Div(q,p)+X_q\ .
$$  
Then, we obtain that
$$
c(q)\ \grad_q\Div(q,p)=\grad_q\Div^*(p,q)\ .
$$
This implies that there exists a function $f:[0,K]\rightarrow\RR^+$ such that $\Div(q,p)=f(\Div^*(p,q))$.  
This proves claim \eqref{D&D*}. 

The relation $\Div^*(q,p)=f^*\left(\Div(p,q)\right)$ can be proved by means of the methods described above by changing the role of the $\nabla^*$-connection with the $\nabla$-connection. \hfill $\square$ 

\vspace{.2cm}

{In the rest of this section, we establish a close relation among the Phi-functions $\varphi$ and $\varphi^*$ and the canonical divergences $\Div$ and $\Div^*$.}
\begin{theorem}
\label{ThDiv&Phi}
{Let $\varphi(q,p)$ be the function given by Eq. \eqref{phifunction} with $p$ and $q$ reversed and $\Div^*(p,q)$ be the dual canonical divergence (resp.\, $\varphi^*(q,p)$ and $\Div(p,q)$) on a statistical manifold $(\Ma,\metric,\nabla,\nabla^*)$.} {Then, there exists a function $h$ (resp. $h^*$) satisfying the conditions $h(0)=0$ and $h^{\prime}(0)>0$ (resp. $h^*(0)=0$ and $h^{*\prime}(0)>0$) such that
\begin{equation}\label{Phi&D*}
\varphi(q,p)=h\left(\Div^*(p,q)\right)\qquad \left(\mbox{resp.}\;\varphi^*(q,p)=h^*\left(\Div(p,q)\right)\right) \ .
\end{equation}}
\end{theorem}
\noindent {\bf Proof.}\; {Consider the Taylor expansion \eqref{Taylor} and then interchange the role of $p$ with $q$:
\begin{align*}
&\varphi(q,p)= \frac{1}{2}\metric_{ij}(q) z^i z^j+\frac{1}{6}\Lambda_{ijk}(q) z^i z^j z^k+O(\|\bz\|^3)\ ,
\\ &\Lambda_{ijk}(q)=2 \stackrel{*}{\Gamma}_{ijk}+\Gamma_{ijk},\quad \mathbf{z}=\bxi_p-\bxi_q\ .
\end{align*} }
{The divergence $\varphi(q,p)$ induces the dual structure $(\metric,\nabla,\nabla^*)$ on $\Ma$ in the same way as the divergence $\Div^*(p,q)$. In fact, we have that
that
\begin{equation*}
\label{Phinabla*}
\left.\partial_i\partial_j\partial_k^{\prime}\varphi_q(p)\right|_{p=q}=-\Gamma^*_{ijk}(p)\ ,
\end{equation*}
where $\partial_i=\frac{\partial}{\partial\xi^i_p}$ and $\partial^{\prime}_i=\frac{\partial}{\partial\xi^i_q}$. On the other side, from Remark \ref{dualPhi}  we know that $\left.\partial_i\partial_j\partial_k^{\prime}\Div^*_p(q)\right|_{p=q}=-\Gamma^*_{ijk}(p)$, as well.  
According to the theory of  minimum contrast geometry by Eguchi \cite{eguchi1992}, we know that {the family of all curves orthogonally intersecting hypersurfaces of constant $\varphi(q,p)$  are  $\nabla^*$-geodesics, up to a suitable choice of the parametrization, from $p$. This proves that $\grad_q\varphi(q,p)$ is parallel to $\dot{\sigma}^*(1)$, where $\sigma^*$ is the $\nabla^*$-geodesic from $p$ to $q$.}

We can find a constant $\tilde{c}(q)$ depending on $q$ such that
$$
\tilde{c}(q)\ \grad_q\varphi_q=\grad_q\Div^*_p\,.
$$
This shows that there exists a function $h:[0,K]\rightarrow\RR^+$ such that $\varphi(q,p)=h\left(\Div^*(p,q)\right)$ as claimed by Eq. \eqref{Phi&D*}. By repeating the same arguments as above by changing the role of the $\nabla^*$-connection with the $\nabla$-connection, we can prove the relation $\varphi^*(q,p)=h\left(\Div(p,q)\right)$, as well. \hfill $\square$}

\vspace{.2cm}

{In the end, a straightforward application of Theorem \ref{ThmD*&D} and Theorem \ref{ThDiv&Phi} leads to the following relation:
\begin{equation}
\label{Div&W}
\varphi(q,p)=\Upsilon\left(\Div(q,p)\right),\qquad \varphi^*(q,p)=\Upsilon^*\left(\Div^*(q,p)\right)\,,
\end{equation} 
where $\Upsilon:= h\circ f^{-1}$ and $\Upsilon^*:= h^*\circ f^{^*-1}$. Here, $f$ and $f^*$ are provided by Theorem \ref{ThmD*&D} whereas $h$ and $h^*$ are supplied by Theorem \ref{ThDiv&Phi}.}

\section{Proof of Theorem 8}\label{sec4}

Let $(\Ma,\metric,\nabla,\nabla^*)$ be a general statistical manifold and $\U\subset\Ma$ be a dually convex set. In this section, we aim to prove  Theorem \ref{FeaturesDiv} of  Section \ref{Overwiev}. Recall that the divergence introduced in \cite{Ay15} is, on self-dual manifolds, the energy of the geodesic connecting $p$ and $q$. Moreover, such a divergence coincides with the Bregman canonical divergence \eqref{BregmanDiv} on dually flat manifolds. For these reason, we intend to prove Theorem \ref{FeaturesDiv} by establishing a close relation between the canonical divergence $\Div(p,q)$ and the divergence \eqref{Aydivergence} of Ay and Amari. Recall that the latter has been defined by path integration of the vector field $\nihat_t(q)=\exp_{\sigma(t)}^{-1}(q)$ along the $\nabla$-geodesic $\sigma(t)$. In particular, when the $\nabla$-geodesic $\sigma$ goes from $p$ to $q$, such divergence assumes the nice form \eqref{Aydivergence}, i.e.
$$
D(p,q)=\int_0^1 t\|\dot{\sigma}(t)\|^2\ \total t\ .
$$

In order to carry out this comparison, let us consider for each $t\in[0,1]$ a loop $\Sigma_t$ based at $p$ and passing by $\sigma(t)$. We may refer to Fig. \ref{actionfunctional} for figuring out its definition, which is given by
\begin{equation}
\label{loopt}
\Sigma_t(s)=\left\{\begin{array}{ll}
\sigma^*_t(2s), & s\in[0,1/2]\\
\\
\sigma_t(2-2s), & s\in [1/2,1]
\end{array}\right. ,
\end{equation}
{where the $\nabla$-geodesic $\sigma_t(s)\, (0\leq s\leq 1)$ and the $\nabla^*$-geodesic $\sigma_t^*(s)\, (0\leq s\leq 1)$ connect $p$ with $\sigma(t)$, i.e. $\sigma_t(0)=p=\sigma_t^*(0)$ and $\sigma_t(1)=\sigma(t)=\sigma_t^*(1)$.} By means of  Lemma \ref{lemma} in [Appendix \ref{StatManifold}] we know that, if $\Sigma_t$ lies in a sufficiently small neighborhood of $p$, then
\begin{equation}
\label{Parallel&Curvature}
\paralleltransport_{\Sigma_t} \nihat_p(t)=\nihat_p(t)+ \Riemann_{\Sigma_t}\left(\nihat_p^*(t),\nihat_p(t)\right),
\end{equation}
where 
\begin{equation}
\label{RiemannSigmat}
\Riemann_{\Sigma_t}\left(\nihat_p^*(t),\nihat_p(t)\right) :=\int_{B_t} \frac{\paralleltransport\left[ \Riemann\left(\nihat^*(t),\nihat(t)\right)\nihat(t)\right]}{\|\nihat_p^*(t)\wedge\nihat_p(t)\|}\ \total A\,,
\end{equation}
and $\paralleltransport_{\Sigma_t}:\tangent_p\Ma\rightarrow\tangent_p\Ma$ is the $\nabla$-parallel transport along the loop $\Sigma_t$.
{Here, $\nihat^*_p(t)=\stackrel{*}{\exp_p}^{-1}(\sigma(t))$ and $\nihat_p(t)=\exp_p^{-1}(\sigma(t))$. Moreover, $\nihat^*(t)$ and $\nihat(t)$ are the $\nabla$-parallel transport of $\nihat_p^*(t)$ and $\nihat_p(t)$, respectively, from $p$ to each point of the disc $B_t$ defined by $\Sigma_t$, along the unique $\nabla$-geodesic joining them. In addition,}  
$\Riemann$ is the curvature tensor of $\nabla$, and $\paralleltransport$ within the integral denotes the $\nabla$-parallel translation from each point in $B_t$ to $p$ along the unique $\nabla$-geodesic segment joining them.

Now, we represent $ \paralleltransport_{\Sigma_t}$ as the composition of the $\nabla$-parallel transport along $\sigma_t^*$ and the $\nabla$-parallel transport along $\sigma_t$. In particular we can write
$$
\paralleltransport_{\Sigma_t} \nihat_p(t)=\left(\paralleltransport_{\sigma_t}^{-1}\circ \paralleltransport_{\sigma_t^*}\right)\nihat_p(t)\ .
$$
Then, from Eq. \eqref{Parallel&Curvature} we get
$$
\paralleltransport_{\sigma^*_t}\nihat_p(t)=\paralleltransport_{\sigma_t}\nihat_p(t)+\paralleltransport_{\sigma_t}\left[\Riemann_{\Sigma_t}\left(\nihat_p^*(t),\nihat_p(t)\right)\right]\,,
$$
where $\paralleltransport_{\sigma_t^*},\,\paralleltransport_{\sigma_t}:\tangent_p\Ma\rightarrow\tangent_{\sigma(t)}\Ma$ are the $\nabla$-parallel transport along the $\nabla^*$-geodesic $\sigma_t^*$ and the $\nabla$-geodesic $\sigma_t$, respectively. We may notice that $\nihat_p(t)=\exp_p^{-1}(\sigma(t))$ is the velocity vector at $p$ of the $\nabla$-geodesic $\sigma_t$. Therefore, we can write 
\begin{equation}
\label{Pgeod}
\paralleltransport_{\sigma_t}\nihat_p(t)=\dot{\sigma}_t(1)=t\ \dot{\sigma}(t)\ .
\end{equation}
Thus, we obtain
\begin{align}
\Pi_t(p) =\paralleltransport_{\sigma_t^*}\nihat_p(t)
&= \paralleltransport_{\sigma_t}\nihat_p(t)+ \paralleltransport_{\sigma_t}\left[\Riemann_{\Sigma_t}\left(\nihat_p^*(t),\nihat_p(t)\right)\right]\nonumber\\
\label{Pi&nihat}
&= t\,\dot{\sigma}(t)+\paralleltransport_{\sigma_t}\left[\Riemann_{\Sigma_t}\left(\nihat_p^*(t),\nihat_p(t)\right)\right]\,.
\end{align}

{Finally, we can plug the last expression of $\Pi_t(p)$ into the definition \eqref{divergence} of $\Div(p,q)$ and then we obtain}
\begin{equation}
\label{decompositiondivergence}
\Div(p,q)=\int_0^1 t\ \|\dot{\sigma}(t)\|^2\ \total t+\int_0^1 \ \langle\paralleltransport_{\sigma_t}\left[\Riemann_{\Sigma_t}\left(\nihat_p^*(t),\nihat_p(t)\right)\right],\dot{\sigma}(t)\rangle_{\sigma(t)}\ \total t\ .
\end{equation}

This decomposition of $\Div(p,q)$ allows us to provide sufficient conditions for the positivity of $\Div(p,q)$ for all $p,q\in\U$ in a dually convex set $\U$.
\begin{pro}\label{positivityplus}
Consider $p,q\in\U$ in a dually convex set $\U\subset\Ma$. Let us assume the following conditions on the Riemannian curvature tensor $\RC$,
\begin{align}\label{sufficientpositivity}
(i)&\qquad \nabla\ \RC\equiv 0\nonumber\\
(ii)&\qquad \RC\left(X,Y,Y,Y\right)\geq 0\quad \forall\ X,\ Y\in\Tau(\Ma)\ .
\end{align}
Then, we have
\begin{equation}
\label{positivityD}
\Div(p,q)\geq 0\quad \forall\ p,q\in\U,\qquad \Div(p,q)=0\ \Longleftrightarrow\ p=q\ .
\end{equation}
\end{pro}
\noindent{\bf Proof}.\ 
In order to prove this statement, let us consider the decomposition \eqref{decompositiondivergence} of the canonical divergence $\Div(p,q)$. By $\nabla\RC\equiv 0$ we know that the curvature tensor is invariant under all parallel translations with respect to the $\nabla$-connection \cite{Helgason}. Therefore, by Eq. \eqref{RiemannSigmat} and by recalling the definition of $\paralleltransport$ below Eq. \eqref{RiemannSigmat} we obtain 
\begin{eqnarray}
\Riemann_{\Sigma_t}\left(\nihat^*_p(t),\nihat_p(t)\right)&=& \int_{B_t} \frac{\paralleltransport\left[ \Riemann\left(\nihat^*(t),\nihat(t)\right)\nihat(t)\right]}{\|\nihat_p^*(t)\wedge\nihat_p(t)\|}\ \total A\nonumber\\
&=& \int_{B_t} \frac{ \Riemann\left(\paralleltransport\nihat^*(t),\paralleltransport\nihat(t)\right)\paralleltransport\nihat(t)}{\|\nihat_p^*(t)\wedge\nihat_p(t)\|}\ \total A\nonumber\\
&=&\int_{B_t} \frac{ \Riemann\left(\nihat_p^*(t),\nihat_p(t)\right)\nihat_p(t)}{\|\nihat_p^*(t)\wedge\nihat_p(t)\|}\ \total A\nonumber\\
\label{nablaRzero}
&=&\varepsilon_t\  \Riemann\left(\nihat^*_p(t),\nihat_p(t)\right)\nihat_p(t)\ ,
\end{eqnarray}
where
$$
\varepsilon_t:=\frac{\mbox{Area}(B_t)}{\|\nihat_p^*(t)\wedge\nihat_p(t)\|}\ .
$$
Moreover, from Eq. \eqref{Pgeod} we have that
\begin{eqnarray*}
&&\int_0^1 \ \langle\paralleltransport_{\sigma_t}\left[\Riemann_{\Sigma_t}\left(\nihat_p^*(t),\nihat_p(t)\right)\right],\dot{\sigma}(t)\rangle_{\sigma(t)}\ \total t \\
&=& \int_0^1 \ \frac{\varepsilon_t}{t}\ \RC\left(\paralleltransport_{\sigma_t}\nihat^*_p(t),\paralleltransport_{\sigma_t}\nihat_p(t),\paralleltransport_{\sigma_t}\nihat_p(t),\paralleltransport_{\sigma_t}\nihat_p(t)\right)\  \total t\\
& \geq & 0
\end{eqnarray*}
because of Condition $(ii)$ in Eq. \eqref{sufficientpositivity}.
Finally, from Eq. \eqref{decompositiondivergence} we arrive at $\Div(p,q)\geq 0$ for all $p,q\in\U$, with $\Div(p,q)=0$ iff $p=q$.\hfill $\square$

\vspace{.2cm}
By replacing $\nabla$ and $\RC$ in Eq. \eqref{sufficientpositivity} by $\nabla^*$ and $\RC^*$, respectively, we obviously obtain that $\Div^*(p,q)\geq 0$ for all $p,q\in\U$ with $\Div^*(p,q)=0$ iff $p=q$, as well.

\vspace{.2cm}

Recalling the representation \eqref{Aydivergence} of the divergence introduced in \cite{Ay15}, Eq. \eqref{decompositiondivergence} can be rewritten as follows:
\begin{equation}
\label{Div+Ay}
\Div(p,q)\,=\,D(p,q)+\int_0^1 \ \langle\paralleltransport_{\sigma_t}\left[\Riemann_{\Sigma_t}\left(\nihat_p^*(t),\nihat_p(t)\right)\right],\dot{\sigma}(t)\rangle_{\sigma(t)}\ \total t\,,
\end{equation}
where $D(p,q)=\int_0^1\,t\,\|\dot{\sigma}(t)\|^2\,\total t$. In the rest of this section, we consider $\Div(p,q)$ on different classes of statistical manifolds which are ordered according to their generality.

\subsection{Divergence in self dual manifolds and dually flat manifolds}\label{SelfDual&DuallyFlat}

A statistical manifold $(\Ma,\metric,\nabla,\nabla^*)$ is said self-dual when $\nabla=\nabla^*$. In this case {the dualistic structure $(\Ma,\metric,\nabla,\nabla^*)$} reduces to a Riemannian manifold $(\Ma,\metric,\Lc)$ endowed with the Levi-Civita connection. Indeed, from
$$
\Lc=\frac{1}{2}\left(\nabla+\nabla^*\right).
$$
{we obtain $\nabla=\nabla^*=\Lc$. In this case, we then have that the exponential map of $\nabla$ coincides with the exponential map of $\nabla^*$. Therefore, we get
$$
\nihat_p(t)=\exp^{-1}_p(\bar{\sigma}(t))=\stackrel{*}{\exp}_p^{-1}(\bar{\sigma}(t))=\nihat^*_p(t),\quad \forall \ t\in[0,1]\,,
$$
where $\bar{\sigma}(t)$ is the $\Lc$-geodesic from $p$ to $q$. Recalling Eq. \eqref{RiemannSigmat}  we obtain that, in the self-dual case, $\Riemann_{\Sigma_t}\equiv 0$. This follows from the the skew-symmetry of the curvature tensor which implies $\Riemann(\nihat(t),\nihat(t))\equiv 0$ for all $t\in[0,1]$.}
Hence, we have that
$$
\int_0^1 \ \langle\paralleltransport_{\sigma_t}\left[\Riemann_{\Sigma_t}\left(\nihat(t),\nihat(t)\right))\right],\dot{\bar{\sigma}}(t)\rangle_{\sigma(t)}\ \total t\ \equiv 0
$$
and from Eq. \eqref{Div+Ay} we obtain
\begin{equation}
\label{divAyself}
\Div(p,q)=D(p,q)=\int_0^1\ t\ \|\dot{\bar{\sigma}}(t)\|^2\ \total t\,.
\end{equation}

\vspace{.2cm}

In addition, {by noticing that $\bar{\sigma}(t)$ is a $\Lc$-geodesic,} we know from classical Riemannian geometry that the term  $\|\dot{\sigma}(t)\|^2$ is constant with respect to the parameter $t$. {Therefore, we can write
$$
\|\dot{\bar{\sigma}}(t)\|^2=\|\dot{\bar{\sigma}}(0)\|^2=\left\langle\overline{\exp}_p^{-1}(q),\overline{\exp}_p^{-1}(q)\right\rangle_p=d(p,q)^2\,,
$$
where $\overline{\exp}_p$ denotes the exponential map of the Levi-Civita connection $\Lc$ and $d(p,q)$ is the Riemannian distance.}
Thus, performing the integration in Eq. \eqref{divAyself}, we can conclude that the new canonical divergence
corresponds to the energy of the $\Lc$-geodesic $\bar{\sigma}(t)$ from $p$ to $q$, that is  
$$
\Div(p,q)=\frac{1}{2}d(p,q)^2\,.
$$

\vspace{.3cm}

The statistical manifold $(\Ma,\metric,\nabla,\nabla^*)$ is called dually flat when the curvature tensors of $\nabla$ and $\nabla^*$ are zero, i.e. $\Riemann(\nabla)=\Riemann^*(\nabla^*)\equiv 0$. Then, {we can see from Eq. \eqref{RiemannSigmat} that} $\Riemann\equiv 0$ implies $\Riemann_{\Sigma_t}\equiv 0$. Thus, from Eq. \eqref{Div+Ay} we immediately get
$$
\Div(p,q)= D(p,q)=\int_0^1\ t\ \|\dot{\sigma}(t)\|^2\ \total t\,,
$$
{where $\sigma(t)$ is the $\nabla$-geodesic connecting $p$ with $q$, i.e. $\sigma(0)=p$ and $\sigma(1)=q$.}
This proves that also in case of dually flat manifolds our divergence coincides with the one of Ay and Amari.

\vspace{.2cm}

\noindent In a dually flat manifold we can {consider a $\nabla$-affine coordinate system $\{\theta^i\}$ and a $\nabla^*$-affine coordinate system $\{\eta_j\}$ such that $\Gamma_{ijk}(\btheta)=0$ and $\Gamma^*_{ijk}(\bieta)=0$. Here, $\Gamma_{ijk}$ and $\Gamma^*_{ijk}$ denote the connection symbols of $\nabla$ and $\nabla^*$, respectively, and $\btheta=(\theta^1,\ldots,\theta^n)$,  $\bieta=(\eta_1,\ldots,\eta_n)$. In addition, $\{\theta^i\}$ and $\{\eta_j\}$ are dual with respect to the metric tensor $\metric$ in the sense that
$$
\metric\left(\partial_j,\partial^j\right)=\delta_i^j,\quad \partial_i=\frac{\partial}{\partial \theta^i},\quad \partial^j=\frac{\partial}{\partial \eta_j}\,.
$$
The coordinates $\btheta$ and $\bieta$ are connected through the following Legendre transform,
$$
\partial_i\phi=\eta_i,\quad \partial^i\phi^*=\theta^i,\quad \phi+\phi^*-\sum_i \theta^i\eta_i=0\,.
$$

In \cite{Ay15}, the authors showed that their divergence coincides to the divergence \eqref{BregmanDiv} of Bregman type. Therefore, we straightforwardly have that
\begin{align}
\Div(p,q)&=\phi(\btheta_p)+\phi^*(\bieta_q)-\sum_i\theta^i(p)\eta_i(q)\,,\\
\Div^*(p,q)&= \phi(\btheta_q)+\phi^*(\bieta_p)-\sum_i\theta^i(q)\eta_i(p)\,.
\end{align}
In addition, we can see that $\Div^*(p,q)=\Div(q,p)$,  which proves that $\nabla$ and $\nabla^*$ give the same canonical divergence except that $p$ and $q$ are interchanged because of the duality. Such a nice property holds when $\Ma$ is dually flat.

\subsection{Divergence in symmetric statistical manifolds}\label{SymmetricStatisticalmanifold}

In \cite{Kobayashi00}, the authors formulated a potential teoretic property of the divergence $W(p\| q)$ given in Eq. \eqref{HenmiDivergence} by introducing a class of statistical manifolds satisfying the following condition $(S)^*$:
\begin{equation}
\label{S*condition}
(S)^*\,=\,\left\{\begin{array}{ll}
(i) &\qquad  \RC^*(X,Y,Y,Y)=0\quad \forall\ X,Y\in\Tau(\Ma)\\
(ii) & \qquad  \nabla^*\,\RC^* =0,
\end{array}\right\}
\end{equation}
where $\RC^*$ denotes the Riemann curvature tensor of $\nabla^*$. The conditions \eqref{Scondition} and \eqref{S*condition} are a statistical geometric analogue of the concept of symmetric spaces in Riemannian geometry. Hence, we refer to $(\Ma,\metric,\nabla,\nabla^*)$ satisfying conditions $(S)$ and $(S)^*$ as a {\it symmetric statistical manifold}. 

For $p,q\in\U$ in a dually convex set $\U\subset\Ma$, the dual function of the divergence $W(p\| q)$ is defined by
\begin{align}
\label{koba}
W^*(p\|q):= &-\int_0^1 \left\langle \dot{\bar{\sigma}}^*(t),\nihat_t(q)\right\rangle_{\bar{\sigma}^*(t)}\, \total t \,,\quad \nihat_t(q):=\exp^{-1}_{\bar{\sigma}^*(t)}(q)\,.
\end{align}
Here,  $\bar{\sigma}^*(t)\,(0\leq t\leq 1)$ is the $\nabla^*$-geodesic which connects $q$ with $p$, i.e. $\bar{\sigma}^*(0)=q$ and $\bar{\sigma}^*(1)=p$. Moreover, the $\nabla$-geodesic $\bar{\sigma}_t(s)\,(0\leq s\leq 1)$ is such that $\bar{\sigma}_t(0)=\bar{\sigma}(t)$ and $\bar{\sigma}(1)=q$.

Under the condition \eqref{Scondition} the divergence $W^*(p\| q)$ turns out to be a potential function of the vector field $\nihat_t(q)$ \cite{Kobayashi00}. Indeed, it has been proved in \cite{Kobayashi00} that the integral curves of $\grad\, W(\cdot\| q)$ coincide with the $\nabla$-geodesics starting from $q$. In particular, it holds: 
\begin{equation}\label{GradientF}
\nihat_t(q)=\grad_{\bar{\sigma}^*(t)}\, W^{*}(\cdot\| q)\,.
\end{equation}
This implies that the value of $W^*(p\| q)$ depends only on the end points $p,\,q$. 

\subsubsection{Canonical divergence and divergence of Ay and Amari}

In this section we show that, under the condition \eqref{Scondition}, the new canonical divergence coincides with the divergence of Ay and Amari. To this purpose, consider the Eq. \eqref{Div+Ay}. Since $\nabla\RC=0$, from Eq. \eqref{nablaRzero} we can write
\begin{equation}\label{RHenmi}
\Riemann_{\Sigma_t}\left(\nihat^*_p(t),\nihat_p(t)\right)= \varepsilon_t \ \Riemann\left(\nihat^*_p(t),\nihat_p(t)\right)\nihat_p(t),\quad \varepsilon_t=\frac{\mbox{Area}(B_t)}{\|\nihat_p^*(t)\wedge\nihat_p(t)\|}\,,
\end{equation}
where $\nihat_p(t)=\exp_p^{-1}(\sigma(t))$ and $\nihat_p^*(t)=\stackrel{*}{\exp}_p^{-1}(\sigma(t))$. Here, $\sigma(t)$ is the $\nabla$-geodesic such that $\sigma(0)=p$ and $\sigma(1)=q$ for $p,q\in\U$ in a dually convex set $\U$. 
Recall that $\nihat_p(t)$ is the velocity vector at $p$ of the $\nabla$-geodesic $\sigma_t(s)\,(0\leq s\leq 1)$ such that $\sigma_t(0)=p$ and $\sigma_t(1)=\sigma(t)$, namely $\nihat_p(t)=\dot{\sigma}_t(0)$. Moreover, since $\sigma_t$ is a $\nabla$-geodesic we can write $\dot{\sigma}_t(1)=\paralleltransport_{\sigma_t}\nihat_p(t)$, where $\paralleltransport_{\sigma_t}:\tangent_p\Ma\rightarrow\tangent_{\sigma(t)}\Ma$ is the $\nabla$-parallel transport along the curve $\sigma_t$.  In addition, we have $\sigma_t(1)=\sigma(t)$. Hence, we can write $\dot{\sigma}_t=t\,\dot{\sigma}(t)$ or, equivalently, $\dot{\sigma}(t)=\frac{\paralleltransport_{\sigma_t}\nihat_p(t)}{t}$. Consider the second term of the right hand side in Eq. \eqref{Div+Ay}. Under the condition \eqref{Scondition} we can perform the following computation:
\begin{eqnarray}
&&\int_0^1 \ \langle\paralleltransport_{\sigma_t}\left[\Riemann_{\Sigma_t}\left(\nihat_p^*(t),\nihat_p(t)\right)\right],\dot{\sigma}(t)\rangle_{\sigma(t)}\ \total t \nonumber \\
&=& \int_0^1 \varepsilon_t\ \langle\Riemann\left(\paralleltransport_{\sigma_t}\nihat_p^*(t),\paralleltransport_{\sigma_t}\nihat_p(t)\right)\paralleltransport_{\sigma_t}\nihat_p(t),\dot{\sigma}(t)\rangle_{\sigma(t)}\ \total t \nonumber\\
&=& \int_0^1 \frac{\varepsilon_t}{t}\ \langle\Riemann\left(\paralleltransport_{\sigma_t}\nihat_p^*(t),\paralleltransport_{\sigma_t}\nihat_p(t)\right)\paralleltransport_{\sigma_t}\nihat_p(t),\paralleltransport_{\sigma_t}\nihat_p(t)\rangle_{\sigma(t)}\ \total t \nonumber\\
&=& \int_0^1 \frac{\varepsilon_t}{t}\ \RC\left(\paralleltransport_{\sigma_t}\nihat_p^*(t),\paralleltransport_{\sigma_t}\nihat_p(t),\paralleltransport_{\sigma_t}\nihat_p(t),\paralleltransport_{\sigma_t}\nihat_p(t)\right)\ \total t \nonumber\\
&&= 0\ ,\label{RSigmat0}
\end{eqnarray}
where the first equality follows by $\nabla\Riemann\equiv 0$, and the last one follows from $(i)$ in \eqref{Scondition}. As a result, we obtain that  the canonical divergence $\Div(p,q)$ is given by
$$
\Div(p,q)=\int_0^1\ t\ \|\dot{\sigma}(t)\|^2\ \total t,
$$
which corresponds to the divergence of Ay and Amari.

\subsubsection{Canonical divergence and the divergence of Henmi and Kobayashi}

For $p,q$ in a dually convex set $\U$, we compare the new canonical divergence $\Div(p,q)$ with the divergence $W^*(p\| q)$ introduced in \cite{Kobayashi00} by Henmi and Kobayashi. Due to the potential theoretic property held by $W^*(p\| q)$ under the condition \eqref{Scondition} and interchanging the role of $p$ and $q$, we can write
\begin{equation}
\label{W*qp}
W_{*}(q\| p) = -\int_0^1 \left\langle\dot{\bar{\sigma}}_t(0),\dot{\sigma}(t)\right\rangle_{\sigma(t)}\, \total t\,,
\end{equation}
where $\sigma(t)$ is the $\nabla$-geodesic such that $\sigma(0)=p$ and $\sigma(1)=q$, and $\bar{\sigma}_t(s)\,(0\leq s\leq 1)$ is the $\nabla$-geodesic such that $\bar{\sigma}_t(0)=\sigma(t)$, $\bar{\sigma}_t(1)=p$. 

Now, according to Eq. \eqref{Pgeod} and Eq. \eqref{Pi&nihat} we can write
\begin{equation}
\label{Pi&Forcefield}
\Pi_t(p) = \dot{\sigma}_t(1)+\paralleltransport_{\sigma_t}\left[\Riemann_{\Sigma_t}(\nihat^*_p(t),\nihat_p(t))\right]\,,
\end{equation}
where $\Riemann_{\Sigma_t}$ is given by Eq. \eqref{RiemannSigmat} and $\sigma_t(s)\,(0\leq s\leq 1)$ is the $\nabla$-geodesic such that $\sigma_t(0)=p$ and $\sigma_t(1)=\sigma(t)$. Recall that the canonical divergence $\Div(p,q)$ is defined in terms of the inner product of the velocity vector $\dot{\sigma}(t)$ with the vector $\Pi_t(p)$ . Under the condition \eqref{Scondition} we can use Eq. \eqref{RSigmat0} and then write
\begin{equation}
\label{Divsigmat}
\Div(p,q)=\int_0^1\,\langle\Pi_t(p),\dot{\sigma}(t)\rangle_{\sigma(t)}\,\total t=\int_0^1\ \langle\dot{\sigma}_t(1),\dot{\sigma}(t)\rangle_{\sigma(t)}\ \total t\,.
\end{equation}

Notice that $\bar{\sigma}_t(s)$ corresponds to the  reversely oriented $\nabla$-geodesic $\sigma_t(s)$. Then we have
\begin{equation}
\label{Gradientsigmadot}
\dot{\bar{\sigma}}_t(0)=-\dot{\sigma}_t(1)\ .
\end{equation}
Thus, from Eq. \eqref{W*qp} and Eq. \eqref{Divsigmat} we conclude that
\begin{equation}
\label{Ay&Koba}
\Div(p,q)=W^*(q\| p)\,.
\end{equation}

\begin{remark}
The divergence $W^*(q\| p)$ is the potential function for the vector field $\dot{\sigma}(1)$, where $\sigma(t)\,(0\leq t\leq 1)$ is the $\nabla$-geodesic such that $\sigma(0)=p$ and $\sigma(1)=q$. Hence, Eq. \eqref{Ay&Koba} implies that the canonical divergence $\Div(p,q)$ is the potential function for $\dot{\sigma}(1)$, as well.

This claim is confirmed also by the orthogonal decomposition \eqref{Piorth} of $\Pi_q(p)$, namely $\Pi_q(p)=\grad_q\,\Div(p,q)+X_q$ with $\langle X_q,\dot{\sigma}(1)\rangle_q=0$. Indeed, under the condition \eqref{Scondition}, consider the decomposition \eqref{Pi&Forcefield} of the vector $\Pi_q(p)$ when we set $t=1$, 
\begin{eqnarray*}
\label{PiSdecomposition}
\Pi_q(p)&=&\dot{\sigma}(1)+\paralleltransport_{\sigma}\Riemann(\nihat_p^*(q),\nihat_p(q))\nihat_p(q)\\
&=&\dot{\sigma}(1)+\Riemann(\paralleltransport_{\sigma}\nihat_q^*(p),\paralleltransport_{\sigma}\nihat_p(q))\paralleltransport_{\sigma}\nihat_p(q)\,,
\end{eqnarray*}
where we used the condition $\nabla\,\Riemann=0$ and $\nihat_p(q)=\dot{\sigma}(0)$. Notice that $\paralleltransport_{\sigma}\nihat_p(q)=\dot{\sigma}(1)$, because $\sigma$ is a $\nabla$-geodesic. Thus, from $(i)$ of the condition \eqref{Scondition} we have that
$$
\left\langle\Riemann(\paralleltransport_{\sigma}\nihat_q^*(p),\paralleltransport_{\sigma}\nihat_p(q))\paralleltransport_{\sigma}\nihat_p(q),\dot{\sigma}(1)\right\rangle_q=\RC\left(\paralleltransport_{\sigma}\nihat_q^*(p),\dot{\sigma}(1),\dot{\sigma}(1),\dot{\sigma}(1)\right)=0\,,
$$
which proves that $\Riemann(\paralleltransport_{\sigma}\nihat_q^*(p),\paralleltransport_{\sigma}\nihat_p(q))\paralleltransport_{\sigma}\nihat_p(q)$ is orthogonal to $\dot{\sigma}(1)$.

Finally, from the uniqueness of the orthogonal decomposition \eqref{Piorth} it follows that, under the condition \eqref{Scondition}, we have
$$
\grad_q\,\Div(p,q)\,=\,\dot{\sigma}(1),\quad X_q=\Riemann(\paralleltransport_{\sigma}\nihat_q^*(p),\dot{\sigma}(1))\dot{\sigma}(1)\,.
$$
This shows that $\nabla$-geodesics starting from $p$ are integral curves of $\grad_q\,\Div(p,\cdot)$,
\begin{equation}
\label{gradDivvsgeodesic}
\grad\,\Div(p,\cdot)=\left.\frac{\total}{\total t}\right|_{s=1}\sigma_t(s)\,,
\end{equation}
where $\sigma_t(s)$ is a $\nabla$-geodesic such that $\sigma_t(0)=p$.
\end{remark}

\begin{remark}
Actually, the function $W^*(q\|p)$ is the dual divergence \eqref{HenmiDivergence} introduced in \cite{Kobayashi00}. There, the following symmetric property has been proved,
\begin{equation}
\label{kobasymmetry}
W^*(q\|p)=\Psi\left(W(p\|q)\right),
\end{equation}
where $\Psi$ is a function such that $\Psi(0)=0$ and $\Psi^{\prime}(0)=1$. Finally, from Eq. \eqref{Ay&Koba} and Eq. \eqref{kobasymmetry} we obtain the following connection between $\Div(p,q)$ and $W(p\| q)$,
\begin{equation}
\label{nicovskoba}
\Div(p,q)=\Psi\left(W(p\|q)\right)\,.
\end{equation}
\end{remark}

\subsubsection{On the failure of the Phi function to be a canonical divergence}\label{PhiinStatMan}

In this section, we show that the Phi-function does not coincide with the divergence of Ay and Amari nor with the one of Henmi and Kobayashi. Let 
 $(\Ma,\metric,\nabla,\nabla^*)$ be a symmetric statistical manifold and $\U\subset\Ma$ be a dually convex set. For $p,q\in\U$, the Phi-function $\varphi(p,q)$ is given by
$$
\varphi(p,q)=\int_0^1\,\langle\Pi_t(p),\dot{\sigma}^*(t)\rangle_{\sigma^*(t)}\,\total t\,,\quad \Pi_t(p)=\paralleltransport_{\sigma_t^*}\exp_p^{-1}(\sigma^*(t))\,.
$$
Here $\sigma^*(t)\,(0\leq t\leq 1)$ is the $\nabla^*$-geodesic such that $\sigma^*(0)=p$ and $\sigma^*(1)=q$ whereas $\sigma_t^*(s)\,(0\leq s\leq 1)$ is the $\nabla^*$-geodesic such that $\sigma_t^*(0)=p$ and $\sigma_t^*(1)=\sigma^*(t)$. Moreover, $\paralleltransport_{\sigma_t^*}:\tangent_p\Ma\rightarrow\tangent_{\sigma^*(t)}\Ma$ is the $\nabla$-parallel transport along $\sigma_t^*$.

Now consider for each $t\in[0,1]$ the loop $\Sigma^*_t$ based at $p$ and passing by $\sigma^*(t)$. This is defined as follows,
\begin{equation}
\label{loopt*}
\Sigma^*_t(s)=\left\{\begin{array}{ll}
\sigma_t(2s), & s\in[0,1/2]\\
\\
\sigma^*_t(2-2s), & s\in [1/2,1]
\end{array}\right. ,
\end{equation}
where  $\sigma_t(s)\, (0\leq s\leq 1)$ is the $\nabla$-geodesic such that  $\sigma_t(0)=p$ and $\sigma_t(1)=\sigma^*(t)$. By means of Lemma \ref{lemma} in the Appendix \ref{StatManifold} we know that, if $\Sigma^*_t$ lies in a sufficiently small neighborhood of $p$, then
\begin{equation}
\label{Parallel&CurvatureSigma*}
\paralleltransport_{\Sigma^*_t} \nihat_p(t)=\nihat_p(t)+ \Riemann_{\Sigma^*_t}\left(\nihat_p(t),\nihat^*_p(t)\right).
\end{equation}
Notice that $\Riemann_{\Sigma^*_t}\left(\nihat_p(t),\nihat^*_p(t)\right)$ is defined by Eq. \eqref{RiemannSigmat} where $\nihat_p(t)=\exp_p^{-1}(\sigma^*(t))$ and $\nihat_p^*(t)=\stackrel{*}{\exp}_p^{-1}(\sigma^*(t))$. Here, $\paralleltransport_{\Sigma^*_t}=\paralleltransport_{\sigma_t^*}^{-1}\circ\paralleltransport_{\sigma_t}$, and $\paralleltransport_{\sigma_t^*},\paralleltransport_{\sigma_t}:\tangent_p\Ma\rightarrow\tangent_{\sigma^*(t)}\Ma$ denote the $\nabla$-parallel transport along $\sigma_t^*$ and the $\nabla$-parallel transport along $\sigma_t$, respectively. Therefore, from Eq. \eqref{Parallel&CurvatureSigma*} we can write
$$
\paralleltransport_{\sigma_t}\nihat_p(t)=\paralleltransport_{\sigma_t^*}\nihat_p(t)+\paralleltransport_{\sigma_t^*}\left[\Riemann_{\Sigma^*_t}\left(\nihat_p(t),\nihat^*_p(t)\right)\right]\,.
$$
Recalling the definition of $\Pi_t(p)\in\Tau(\sigma^*)$ we get
$$
\Pi_t(p)=\paralleltransport_{\sigma_t^*}\nihat_p(t)=\dot{\sigma}_t(1)-\paralleltransport_{\sigma_t^*}\left[\Riemann_{\Sigma^*_t}\left(\nihat_p(t),\nihat^*_p(t)\right)\right]\,,
$$
where we used $\dot{\sigma}_t(1)=\paralleltransport_{\sigma_t}\nihat_p(t)$ since $\nihat_t(p)=\dot{\sigma}_t(0)$ and $\sigma_t(s)$ is the $\nabla$-geodesic connecting $p$ with $\sigma^*(t)$. Notice that from Eq. \eqref{gradDivvsgeodesic} we can write
$$
\dot{\sigma}_t(1)=\grad_{\sigma^*(t)}\,\Div_p(\sigma^*(t))\,.
$$
In addition, since $\nabla\,\Riemann=0$, from Eq. \eqref{RHenmi} we have 
\begin{eqnarray*}
\paralleltransport_{\sigma_t^*}\left[\Riemann_{\Sigma^*_t}\left(\nihat_p(t),\nihat^*_p(t)\right)\right] &=& \varepsilon_t \paralleltransport_{\sigma_t^*}\left[\Riemann(\nihat_p(t),\nihat_p^*(t))\nihat_p(t)\right]\\
&=&-\varepsilon_t\paralleltransport_{\sigma_t^*}\left[\Riemann(\nihat^*_p(t),\nihat_p(t))\nihat_p(t)\right]\,.
\end{eqnarray*}
Thus, we can write
\begin{equation}
\label{PionSigma*}
\Pi_t(p)=\grad_{\sigma^*(t)}\,\Div_p(\sigma^*(t))+\varepsilon_t\paralleltransport_{\sigma_t^*}\left[\Riemann(\nihat_p^*(t),\nihat_p(t))\nihat_p(t)\right]\,.
\end{equation}
Notice that, since both, $\sigma_t^*$ and $\sigma^*$, are $\nabla^*$-geodesics and $\sigma_t^*(1)=\sigma^*(t)$, we also have 
$$
t\,\dot{\sigma}^*(t)=\dot{\sigma}_t^*(1)=\paralleltransport^*_{\sigma^*_t}\dot{\sigma}_t^*(0)=\paralleltransport^*_{\sigma^*_t}\nihat_p^*(t)\,.
$$ 
Here, $\paralleltransport^*_{\sigma^*_t}:\tangent_p\Ma\rightarrow\tangent_{\sigma^*(t)}\Ma$ is the $\nabla^*$-parallel transport along the $\sigma_t^*$.
By plugging Eq. \eqref{PionSigma*} into the definition Eq. \eqref{phifunction} of $\varphi(p,q)$ we then obtain
\begin{align}
\label{Phi&DivStatManifold}
\varphi(p,q)&=\Div(p,q)+\int_0^1\frac{\varepsilon_t}{t}\left\langle\paralleltransport_{\sigma_t}\left[\Riemann(\nihat_p^*(t),\nihat_p(t))\nihat_p(t)\right],\paralleltransport^*_{\sigma_t^*}\nihat_p^*(t)\right\rangle_{\sigma^*(t)}\,\total t\nonumber\\
&=\Div(p,q)+\int_0^1\frac{\varepsilon_t}{t}\left\langle\Riemann(\nihat_p^*(t),\nihat_p(t))\nihat_p(t),\nihat_p^*(t)\right\rangle_p\total t\,,
\end{align}
where we employed the invariance of the inner product under the parallel transport with respect to dual connections. {Note that in the case of dually flat statistical manifold, the curvature tensor $\Riemann$, and therefore the integral on the right hand side of \eqref{Phi&DivStatManifold}, vanishes. This shows that in the dually flat case $\Div$ and $\varphi$ coincide. {Moreover, we may observe that in the self-dual case, i.e. when the two linear connections coincide with the Levi-Civita connection $\nabla=\nabla^*=\Lc$, also the vector fields $\nihat_p(t)$ and $\nihat_p^*(t)$ coincide. Consequently, the curvature tensor $\Riemann$ is zero and we can see from  \eqref{Phi&DivStatManifold}  that  $\Div$ and $\varphi$ coincide, as well.  However, in general they will be different. Finally, we can conclude that, on a symmetric statistical manifold, the Phi-function $\varphi(p,q)$ does not correspond  to the divergence of Ay and Amari nor to the divergence of Henmi and Kobayashi, unless special cases when the  integral on the right hand side of \eqref{Phi&DivStatManifold} vanishes.}

\begin{remark}
From Eq. \eqref{Phi&DivStatManifold}, we can see that 
 a sufficient condition to getting $\varphi(p,q)=\Div(p,q)$ on a symmetric statistical manifold is obtained by requiring that
\begin{equation}
\label{SectionalCurv}
\left\langle\Riemann(\nihat_p^*(t),\nihat_p(t))\nihat_p(t),\nihat_p^*(t)\right\rangle_p=\RC(\nihat_p^*(t),\nihat_p(t),\nihat_p(t),\nihat_p^*(t))=0\,.
\end{equation}
Such a condition amounts to require that the sectional curvature of the plane generated by $\nihat_p(q)$ and $\nihat_p^*(q)$ is zero.

Furthermore, if $\varphi(p,q)=\Div(p,q)$ we obtain from Eq. \eqref{r&Div*&Phi}
$$
r(p,q)=\Div(p,q)+\Div^*(p,q)\,.
$$
We conjecture that Eq. \eqref{SectionalCurv} plus $\nabla\,\Riemann\equiv 0$ are sufficient conditions for the canonical divergence $\Div(p,q)$ to achieve the symmetry $\Div(p,q)=\Div^*(q,p)$ as in the case of the Bregman canonical divergence on dually flat manifolds.
\end{remark}

 \section{Discussions and concluding Remarks}\label{Conclusions}

In this work we introduced a new divergence by resorting to an extensive investigation of the geodesic geometry in Information Geometry. Here, the natural object of study is a smooth manifold $\Ma$ endowed with a dualistic structure $(\metric,\nabla,\nabla^*)$. This is given in terms of a metric tensor $\metric$ and a couple of linear torsion-free connections $\nabla$ and $\nabla^*$ on the tangent bundle $\tangent\Ma$ which are dual with respect to $\metric$ in the sense that Eq. \eqref{dualconnection} holds true. Inspired by the classical theory built in Riemannian geometry around the Gauss Lemma, we introduced a pseudo-squared-distance which is obtained from the pseudo-energy \eqref{pseudoenergy} of $\nabla^*$ geodesics as well as from the pseudo-energy \eqref{pseudoenergy*} of $\nabla^*$ geodesics (see Section \ref{GeodGeometry} for more details). Then, the analysis of the first geodesic variation of the pseudo-energy $\Lag(\sigma^*)$ has led to Theorem \ref{thmorthog}. This supplies the extension of the Gauss Lemma to Information Geometry. In particular, for every $p,q$ in a dually convex set $\U\subset\Ma$ we can write down the pseudo-squared-distance $r(p,q)$ between them. Hence, we have proved that for all $q$ in the hypersurface $\hyperm_p$ of constant pseudo-squared-distance $r(p,q)$ centered at $p$, the orthogonal ray of $\hyperm_p$ at $q$  is generated by the sum $\Pi_q+\Pi_q^*$. The tangent vectors $\Pi_q(p)$ and $\Pi_q^*(p)$ are defined in terms of $\nabla$ and $\nabla^*$ geodesics and $\nabla$ and $\nabla^*$ parallel transports (see Eq. \eqref{P} and Eq. \eqref{P*} for more details). Crucially based on Theorem \ref{thmorthog}, we have proved Theorem \ref{thmgradient} which asserts that the pseudo-squared-distance $r(p,q)$ is the potential function of the sum $\Pi+\Pi^*$. From the claim $\grad_q\,r(p,q)=\Pi_q(p)+\Pi^*_q(p)$ we can indeed write
$$
r(p,q)=\int_0^1\ \langle \Pi_t(p),\dot{\gamma}(t)\rangle_{\gamma(t)}\,\total t+\int_0^1\  \langle \Pi_t^*(p),\dot{\gamma}(t)\rangle_{\gamma(t)}\,\total t\,,
$$
for any arbitrary path $\gamma:[0,1]\rightarrow\U$ such that $\gamma(0)=p$ and $\gamma(1)=q$. Here, $\Pi_t(p),\Pi_t^*(p)\in\Tau(\gamma)$ are vector fields on $\gamma$ defined by Eq. \eqref{vectorfieldPiI} and Eq. \eqref{vectorfieldPI*I}. At this point we addressed the definition of the novel divergence $\Div(p,q)$ and its dual function $\Div^*(p,q)$. First of all, we resorted to a classical result in Information Geometry which shows that the sum of the canonical divergence of Bregman type \eqref{BregmanDiv} and its dual function corresponds, on dually flat manifolds, to the function $r(p,q)$ \cite{Amari00}. By combining this classical statement to the potential property of $r(p,q)$, we then introduced $\Div(p,q)$ as the $\nabla$-geodesic path integral of $\Pi$ and $\Div^*(p,q)$ as the $\nabla^*$-geodesic path integral of $\Pi^*$ aiming to obtain  that $\grad\,\Div=\Pi$ and $\grad\,\Div^*=\Pi^*$.  Unfortunately, it turns out that, in general, $\grad\,\Div\neq \Pi$ as well as $\grad\,\Div^*\neq \Pi^*$. However, we succeeded to supply orthogonal decompositions of $\Pi$ and $\Pi^*$ in terms of $\grad\,\Div$ and $\grad\,\Div^*$, respectively (see Theorem \ref{Pi&Pi*orthogonal}).

In order to prove Theorem \ref{Pi&Pi*orthogonal} we introduced two further functions $\varphi$ and $\varphi^*$ (Phi-functions herein the paper) and we proved that these functions are intrinsically characterized by the local decomposition of $\Pi$ and $\Pi^*$. In particular, $\varphi(p,q)$ is obtained very naturally by the decomposition of $\Pi$ in terms of a gradient vector field and another vector field which is orthogonal to $\nabla^*$-geodesics whereas $\varphi^*(p,q)$ is characterized by the decomposition of $\Pi^*$ in terms of a gradient vector field and another vector field which is orthogonal to $\nabla$-geodesics (see Theorem \ref{localdecompositionPiI}).

Both, the canonical divergences $\Div\,,\Div^*$ and the Phi-functions $\varphi\,,\varphi^*$, are consistent  with the dualistic structure $(\metric,\nabla,\nabla^*)$ according to Eq. \eqref{metricfromdiv} and Eq. \eqref{dualfromDiv} (see Theorem \ref{positivityI}). Moreover, they are complementary in the sense that we can write $r(p,q)=\Div(p,q)+\varphi^*(p,q)=\varphi(p,q)+\Div^*(p,q)$. Such a complementarity is finally exploited to accomplish the proof of Theorem \ref{Pi&Pi*orthogonal}.

In this article, the divergence $\Div(p,q)$ is highlighted as a suitable one towards tha definition of a canonical divergence on general statistical manifolds. Firslty, we showed that, in the self-dual case where $\nabla=\nabla^*$ coincides with the Levi-Civita connection, $\Div(p,q)$ corresponds to the energy of the geodesic connecting $p$ with $q$. Moreover, when $(\metric,\nabla,\nabla^*)$ is dually flat, the new divergence $\Div(p,q)$ reduces to the canonical divergence \eqref{BregmanDiv} of Bregman type. We then used the divergence introduced in \cite{Ay15} and the divergence defined in \cite{Kobayashi00} as benchmark of our proposal. In particular, we succeeded to prove that, on the class of symmetric statistical manifolds, which constitutes a statistical geometric analogue of the concept of symmetric spaces in Riemannian geometry, $\Div(p,q)$ coincides with the divergence of Ay and Amari. Furthermore, it also coincides with the dual of the potential function  by Henmi and Kobayashi but with the role of $p,\,q$ interchanged. Unlike the function $\Div(p,q)$, we remarked that the Phi-function $\varphi(p,q)$ does not coincide nor with the divergence of Ay and Amari nor with the divergence of Henmi and Kobayashi. This put the Phi-functions in the background with respect to $\Div(p,q)$ to provide a general definition of canonical divergence. For these reasons, we may select the function $\Div(p,q)$, { instead of $\varphi(p,q)$, as the suitable one towards the definition of a} canonical divergence on a general statistical manifold.

In this article, we also addressed our investigation to the symmetry property of $\Div(p,q)$ along the line put forward by the canonical divergence of Bregman type originally introduced on dually flat manifolds (see Eq. \eqref{duallyflatsymmetry} for more details). However, we succeeded to prove only a weak version of this symmetry property and the general problem is still open. According to our analysis, we may conjecture that this issue is closely related to the decomposition of the pseudo-squared-distance $r(p,q)$ in terms of $\Div(p,q)$ and $\Div^*(p,q)$. Once again, we can support this conjecture from the classical theory previously developed by Amari and Nagaoka in \cite{Amari00}. This will constitute the object of study of a forthcoming investigation.

\vspace{.2cm}

\noindent Several examples of divergences can be found in the literature arising from a wide range of physical sciences.  In \cite{Ciaglia17} a divergence is defined as the solution of the Hamilton-Jacobi problem associated with a canonical Lagrangian defined in $\tangent\Ma$. In \cite{Zhang17} by resorting to the dual structure of the Hamiltonian and Lagrangian formulation of mechanics in  $\tangent^*\Ma$ and $\tangent\Ma$, it is established that the divergence function agrees with the exact discrete Lagrangian up to third order if and only if $\Ma$ is a Hessian manifold. The new divergence introduced in this manuscript is  based on an extensive analysis of the geodesic geometry of a general statistical manifold. For this reason, it turns out to be intrinsically related to the dualistic structure $(\metric,\nabla,\nabla^*)$. Moreover, we proved that the new divergence satisfies all the basic requirements to be a canonical divergence according to pattern laid down in \cite{Ay15} and \cite{Ay17}. In addition, it coincides, in some cases, with the divergence of Ay and Amari and it is closely related to  the divergence of Henmi and Kobayashi. 

{ We conclude our paper with an important note.
Even though we refer to our new proposal as a canonical divergence it is meant to be one candidate which we propose in addition to a number of already existing candidates within the search for the most natural divergence on a general statistical manifold. The main contribution of this paper is to highlight a particular geometric perspective within that search, which originates from a refined analysis of our information-geometric generalisation of the celebrated Gauss Lemma. We presented a number of observations and useful derivations which are interesting in their own right. However, along the chain of arguments we were faced with various degrees of freedom which required  particular choices among equally natural possibilities. Therefore, we do not claim that our proposal represents the most natural divergence. However, we do believe that the presented derivations will help moving forward within the search for a general canonical divergence. }

\appendix

\section{Differential Geometry of Statistical Manifolds}\label{StatManifold}

In this section we review useful tools of differential geometry of statistical manifolds which are of relevance to our work. We describe classes of statistical manifolds in terms of curvature tensor features of them and give description of {\it conjugate symmetric} and {\it dually flat} statistical manifolds. For a more detailed presentation we refer to \cite{Lauritzen87}, \cite{Amari16} and \cite{Ay17}. A statistical manifold $\esse=(\Ma,\metric,\nabla,\nabla^*)$ is the datum of a $C^{\infty}$ manifold $\Ma$, a metric tensor $\metric$ and two torsion-free affine connections  $\nabla$ and $\nabla^*$ such that Eq. \eqref{dualconnection} holds true. Let us recall that an affine connection $\nabla$ on $\Ma$ is a linear connection on the tangent bundle $\tangent\Ma$, 
$$\nabla: \Tau(\Ma)\times \Tau(\Ma)\rightarrow \Tau(\Ma),\quad(X,Y)\mapsto \nabla_X Y\ ,$$
such that
\begin{eqnarray*}
 \nabla_{f X_1+g X_2} X=f\nabla_{X_1}X+g\nabla_{X_2}X && \forall X_1,X_2,X\in\Tau(\Ma) \ \mbox{and}\ f,g\in C^{\infty}(\Ma)\\
 \nabla_{X}( a X_1+b X_2)=a\nabla_{X}X_1+b\nabla_{X}X_2 && \forall X_1,X_2,X\in\Tau(\Ma) \ \mbox{and}\ a,b\in \RR\\
 \nabla_X (f Y)=f\nabla_X Y+X(f)Y \quad  && \forall X, Y\in\Tau(\Ma) \ \mbox{and}\ f\in C^{\infty}(\Ma)\ .
\end{eqnarray*}
Roughly speaking, an affine connection is directional derivative of vector fields. In particular, $\nabla_X Y$ is the change
of $Y$ in the direction of $X$. The rule for comparing vectors in two distinct tangent spaces $\tangent_p\Ma$ and
$\tangent_q\Ma$ is established by the notion of parallel transport.

Let us now introduce such a notion by relying on a smooth curve $\gamma:[0,1]\rightarrow\Ma$ on $\Ma$. A vector field along $\gamma$ is a smooth map $V:[0,1]\rightarrow\tangent\Ma$ such that $V(t)\in\tangent_{\gamma(t)}\Ma$ for all $t\in [0,1]$. Let $\Tau(\gamma)$ the space of all vector fields along $\gamma$, then the {\it covariant derivative} $\nabla_t:\Tau(\gamma)\rightarrow\Tau(\gamma)$ of $V\in\Tau(\gamma)$ along $\gamma$ is defined in terms of the connection $\nabla$ as $\nabla_t V(t):=\nabla_{\dot{\gamma}(t)}\widetilde{V}$, where $\widetilde{V}$ is the extension of $V$ to $\Tau(\Ma)$. A vector field $V\in\Tau(\gamma)$ is said to be {\it parallel} along $\gamma$ with respect to $\nabla$ if $\nabla_tV(t)\equiv 0$ for all $t\in[0,1]$. In this case, a basic result in Calculus allows us to consider the isomorphism
\begin{equation}
\label{parallel}
\paralleltransport_{\gamma}:\tangent_{\gamma(t_0)}\Ma\rightarrow \tangent_{\gamma(t)}\Ma,\quad V \mapsto \paralleltransport_{\gamma}(V)
\end{equation}
where $\paralleltransport_{\gamma}(V):= V(t)$ and $V\in\Tau(\gamma)$ is the unique parallel vector along $\gamma$ such that $V(t_0)\equiv V$. Likewise, we have the parallel transport with respect to the $\nabla^*$-connection,
\begin{equation}
\label{parallel*}
\paralleltransport^*_{\gamma}:\tangent_{\gamma(t_0)}\Ma\rightarrow \tangent_{\gamma(t)}\Ma,\quad V \mapsto \paralleltransport^*_{\gamma}(V)\ .
\end{equation}
For the sake of simplicity, from here on, we only refer to the $\nabla$ connection. All the concepts that we will be describing can be naturally passed to the $\nabla^*$-connection. 

The expression of $\nabla$ connection in local coordinates ${\xi_p^i}$ at $p\in\Ma$ is given in terms of the local basis $\{\partial_i\}_p$ ($\partial_i=\partial\slash\partial\xi_p^i$) of the tangent space $\tangent_p\Ma$ by means of the Christoffel's symbols $\Gamma_{ij}^k$,
$$
\nabla_{\partial_i}\partial_j=\Gamma_{ij}^k\partial_k,
$$
where we adopted Einstein convention according
to which whenever an index appears in an expression
as upper and lower index, we sum over that index.  The same applies to the $\nabla^*$-connection, i.e. $\nabla^*_{\partial_i}\partial_j={\Gamma}_{ij}^{*k}  \partial_k $.

By relying on local coordinates $\{\xi\}$, we can also give the local expression of the parallel transport $\paralleltransport$. Consider $X\in\Tau(\gamma)$, then we have $X(t)=X^i(t)\partial_i(t)$, where $\{\partial_i(t)\}$ is a local frame $\gamma(t)$. Then, we have that
\begin{equation}
\label{localparallel}
\frac{\total X^k(t)}{\total t}+\Gamma^k_{ij}(\gamma(t))\dot{\gamma}^i(t) X^j(t)=0\ .
\end{equation}
It is clear from Eq. \eqref{localparallel} that whenever we specify one initial condition $X^i(0)=X_p^i\in\tangent_p\Ma$, we get one solution of the differential equation and then we can define the isomorphism \eqref{parallel}.

A {\it geodesic} of $\nabla$ is a curve with parallel tangent vector field,
\begin{equation}
\label{geodesic}
\nabla_t\dot{\gamma}\equiv 0\ ,
\end{equation}
which in local coordinates  reads as
\begin{equation}
\label{localgeod}
\ddot{\gamma}^k+\Gamma_{ij}^k\dot{\gamma}^i\dot{\gamma}^j=0\ .
\end{equation}
For all $p\in\Ma$ and $X_p\in \tangent_{p}\Ma$
there is a unique geodesic $\gamma_{X_p}$ such that,
\begin{equation}
\gamma_{X_p}(0)=p\text{, and }\dot{\gamma}_{X_p}(0)=X_p\text{.}
\label{initialcond}%
\end{equation}
Hence, by defining for $X_p\in \tangent_{p}\Ma$,
\begin{equation}
\exp_p(X_p):=\gamma_{X_p}(1)\text{,}%
\end{equation}
we obtain the {\it exponential map} at $p$. The exponential map is in general
well-defined at least in a neighborhood of zero in $\tangent_p\Ma$ and,
moreover, can be globally defined in special cases.

\subsection{Conjugate Symmetric Statistical Manifolds}

To the affine connection $\nabla$ we can associate two tensors, the {\it torsion} and the {\it curvature}. They are given by
\begin{eqnarray}
&&{\rm Tor}(X,Y)=\nabla_X Y-\nabla_Y X-[X,Y]\label{torsion}\\
\nonumber\\
&&\Riemann(X,Y)Z=\nabla_X\nabla_Y Z-\nabla_Y\nabla_X Z-\nabla_{[X,Y]}Z,\label{curvature}
\end{eqnarray}
where $X,Y,Z\in\Tau(\Ma)$ and $[X,Y]=XY-YX$ is the Lie bracket of $X$ and $Y$.

Analogously, we can associate two tensors to the dual connection: the torsion tensor and the curvature tensor of $\nabla^*$,
\begin{eqnarray}
&&{\rm Tor}^*(X,Y)=\nabla^*_X Y-\nabla^*_Y X-[X,Y]\label{torsion^*}\\
\nonumber\\
&&\Riemann^*(X,Y)Z=\nabla^*_X\nabla^*_Y Z-\nabla^*_Y\nabla^*_X Z-\nabla^*_{[X,Y]}Z\ . \label{curvature*}
\end{eqnarray}

Then ${\cal S}=(\Ma,\metric,\nabla,\nabla^*)$ is called a statistical manifold when both the connections $\nabla$ and $\nabla^*$ are {\it torsion free}, i.e ${\rm Tor}\equiv 0$ and ${\rm Tor}^*\equiv 0$. From this, it follows that the curvature $\Riemann$ satisfies the first  Bianchi identity,
\begin{eqnarray}
&& \Riemann(X,Y)Z+\Riemann(Y,Z)X+\Riemann(Z,X)Y=0\label{Bianchi1}
\end{eqnarray}
for all $X,Y,Z\in\Tau(\Ma)$. The same holds true for the curvature tensor $\Riemann^*$. 

Given the metric structure on $\Ma$, we can also consider the Riemann curvature tensor of $\nabla$ that is defined as follows
\begin{equation}
\label{RiemC}
\RC(X,Y,Z,W):=\metric\left(\Riemann_{XY}Z,W\right)\ .
\end{equation}
From Eq. \eqref{curvature} it immediately follows that $\RC(X,Y,Z,W)=-\RC(Y,X,Z,W)$ and in particular $\RC(X,X,Z,W)=0$. Moreover, from the first Bianchi identity \eqref{Bianchi1} we have that 
$$
\RC(X,Y,Z,W)+\RC(Y,Z,X,W)+\RC(Z,X,Y,W)=0\ .
$$
Analogously we can define the Riemann curvature tensor of $\nabla^*$,  
\begin{equation}
\label{RC*}
\RC^*(X,Y,Z,W):=\metric\left(\Riemann^*_{XY}Z,W\right)\ ,
\end{equation}
and same equalities  as $\RC$ hold true as well.

Consider now the Riemann curvature tensors $\RC$ and $\RC^*$ both together. We have the following result \cite{Lauritzen87},
\begin{pro}\label{}
If $\RC$ is the Riemann curvature tensor of $\nabla$ and $\RC^*$ the one of $\nabla^*$ we have that
\begin{equation}\label{alter}
\RC(X,Y,Z,W)=-\RC^*(X,Y,W,Z)\ .
\end{equation}
\end{pro}
\noindent{\bf Proof.} By considering $X,Y$ as part of an orthonormal frame on the tangent bundle $\tangent\Ma$ we can assume that $[X,Y]=0$. Owing to this consideration we have that
\begin{eqnarray*}
&&XY\metric\left(Z,W\right)= X\left(Y\metric(Z,W)\right)\\
&& = X\left(\metric\left(\nabla_Y Z,W\right)+\metric\left(Z,\nabla^*_YW\right)\right)\\
&&=\metric\left(\nabla_X\nabla_Y Z,W\right)+\metric\left(\nabla_Y Z,\nabla^*_X W\right)+\metric\left(\nabla_X Z,\nabla^*_YW\right)+\metric\left(Z,\nabla_X^*\nabla^*_Y W\right)
\end{eqnarray*}
By alternating $X$ and $Y$ we arrive at
\begin{eqnarray*}
&&0=[X,Y]\metric(Z,W)=XY\metric(Z,W)-YX\metric(Z,W)\\
&&=\RC(X,Y,Z,W)+\RC^*(X,Y,Z,W) \hspace{5cm} \square
\end{eqnarray*}
As direct consequence we can state the following result.
\begin{corollary}\label{coraltern}
The following conditions are equivalent,
\begin{enumerate}
\item $\RC\equiv\RC^*$
\item $\RC(X,Y,Z,W)=-\RC(X,Y,W,Z)$
\end{enumerate}
\end{corollary}
From the second condition in Cor. \ref{coraltern} we trivially have that
\begin{equation}
\label{Kobay1}
\metric\left(\Riemann(X,Y) Z,Z\right)\equiv 0\quad \mbox{for all}\ X,Y,Z\in\Tau(\Ma)\ .
\end{equation}
Another consequence of Cor. \ref{coraltern} is that $\nabla$ is flat if and only if $\nabla^*$ is flat. Let us now briefly discuss about the second condition of Cor. \ref{coraltern}, or equivalently the Eq. \eqref{Kobay1}, and see for which classes of statistical manifolds it holds true.

Given the dual structure $(\metric,\nabla,\nabla^*)$, we can obtain the Levi-Civita connection as follows \cite{Amari16},
\begin{equation}
\label{LC}
\overline{\nabla}:=\frac{1}{2}\left(\nabla+\nabla^*\right).
\end{equation}
In addition, we can define a totally symmetric cubic tensor $T$ \cite{Lauritzen87},
\begin{equation}
\label{cubic}
T(X,Y,Z):=\metric\left(\widetilde{T}(X,Y),Z\right),\quad \mbox{where}\ \widetilde{T}(X,Y):=\nabla_X Y-\nabla^*_X Y\ .
\end{equation}

Let us now define a $1$-parameter family of $\alpha$-connections on $\Ma$ as follows,
\begin{equation}
\label{alphaconnection}
\stackrel{\alpha}{\nabla}_X Y:=\overline{\nabla}_X Y-\frac{1}{2} \widetilde{T}(X,Y) .
\end{equation}
From the torsion-freeness of the statistical manifold $\esse$ and the symmetry of $\widetilde{T}$ we have that
\begin{equation}
\left(\stackrel{\alpha}{\nabla}\right)^*=\stackrel{-\alpha}{\nabla}\quad \mbox{and}\quad \stackrel{1}{\nabla}=\nabla,\ \stackrel{-1}{\nabla}=\nabla^*\ .
\end{equation}

We say that a statistical manifold $\esse=(\Ma,\metric,\nabla,\nabla^*)$ is {\it conjugate symmetric} if for all $\alpha$ the curvature tensor $\stackrel{\alpha}{\RC}$ fulfils the following relation,
\begin{equation}
\label{alphaconjsym}
\stackrel{\alpha}{\RC}\equiv\stackrel{-\alpha}{\RC}\ .
\end{equation}

Therefore, by means of Cor. \ref{coraltern} we have
\begin{pro}
\label{suffconjsym}
Sufficient conditions for a statistical $\esse=(\Ma,\metric,\nabla,\nabla^*)$ being conjugate symmetric are 
\begin{enumerate}
\item There exists $\alpha\neq 0$ such that $\stackrel{\alpha}{\RC}\  \equiv\ \stackrel{-\alpha}{\RC}$.
\item There exists $\alpha\neq 0$ such that $\stackrel{\alpha}{\RC}\equiv 0$, i.e. $\esse$ is $\alpha$-flat.
\end{enumerate}
\end{pro}
Finally, in a conjugate symmetric manifold, the Riemann curvature tensor satisfies all the identities as the Riemann curvature tensor of the Levi-Civita connection, i.e.
\begin{align}\label{RCidenty1}
&\stackrel{\alpha}{\RC}(X,Y,Z,W)=-\stackrel{\alpha}{\RC}(Y,X,Z,W);\\\label{RCidenty2}
&\stackrel{\alpha}{\RC}(X,Y,Z,W)+\stackrel{\alpha}{\RC}(Y,Z,X,W)+\stackrel{\alpha}{\RC}(Z,X,Y,W)=0;\\\label{RCidenty3}
&\stackrel{\alpha}{\RC}(X,Y,Z,W)=-\stackrel{\alpha}{\RC}(X,Y,W,Z);\\\label{RCidenty4}
&\stackrel{\alpha}{\RC}(X,Y,Z,W)=\stackrel{\alpha}{\RC}(Z,W,X,Y).
\end{align}

\vspace{.5cm}

The statistical manifold $\esse=(\Ma,\metric,\nabla,\nabla^*)$ is called {\it dually flat} if $\RC\equiv 0\equiv \RC^*$. Then, according to Eq. \eqref{alphaconjsym} and Pro. \ref{suffconjsym} we can say that a dually flat manifold is conjugate symmetric. In this particular case, there exists $\alpha_0$ such that $\stackrel{\alpha_0}{\RC}\equiv 0$ and then the statistical manifold $\esse$ is often referred as equivalent to {\it dually flat} manifold \cite{Amari16}. In this particular case, we can rely on two sets of local coordinates $\{\theta^i\}$ and $\{\eta_i\}$ such that
$$
\Gamma_{ijk}(\theta)=0,\quad \mbox{and}\ \stackrel{*}{\Gamma}_{ijk}(\eta)=0\ .
$$
Here, $\Gamma_{ijk}$ and $\stackrel{*}{\Gamma}_{ijk}$ are the connection symbols of $\nabla$ and $\nabla^*$, respectively. In local coordinates they are expressed by
\begin{equation}
\label{nablasymb}
\Gamma_{ijk}=\metric_{il}\Gamma_{jk}^l,\quad \stackrel{*}{\Gamma}_{ijk}=\metric_{il} {^*\Gamma}_{jk}^{*l}\ ,
\end{equation}
where $\Gamma_{jk}^l$ and ${\Gamma}_{jk}^{*l}$ are the Christoffel's symbols of $\nabla$ and $\nabla^*$, respectively. 
Additionally, if we consider the tangent vectors $\{\partial_i\}$ and $\{\partial^i\}$ of the local coordinates $\{\theta^i\}$ and $\{\eta_i\}$ we have that
\begin{equation}
\label{tangentdualflat}
\metric_{ij}\partial_i\partial^j=\delta_i^j,
\end{equation}
meaning that these tangent vectors are reciprocal orthogonal with respect to the metric tensor $\metric$.

\subsection{Parallel transport and curvature tensor}

Now we describe the connection between the parallel transport and the curvature tensor of the connection $\nabla$. Obviously, the same is for $\nabla^*$ connection. Roughly speaking, parallel transport along a loop $\Sigma$ based at $p\in\Ma$ provides the Lie group of rotations $\paralleltransport_{\Sigma}:\tangent_p\Ma\rightarrow\tangent_p\Ma$. This is called the {\it holonomy group} of $\nabla$ at $p$. Then, the Lie algebra of it is spanned by the curvature tensor of $\nabla$.

Given $p\in \Ma$, let
\begin{equation}
\label{loop}
\Loop_p:=\left\{\Sigma:[0,1]\rightarrow \Ma\ |\ \Sigma(0)=\Sigma(1)=p\right\}
\end{equation} 
be the set of piecewise smooth loop based on $p$ and assume that $\Ma$ is simply connected. Then, each $\Sigma\in\Loop_p$  is homotopic to the trivial loop.

Therefore the {\it holonomy} of $\nabla$ at $p\in \Ma$ is defined as the subset of $\Aut(\tangent_p\Ma)$, i.e. the automorphisms of $\tangent_p \Ma$,
\begin{equation}
\label{holonomy}
\Hol_p:=\left\{\paralleltransport_\Sigma\in \Aut(\tangent_p\Ma)\ |\ \Sigma\in\Loop_p \right\}\ .
\end{equation}
Basics properties of $\Hol_p$ are listed in the following proposition.
\begin{pro}
The following basic properties of $\Hol_p$ hold true:
\begin{enumerate}
\item $\Hol_p$ is a closed Lie subgroup of $\Aut(\tangent_p\Ma)$ and its Lie algebra $\mathfrak{hol}_p \subset \End(\tangent_p\Ma)$ is called the {\it holonomy algebra} at $p$.
\item Given $\Sigma^{\prime}:[0,1]\rightarrow \Ma$ such that $\Sigma^{\prime}(0)=p$ and $\Sigma^{\prime}(1)=q$. Let $P_{\Sigma^{\prime}}:\tangent_p\Ma\rightarrow \tangent_q \Ma$ the parallel transport along $\Sigma^{\prime}$. Then
$$
P_{\Sigma^{\prime}}\circ \Hol_p\circ {P_{\Sigma^{\prime}}}^{-1}=\Hol_q.
$$
\end{enumerate}
\end{pro}
From the second property in the latter Proposition,  it follows that the holonomy groups are independent of the base point. 

Since $\nabla$ is torsion free, the {\it Ambrose-Singer Holonomy Theorem} \cite{Schwach99} supplies a very remarkable connection between the curvature tensor $\Riemann$ and the holonomy algebra $\hol_p(\nabla)$ of $\nabla$. It states that $\hol_p(\nabla)$ is generated by operators $\Riemann_{\Sigma}(x,y):=\paralleltransport_\alpha\circ \Riemann({\paralleltransport_{\Sigma}}^{-1}x,{\paralleltransport_{\Sigma}}^{-1}y)\circ {\paralleltransport_{\Sigma}}^{-1}$,
\begin{equation}
\label{ASTh}
\hol_p=\langle\left\{( \Riemann)_{\Sigma}(x,y)\ |\ x,y\in \tangent_p\Ma, \Sigma\ \mbox{a loop at}\ p \right\}\rangle.
\end{equation}
Eq. \eqref{ASTh} shows that $\hol_p(\nabla)$ is the vector subspace of $\End(\tangent_p\Ma)$
spanned by the endomorphisms $\Riemann_{\Sigma}(x,y)$. Thus, $\Riemann$ determines $\hol_p(\nabla)$ and, hence $\Hol(\nabla)$.
Therefore, if we consider the case of a flat manifold where $\Riemann\equiv 0$ we have$\hol_p(\nabla)=0$, from
which follows that $\Hol(\nabla) = Id$.

\vspace{.5 cm}

For the purpose of the present manuscript, the previous theoretical setting, that highlights the connection between holonomy and curvature tensor, is performed into the following result.

\begin{lemma}\label{lemma}
Let $B$ be a smooth closed $2$-disk such that $p\in\partial B$ and $B$ is foliated by connecting $\nabla$-geodesics segment starting from $p$. Then
\begin{equation}
\paralleltransport_{\Sigma}Z_{p}-Z_{p}=\int_B \frac{\paralleltransport\left(\Riemann(X,Y)Z\right)}{\|X\wedge Y\|} dA,
\end{equation}
where
\begin{itemize}
\item $dA$ is the surface area measure on $B$ induced by the Riemannian metric tensor $\metric$ on $\Ma$.
\item $X$ and $Y$ are linearly independent vector fields on $B$.
\item $\Sigma:[0,1]\rightarrow \partial B$  is a parametrization of $\partial B$ such that $\Sigma(0)=\Sigma(1)=p$ and, given
any inward pointing vector $X\in\tangent_{p}B$, the orientation of $(\dot{\Sigma},X)$ is the same as
$(X,Y)$.
\item  $Z_{p}\in\tangent_{p}\Ma$ and $Z$ is defined by parallel translating $Z_{p}$ first along the parametrized
curve $\Sigma$ and then, for each $0\leq s\leq 1$, along the unique $\nabla$-geodesic segment going from
$\Sigma(s)\in\partial B$ to $B$.
\item $\paralleltransport$ is parallel translation from each point in $B$ to $p$ along the unique $\nabla$-geodesic segment joining them.
\end{itemize}
\end{lemma}
\noindent {\bf Proof.} The proof of this result is provided in \cite{Yang}. However, for its relevance to our work, we report it here. 

Consider a map $H:[0,1]\times [0,1]\rightarrow B$ such that $H(1,t)=\Sigma(t)$ and $H(\cdot,t)$ is the $\nabla$-geodesic connecting $p$ to $\Sigma(t)$, for all $t\in[0,1]$. Let us denote
$$
S(s,t)=\frac{\partial H}{\partial s},\quad T(s,t)=\frac{\partial H}{\partial t}
$$
then we have $[S,T]=0$ since $\nabla$ is torsion free. Let us observe that $T$ is the Jacobi vector field along each geodesic $H(\cdot,t)$. Define now $$J=T-\frac{\metric\left(S,T\right)}{\|S\|^2}S.$$
Clearly $J$ is orthogonal to $T$. Then we have that
$$
dA=\|T\wedge S\|ds dt=\| S\| \|T\| ds dt.
$$

Let $\{e_i\}\subset\tangent_p\Ma$ be an orthonormal frame and extend it by parallel transport along each $\nabla$ geodesic $H(\cdot,t)$. In particular we have that $$ \nabla_t e_i(0,t)=0, \quad \nabla_s e_i(s,t)=0,$$
for all $(s,t)\in [0,1]\times [0,1]$. In addition, we also have that $$\nabla_t Z(1,t)=0, \quad \nabla_s Z(s,t)=0,$$
for all $(s,t)\in [0,1]\times [0,1]$. Let us now note that $Z_p=Z(0,0)=Z(1,0)$ and that $\paralleltransport_{\Sigma} Z_p=Z(1,1)=Z(0,1)$. Then we have that
\begin{eqnarray*}
\langle e_i(p),(\paralleltransport_{\Sigma}Z_p-Z_p)\rangle_p &=& \langle e_i(0,1),Z(0,1)\rangle_p-\langle e_i(0,0)
,Z(0,0)\rangle_p\\
&=&\int_0^1 \partial_t\langle e_i(0,t), Z(0,t)\rangle dt\\
&=& \int_0^1 \langle e_i, \nabla_t Z(0,t)\rangle dt\\
&=& \int_0^1\left[\langle e_i,\nabla_t Z(1,t)\rangle-\int_0^1\partial_s\langle e_i, \nabla_t Z(s,t)\rangle ds\right]dt\\
&=&-\int_0^1 \int_0^1\langle e_i, \nabla_s\nabla_t Z(s,t)\rangle ds dt\\
&=&-\int_0^1 \int_0^1\langle e_i, \Riemann(S,T) Z(s,t)\rangle ds dt\\
&=&-\int_0^1 \int_0^1\langle e_i, \Riemann(S,J) Z(s,t)\rangle ds dt\\
&=&-\int_0^1 \int_0^1\langle e_i, \Riemann(\sigma,\tau) Z(s,t)\rangle \|S\| \|T\| ds dt\\
&=& \int_0^1 \frac{\langle e_i, \Riemann(X,Y) Z\rangle}{\|X\wedge Y\|} dA,
\end{eqnarray*}
where $\sigma=S\slash\|S\|$, $\tau=T\slash\|T\|$ form an orthonormal frame on $B$. Finally the result follows from
$$
\paralleltransport_{\Sigma} \left(\Riemann (X,Y)Z\right)=\sum_i e_i(0,0)\langle e_i(s,t),\Riemann(X,Y)Z\rangle \ .$$

\end{document}